\documentclass[11pt]{article}

\usepackage[utf8]{inputenc}
\usepackage[T1]{fontenc}
\usepackage[english]{babel}
\usepackage{amsmath}
\usepackage{amssymb}
\usepackage{amsfonts}
\usepackage{graphicx}
\usepackage{booktabs}
\usepackage{float}
\usepackage{caption}
\usepackage{subcaption}
\usepackage[margin=1in]{geometry}
\usepackage{url}
\usepackage{hyperref}
\usepackage{authblk}
\usepackage{setspace}
\usepackage{indentfirst}
\usepackage{longtable}
\usepackage{array}
\usepackage{multirow}
\usepackage{enumitem}
\usepackage{booktabs}
\usepackage{tabularx}
\usepackage{xcolor}
\usepackage{orcidlink}
\usepackage{lineno}
\usepackage[section]{placeins}

\newcommand{\rev}[1]{{\color{black}#1}}
\title{A Trajectory-Based Mathematical Framework for Wildfire Front Propagation in Heterogeneous Landscapes}

\author[]{Hengameh R. Dehkordi\,\orcidlink{0000-0002-1738-3373}}

\affil[1]{Center of Mathematics, Computing, and Cognition (CMCC),
	Federal University of ABC (UFABC), 09210-580, Santo Andr\'e, SP, Brazil\\
}

\affil[]{hengameh.r@ufabc.edu.br}
\date{\today}

\begin{document}
	
	\maketitle
	
\begin{abstract}
We present a mathematical formulation based on trajectories for modeling the propagation of wildfire fronts in heterogeneous landscapes and apply it to the Brazilian Pantanal. The formulation is made up of a direction-dependent propagation law and Finsler orthogonality along a moving boundary. This results in a nonlinear local system for outward fire front trajectories.
We provide a strictly convex translated ellipse representation at each point of the evolving boundary which guaranties uniqueness and existence of the chosen outward velocity.
The effective propagation coefficients that are dependent on space and the hourly wind forcing permit the perimeter to evolve irregularly without imposing a fixed global front geometry. Computed trajectories were georeferenced and isotropically spatially registered with VIIRS active-fire observations. ERA5 wind and ERA5-Land vegetation data were used to describe the environmental forcing and landscape heterogeneity. The full-run Dice coefficients for the observation-guided spatial comparisons for the exploratory experiments for September 3-5, 2020 are 0.578, 0.571, and 0.467 respectively. One-at-a-time sensitivity perturbations of the spatial propagation coefficient, wind speed magnitude and rate-law exponents quantify the governing influence of the environmental forcing and parameter variations on the spatial extent of the calculated fire advance. These observed patterns also guided calibration and, in the absence of exact timestamps of satellite acquisition, the results show within-event spatial agreement rather than independent predictive validation. The framework provides a mathematically rigorous basis for subsequent uncertainty quantification, independent validation over multiple events and prospective sector-wise assessments of simulated fire spread.
\end{abstract}


	
	\section{Introduction}
	
	The Pantanal is one of the world's largest tropical wetland systems. It is
	characterized by extensive floodplains, generally low terrain gradients, and
	pronounced seasonal variations in water availability. In Brazil, it extends
	across the states of Mato Grosso and Mato Grosso do Sul within the Upper
	Paraguay River Basin. Its hydrological dynamics strongly influence the spatial
	and temporal distribution of vegetation
	\cite{pott_pantanal_2011,libonati2020rescue}. The Pantanal is a transitional
	region influenced by the Cerrado, Chaco, Atlantic Forest, and Amazon
	Rainforest. Consequently, it contains a heterogeneous mosaic of grasslands,
	shrublands, forests, and wetland vegetation
	\cite{pott_pantanal_2011,fuller2014recent,silva2021caracteristicas}.
	
	This environmental heterogeneity has important implications for wildfire
	behavior. Seasonal flooding favors aquatic and low-growing vegetation. During
	the dry season, however, reductions in water availability can increase the
	susceptibility of vegetation and accumulated fuels to fire. Spatial
	differences in vegetation structure, fuel conditions, and moisture may
	therefore produce substantial variations in fire propagation, even within
	relatively small areas \cite{pelissari2023dynamics}. Fire spread is also
	strongly influenced by wind and fuel characteristics
	\cite{andrews2018rothermel}. Accordingly, the present study focuses on spatial
	vegetation heterogeneity, time-dependent wind, and prescribed fuel moisture.
	Topographic effects are not included in the proposed model. The Pantanal is
	characterized by generally low-relief terrain, and slope effects are assumed
	to be negligible for the present case study.
	
	A large body of wildfire research has addressed the estimation of fire rate
	of spread through physical and empirical models. Among the most influential
	contributions is the surface-fire model developed by Rothermel
	\cite{rothermel1972mathematical}. This model relates fire spread to fuel-bed
	properties, fuel moisture, wind, and terrain conditions. It has provided an
	important foundation for subsequent fire-behavior modeling and has been
	extensively examined and applied \cite{andrews2018rothermel}. Empirical
	studies have likewise demonstrated the important roles of wind, fuel
	characteristics, and moisture in surface-fire propagation
	\cite{xavier1998forest,fuller2014recent}. Nevertheless, pointwise
	application of a complete physical or semi-empirical fire-spread formulation
	may require numerous spatially distributed fuel parameters. Such information
	is not always available at the spatial resolution required for local
	propagation calculations. In the present study, the Rothermel framework is
	therefore used as a physical reference for the dependence of fire spread on
	environmental conditions. It is not used as the direct operational
	rate-of-spread equation.
	
	Advanced mathematical approaches have also been developed to represent
	direction-dependent wildfire propagation. Markvorsen \cite{Markvorsen2016}
	introduced a geodesic-spray paradigm for wildfire spread.
	Geometric formulations, including Randers-type structures, were investigated
	for non-isotropic wildfire propagation in
	\cite{dehkordi2019,dehkordi2022applications,dehkordi2026satellite}. Further developments have
	addressed the mathematical construction and evolution of wildfire fronts
	\cite{dehkordi2023wildfire}. Related formulations have also been extended to
	incorporate environmental directional effects \cite{angel2023general}.
	Collectively, these studies demonstrate the potential of advanced
	mathematical techniques for representing anisotropic and spatially
	differentiated fire propagation.
	
	Operational fire-growth simulators provide another important class of
	approaches. FARSITE, for example, propagates fire fronts using Huygens'
	principle and elliptical growth rules. Its local rates of spread are derived
	from the Rothermel formulation
	\cite{richards1990elliptical,finney1998farsite,andrews2018rothermel}.
	Such systems provide mature frameworks for landscape-scale fire-growth
	simulation. Their application, however, generally relies on detailed spatial
	information describing fuels and environmental conditions. Obtaining such
	information at sufficient spatial resolution and consistency may be
	challenging in heterogeneous or data-limited environments.
	
	The approach proposed in this study follows a different mathematical
	construction. The perimeter is discretized into points rather than evolved
	through a prescribed geometric fire-growth rule. At each point, a local
	nonlinear system is solved to determine an admissible outward propagation
	trajectory. The endpoints of these trajectories define the subsequent
	perimeter. Local propagation rates are described by a compact, physically
	interpretable relation. This relation incorporates wind, prescribed fuel
	moisture, and a spatially varying effective propagation coefficient. The
	resulting propagation response can vary along the perimeter without imposing
	a uniform rate of spread or a fixed global propagation shape. The formulation
	also provides a mechanism for incorporating observational information into
	the local propagation field when detailed spatial fuel-bed parameters are
	unavailable.
	
	In parallel with advances in fire-spread modeling, satellite remote sensing
	has substantially improved the observation of wildfire activity. MODIS
	active-fire products have been widely used for global fire detection
	\cite{justice2002modis}. The VIIRS (Visible Infrared Imaging Radiometer Suite)
	375~m active-fire product provides finer-resolution observations and improves
	the detection and spatial characterization of fire activity
	\cite{schroeder2014viirs}. Reanalysis products provide complementary
	environmental information. ERA5 supplies temporally resolved atmospheric
	variables \cite{poblet2020era5}, whereas ERA5-Land provides land-surface
	information, including low- and high-vegetation cover and vegetation type.
	Satellite observations and environmental data have also been used to
	investigate wildfire dynamics and changing fire conditions in the Pantanal
	\cite{pelissari2023dynamics,gomes2022projections}.
	
	Despite these developments, an important modeling challenge remains. Observed
	fire evolution and spatially heterogeneous environmental information must be
	linked to a mathematical rule capable of propagating an irregular fire
	perimeter. Detailed fire-spread formulations may require spatial fuel
	information that is unavailable or insufficiently resolved. Satellite
	active-fire detections, in contrast, provide valuable information about the
	location and evolution of fire activity. However, they do not directly
	specify a continuous local propagation law. Mathematical perimeter
	propagation requires local information describing both the magnitude and
	direction of fire-front advancement. Bridging these sources of information
	therefore requires a framework that connects observational and environmental
	data with a well-defined mathematical propagation procedure.
	
	Motivated by this gap, the present study develops a satellite-assisted
	trajectory framework for modeling spatially heterogeneous wildfire
	propagation in the Brazilian Pantanal. VIIRS active-fire observations are
	used to characterize the observed fire evolution and construct the initial
	fire boundaries. ERA5 provides hourly wind magnitude and direction.
	ERA5-Land provides information on spatial vegetation heterogeneity. These
	data are coupled with a spatially varying effective propagation coefficient,
	$\gamma(x,y)$. Fire-front evolution is then described through a
	trajectory-based mathematical formulation. A nonlinear system is solved
	locally at each discretized perimeter point. Repeated application of this
	procedure generates the evolving fire perimeter. The local propagation
	behavior can vary in space, while the wind forcing can vary in time. The
	framework therefore provides a mathematically explicit connection among
	observed fire development, environmental heterogeneity, and local fire-front
	propagation. At the same time, it does not require a complete pointwise
	parameterization of the fuel bed.
	
	The objective of this first-stage study is to establish the mathematical and
	computational feasibility of the proposed integration and to assess its ability
	to represent the principal spatial characteristics of an observed wildfire
	event. {For each daily experiment, the available VIIRS pattern informs estimation of the spatial coefficient $\gamma(x,y)$ through propagation calculated over the first 12 model steps. The field is then held fixed during steps 13--24. Since the retained satellite archive contains dated snapshots rather than verified observations at the 12-h separation, this division is computational and does not constitute independent temporal validation. In particular, the same observed pattern can inform calibration and the subsequent spatial comparison. The experiments are therefore presented as exploratory, observation-guided demonstrations of the formulation.}
	
	{The choice of one well-characterized event serves the methodological objective of examining the trajectory construction under spatially heterogeneous conditions while keeping the case-specific environmental assumptions consistent. Independent-event validation and systematic comparisons with alternative propagation models would require harmonized initial perimeters, observation times, spatial fuel information, calibration protocols, and reference datasets. Such comparisons are necessary to establish generalizability, but they are not claimed by this first-stage study.}
	
	The main contributions of this study are {fivefold}:
	\begin{enumerate}
		\item A trajectory-based mathematical framework is developed for
		propagating irregular wildfire perimeters. A local nonlinear system is
		solved at each discretized boundary point to determine its outward
		propagation trajectory, allowing the fire front to evolve without
		imposing a uniform rate of spread or a fixed global propagation shape.
		
		\item Spatially heterogeneous and temporally varying fire conditions are
		incorporated by combining a local effective propagation coefficient
		$\gamma(x,y)$ with time-dependent wind forcing. VIIRS observations,
		ERA5 wind data, and ERA5-Land vegetation information are integrated to
		connect observed fire development and environmental heterogeneity with
		the mathematical propagation model.
		
		\item An observation-guided calibration and propagation protocol is
		introduced for $\gamma(x,y)$. For each daily experiment, the available
		satellite-observed fire pattern guides the estimation of the coefficient
		field over the first 12 model hours. The resulting field is then held fixed
		during hours 13--24, providing a {fixed-parameter continuation rather than independent temporal validation} without further adjustment of $\gamma(x,y)$.
		
		\item The calculated trajectories are rasterized at the VIIRS spatial
		resolution for direct spatial comparison with satellite observations.

		\item {The framework provides a basis for prospective sector-wise indicators of modeled fire spread: the fraction of simulated advance assigned to each initial-front sector could be compared across hypothetical risk-management scenarios. This capability is proposed as an analytical use of the computed trajectories and is not operationally validated here.}
	\end{enumerate}
	
	The remainder of this paper is organized as follows. Section~2 describes the
	study area and the observational and environmental data used in the Pantanal
	case study. Section~3 presents the proposed methodology. It includes the
	rate-of-spread formulation and the trajectory-based mathematical model for
	fire-front propagation. Section~4 presents the numerical results and model
	evaluation. Section~5 discusses the principal findings, assumptions, and
	limitations of the proposed framework. Finally, Section~6 summarizes the
	main conclusions and outlines directions for future research.

	\section{Data and Study Area}
	
	This section describes the study area and the data sources used to characterize the wildfire event and the environmental conditions relevant to its propagation. Satellite observations and reanalysis products provide complementary information on fire activity, wind conditions, and spatial vegetation heterogeneity. These data constitute the observational and environmental basis for the wildfire propagation model developed in the subsequent section.

	\subsection{Study Area: The Brazilian Pantanal}
	
\begin{figure}[htbp]
		\centering
		\begin{subfigure}[t]{0.4\textwidth}
			\centering
			\includegraphics[width=\linewidth]{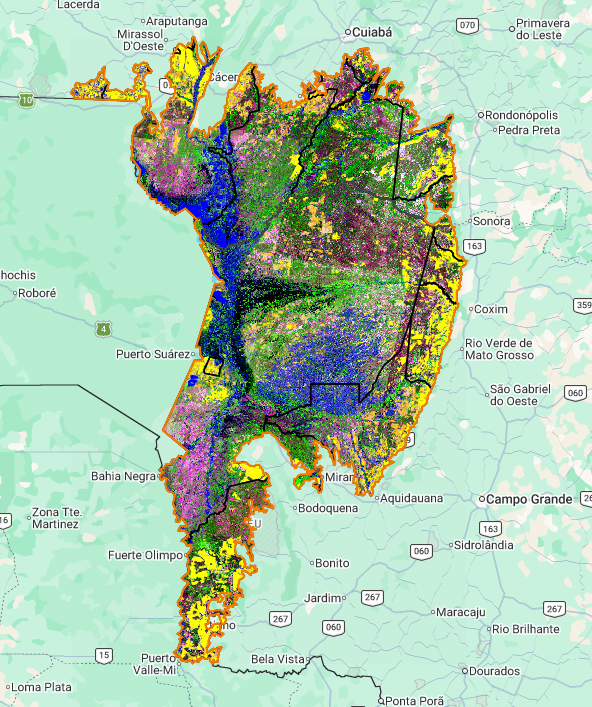}
			\caption{Brazilian Pantanal map}
		\end{subfigure}
		\hfill
		\begin{subfigure}[t]{0.5\textwidth}
			\centering
			\includegraphics[width=\linewidth]{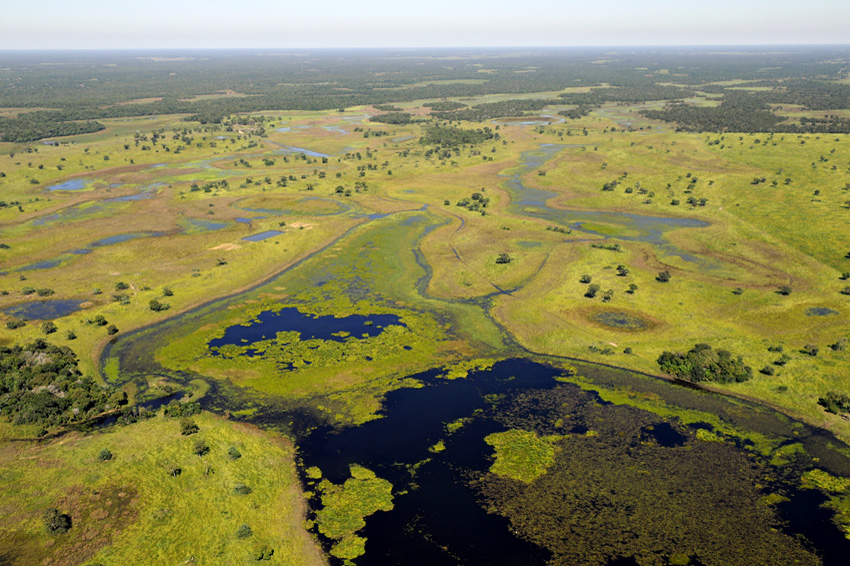}
			\caption{ Pantanal is the world's largest tropical wetland area}
		\end{subfigure}
		\caption{ Sources: TerraBrasilis \href{https://terrabrasilis.dpi.inpe.br/app/map/deforestation?hl=en}{[1]}; \rev{WWF image source:} \href{https://www.wwf.org.br/?32263/Pantanal-wetlands-Images}{\rev{website}}. }
		\label{pan-map}
	\end{figure}

	The study area is located in the Brazilian Pantanal, within the states of Mato Grosso and Mato Grosso do Sul (Fig.~\ref{pan-map}). The Pantanal is characterized by extensive floodplains, low topographic gradients, and strong seasonal hydrological dynamics. Its transitional location among the Cerrado, Chaco, Atlantic Forest, and Amazon Rainforest biomes contributes to a complex vegetation mosaic consisting of grasslands, shrublands, forests, and wetland vegetation \cite{pelissari2023dynamics}. Such spatial heterogeneity in vegetation structure may produce substantial variations in fuel characteristics across the study area, resulting in corresponding spatial differences in wildfire propagation behavior.
	
	The case study considers a wildfire event in the Brazilian Pantanal for which successive satellite observations are used to characterize the initial fire perimeter and its subsequent evolution. The analysis focuses primarily on the effects of spatial vegetation heterogeneity and time-dependent wind conditions on fire propagation, while fuel moisture is incorporated as a prescribed model parameter.

	\subsection{Satellite Fire Detection Data}
	
	Active-fire information for the case study is obtained from the Visible Infrared Imaging Radiometer Suite (VIIRS) 
	375~m active-fire product 
	distributed through the NASA Fire Information for Resource Management System 
	(FIRMS). The relatively fine spatial resolution of VIIRS makes these 
	observations suitable for characterizing the spatial evolution of the 
	wildfire considered in this study \cite{schroeder2014viirs}. The VIIRS 
	detections are used to identify the observed fire activity at successive 
	times and to support the characterization of the initial fire perimeter and 
	the subsequent evolution of the fire front.
	
	The satellite observations provide discrete snapshots of active-fire
	locations rather than a continuous record of fire-front evolution. They are
	used within the estimate--freeze--propagate procedure adopted in this study.
	For each daily experiment, the available satellite-observed fire pattern is
	used to guide the estimation of the event-specific effective propagation
	coefficient $\gamma(x,y)$ by evaluating the calculated propagation over the
	first 12 model hours. Thus, the 12-h interval defines the model calibration
	window and should not be interpreted as a period of continuous satellite
	observation. After this calibration stage, $\gamma(x,y)$ is held fixed and
	no further adjustment is made during hours 13--24. {Because the available observed pattern also informs calibration, the latter interval represents a fixed-parameter model continuation and cannot be interpreted as independently observed or out-of-sample validation.}

	\subsection{Wind Data from ERA5 Reanalysis}
	
	Hourly near-surface wind information is obtained from the ERA5 reanalysis produced by the European Centre for Medium-Range Weather Forecasts (ECMWF) \cite{poblet2020era5}. ERA5 provides eastward and northward wind components at 10~m above the surface, denoted by $u_{10}$ and $v_{10}$, respectively. The corresponding wind-speed magnitude is computed as
	\begin{equation}
		U(x,y,t)
		=
		\sqrt{u_{10}^{2}(x,y,t)+v_{10}^{2}(x,y,t)}.
		\label{wind_speed}
	\end{equation}
	
	In the present implementation, the ERA5 wind field is used directly as the environmental wind forcing. No fuel-dependent wind adjustment factor is introduced to reinterpret the reanalysis wind as a midflame wind speed. This choice is consistent with the phenomenological nature of the local propagation model developed below and avoids introducing an additional uncertain conversion that would itself depend on vegetation and fuel structure.
	
	When wind direction is expressed using the conventional meteorological convention, $\theta_{\mathrm{meto}}$ denotes the direction \emph{from which} the wind originates, measured clockwise from North. The propagation model instead requires the  mathematical convention  $\theta_{\mathrm{math}}$, that is the direction \emph{toward which} the wind vector acts, measured counterclockwise from East. The required direction is therefore
	\begin{equation}
		\theta_{\mathrm{math}}
		=
		\operatorname{mod}
		\left(
		270^\circ-\theta_{\mathrm{meto}},
		360^\circ
		\right).
		\label{eq:wind_direction_conversion}
	\end{equation}
	
	This transformation simultaneously changes the reference axis, reverses the meteorological ``from'' direction to the propagation ``toward'' direction, and converts the clockwise meteorological convention to the counterclockwise mathematical convention. The modulo operation ensures that the resulting angle lies in $[0^\circ,360^\circ)$.
	
	\subsection{Vegetation and Fuel Data}
	\label{subsec:vegetation_data}
	
	The vegetation of the Pantanal is highly heterogeneous and cannot be represented adequately by a single fuel type. Earlier ecological and fire-behavior descriptions provide useful qualitative analogies for the fine and wetland fuels encountered in the region. For example, grass-dominated Pantanal fuels share some structural characteristics with the Western-grass categories used in classical fuel-model descriptions, while mixed grass--woody areas can be compared qualitatively with savanna-type fuel complexes. Wetland vegetation also motivates comparison with Sawgrass-like environments \cite{shantz_great_plains_1923,keyser_savannah_woodland,fuller2014recent}. These analogies are used only as qualitative references; they are not taken to imply that Pantanal vegetation is identical to the corresponding North American fuel models.
	
	For the operational case study, vegetation heterogeneity is characterized directly from ERA5-Land fields. Six variables are considered: low-vegetation cover, high-vegetation cover, type of low vegetation, type of high vegetation, leaf area index of low vegetation, and leaf area index of high vegetation. The low- and high-vegetation cover variables provide fractional information within each grid cell, the corresponding vegetation-type variables identify the dominant class assigned to each vegetation stratum, and the leaf area index (LAI) variables provide an additional quantitative description of vegetation density.
	
	Analysis of the ERA5-Land vegetation data confirms that the study domain cannot be represented by a single homogeneous fuel or vegetation class. The low-vegetation component is predominantly associated with short-grass and tall-grass classes, whereas the high-vegetation component is mainly represented by evergreen broadleaf vegetation and interrupted forest. Moreover, the fractional contributions of low and high vegetation vary substantially across the study domain. For example, in a grid cell located close to the wildfire region at approximately $(-17.25^\circ,-56.50^\circ)$, high vegetation accounts for about $74.5\%$ of the grid cell and is classified as evergreen broadleaf vegetation, whereas low vegetation accounts for about $24.7\%$ and is classified as short grass. In a nearby grid cell at approximately $(-17.25^\circ,-56.25^\circ)$, the distribution changes to approximately $47.2\%$ high vegetation and $52.8\%$ low vegetation, represented by interrupted forest and tall grass, respectively (Fig.\ref{fig:vegetation_spatial}).
	
	The LAI fields provide further evidence of spatial differences in vegetation structure. Across the analyzed grid cells, high-vegetation LAI varies from approximately $2.4$ to $5.1~\mathrm{m^2\,m^{-2}}$, while low-vegetation LAI is approximately $1.7$--$2.2~\mathrm{m^2\,m^{-2}}$. Thus, neighboring locations may differ not only in their dominant vegetation class and fractional cover, but also in the density of the corresponding vegetation strata. Taken together, vegetation type, fractional cover, and LAI indicate pronounced spatial heterogeneity in the environmental conditions surrounding the wildfire.
	
	These data demonstrate that even within the relatively limited area surrounding the wildfire, the vegetation structure changes considerably from one location to another. Consequently, assigning a single uniform fuel class or a constant baseline propagation coefficient to the complete fire perimeter would not adequately represent the observed environmental heterogeneity. This data analysis provides the empirical motivation for the spatially varying propagation coefficient introduced in Section~\ref{subsec:spread_rates}.
	
	It is important to distinguish this environmental characterization from the numerical estimation of $\gamma(x,y)$. The ERA5-Land vegetation variables are not converted directly into values of $\gamma$ and are not used as a deterministic fuel-to-spread lookup table. Rather, they provide independent evidence that the baseline propagation behavior should be allowed to vary spatially. The numerical values of $\gamma(x,y)$ are instead estimated through an
	observation-guided calibration in which the available satellite fire pattern
	is compared with the calculated propagation over the first 12 model hours,
	as described in Sections~3.1 and 4.2.
	
	\subsection{Data Provenance and Reproducibility}
	\label{subsec:data_provenance}
	
	The fire observations are from the VIIRS 375-m active-fire product distributed
	through NASA's Fire Information for Resource Management System (FIRMS),
	available at \url{https://firms.modaps.eosdis.nasa.gov/}. The available records identify the three observations by date,
	September~3--5, 2020, but the exact satellite overpass timestamps and the
	original FIRMS query identifiers are not available in the retained archive.
	Accordingly, the observations are treated in this study as dated satellite
	snapshots rather than as hourly measurements. No hourly VIIRS acquisition
	sequence is assumed in the calibration or evaluation procedures. {The exact alignment between the 24 model steps and true satellite acquisition times is therefore uncertain; the starting time associated with the ERA5 sequence cannot be independently verified from the retained VIIRS records. A temporal-offset sensitivity analysis would be needed to quantify this uncertainty.}
	
	Hourly wind speed and direction were obtained from ERA5 through the
	Copernicus Climate Data Store (\url{https://cds.climate.copernicus.eu/}).
	The $10$-m eastward and northward wind components were converted to speed
	and meteorological direction before the transformation in
	Eq.~\eqref{eq:wind_direction_conversion}. ERA5 is hourly with a nominal
	$0.25^\circ$ grid for the fields used here. Vegetation information was
	obtained from ERA5-Land through the same service. Low- and high-vegetation
	cover and vegetation-type variables were used to describe heterogeneity at
	the product's nominal approximately $0.1^\circ$ grid.
	
	Preprocessing consisted of extracting fire-perimeter vertices, retaining
	decimal-degree coordinates, converting them to a local Cartesian system in
	kilometres, grouping repeated wind, direction, and moisture coordinates into
	chronological hourly records, and interpolating environmental values with
	inverse-distance weighting of exponent two. Wind direction was interpolated
	through sine and cosine components to avoid angular discontinuities. Moisture
	was prescribed as $M=10\%$ in all three case-study workbooks.

	\section{Methodology}
	
	The methodology implements the flexible data-assisted propagation
	framework: it couples satellite observations, environmental data, and
	established fire-spread principles within a trajectory-based formulation
	whose inputs are the available, empirical, observational and environmental data
	streams. VIIRS active-fire observations define the initial perimeter,
	provide the information from which the local spread rates are calibrated,
	and supply the comparison data; hourly ERA5 wind provides the atmospheric
	forcing; and ERA5-Land vegetation characterizes the spatial heterogeneity
	that motivates a nonuniform coefficient field. At each perimeter point, a
	coupled system of nonlinear equations is solved for the outward
	propagation trajectory, and the endpoints of these trajectories define
	the perimeter at the next hourly step; recursive application generates
	the evolving fire front.

	\subsection{Spatially Varying Rates of Fire Spread}
	\label{subsec:spread_rates}
	
	Laboratory and field experiments consistently show that wildfire spread
	increases with wind intensity and decreases with fuel moisture, with
	magnitudes that depend strongly on vegetation and fuel structure
	\cite{andrews2018rothermel,xavier1998forest,fuller2014recent}. 
	Motivated by these observations, the local head-fire rate of spread at
	location $(x,y)$ and time $t$ is represented as
	\begin{equation}
		R_H(x,y,t)
		=
		\gamma(x,y)
		\left[1+U^*(x,y,t)\right]^{\alpha}
		\exp\left[-\beta M(x,y,t)\right],
		\label{eq:RH}
	\end{equation}
	where $\gamma(x,y)$ is a spatially varying baseline propagation
	coefficient, $U^*(x,y,t)=U(x,y,t)/U_{\mathrm{ref}}$ is the dimensionless
	normalized wind speed with $U_{\mathrm{ref}}=1~\mathrm{m\,s^{-1}}$, and
	$M(x,y,t)$ denotes the local fuel-moisture condition. The exponent
	$\alpha$ controls the sensitivity of the head-fire rate to wind speed,
	whereas $\beta$ determines the moisture attenuation; since the wind and
	moisture terms are dimensionless, $\gamma(x,y)$ has the same units as
	$R_H$. \rev{In the numerical experiments, $M=10$ denotes moisture expressed in percentage points (10\%), rather than the fractional value $0.10$; accordingly, $\beta$ is expressed per percentage point.}
	
	The values of $\alpha$ reported in empirical fire-spread studies vary
	according to fuel type. Representative intervals reported in the literature
	are approximately
	\begin{equation}
		1.0 \leq \alpha \leq 1.5
		\qquad \text{for grass fuels},
		\label{eq:alpha_grass}
	\end{equation}
	\begin{equation}
		0.8 \leq \alpha \leq 1.2
		\qquad \text{for shrub fuels},
		\label{eq:alpha_shrub}
	\end{equation}
	and
	\begin{equation}
		0.3 \leq \alpha \leq 0.8
		\qquad \text{for forest-litter fuels}
		\label{eq:alpha_forest}
	\end{equation}
	\cite{andrews2018rothermel}. These intervals show that the wind sensitivity
	of the spread rate differs considerably among fuel environments, so an
	appropriate value of $\alpha$ is selected within the literature-supported
	intervals according to the vegetation characteristics of the experiment.
	
	The moisture-sensitivity coefficient $\beta$ is likewise fuel-dependent:
	values of approximately $0.05 \leq \beta \leq 0.07$ are reported for
	\textit{Pinus pinaster} needle fuels, and representative intervals extend to
	\begin{equation}
		0.05 \leq \beta \leq 0.10
		\qquad \text{for relatively dry fuels},
		\label{eq:beta_dry}
	\end{equation}
	\begin{equation}
		0.10 \leq \beta \leq 0.15
		\qquad \text{for moderately moist fuels},
		\label{eq:beta_moderate}
	\end{equation}
	and
	\begin{equation}
		0.15 \leq \beta \leq 0.25
		\qquad \text{for live or relatively moist fuels}.
		\label{eq:beta_live}
	\end{equation}

	The fuel-moisture variable $M(x,y,t)$ is retained explicitly in
	Eq.~\eqref{eq:RH} because moisture conditions can substantially influence
	wildfire propagation. No fixed value or universal interval is imposed on
	$M$ in the general formulation. For the Pantanal case study, fuel moisture
	is treated as a prescribed model parameter, and its specific value is
	introduced together with the numerical experiments.
	
	The coefficient $\gamma(x,y)$ plays a different role: it represents an
	effective local baseline capability of the landscape to sustain fire
	propagation, accounting for vegetation, fuel, and other unresolved local
	effects not captured by the wind and moisture terms. The vegetation
	analysis of Section~\ref{subsec:vegetation_data} justifies allowing
	$\gamma$ to vary in space, but does not determine its numerical values.
	
	\subsection{Calibration of the Spatial Propagation Coefficient}
	\label{subsec:gamma_calibration}
	
	The spatially varying baseline coefficient $\gamma(x,y)$ was estimated separately for each daily experiment using an observation-guided, sector-wise procedure. The initial perimeter and additional spatial control locations were grouped according to their position relative to the principal zones of observed fire advance. For September 3, the perimeter contained 18 primary vertices, supplemented by 30 control locations distributed across eight directional sectors (north, northeast, east, southeast, south, southwest, west, and northwest). For September 4, control locations were classified by their proximity to VIIRS active-fire pixels and subsequently refined in selected northern, western, eastern, and southern sectors. For September 5, control locations distinguished advancing sectors from relatively quiescent boundary segments. The final fields were defined from 48, 35, and 35 control locations, respectively.
	
	Calibration used the first 12 hourly propagation steps. Candidate coefficients were adjusted iteratively by comparing the simulated front envelope with the available VIIRS active-fire pattern. Values were reduced in sectors where the simulated front extended excessively beyond detected activity, retained or increased where the modeled advance was insufficient, and set to zero at selected quiescent control locations. The adjustments were case-specific rather than the result of a single global optimization. The exact control-point coordinates and corresponding calibrated values
	$\gamma_i$ are archived in the input workbooks and are available as
	described in the Data Availability statement. For September 3, tested local coefficients included values between $0$ and $0.350\,\mathrm{km\,h^{-1}}$; for September 4, proximity-dependent factors $\{1.00,0.60,0.25,0.08\}$ and sectoral damping factors were considered; and for September 5, a reduction of the initially estimated spread rates and zero-valued controls in quiescent sectors were used. Pixel-level overlap measures, including the intersection-over-union and Dice coefficients, were used to assess spatial agreement, alongside inspection of sectoral over- and under-propagation. The reported global multipliers of the calibrated field ($0.5$, $0.8$, $1.0$, $1.2$, and $2.0$) belong to the subsequent sensitivity assessment and should not be confused with the sector-wise calibration candidates.
	
	Given the calibrated control values $\gamma_i$ at locations $\mathbf{x}_i=(x_i,y_i)$, the continuous field was constructed by inverse-distance weighting with exponent two:
	\begin{equation}
		\gamma(\mathbf{x})=
		\frac{\displaystyle\sum_{i=1}^{N_c}\gamma_i\,d_i(\mathbf{x})^{-2}}
		{\displaystyle\sum_{i=1}^{N_c}d_i(\mathbf{x})^{-2}},
		\qquad
		d_i(\mathbf{x})=\|\mathbf{x}-\mathbf{x}_i\|,
		\label{eq:gamma_idw}
	\end{equation}
	with $\gamma(\mathbf{x}_i)=\gamma_i$ at a control location. The resulting field was sampled at the evolving perimeter vertices and held fixed during propagation steps 13--24. \rev{At a control point with $\gamma_i=0$, the local head- and backing-fire rates vanish. Mathematically, such points require a stationary update $\mathbf v=\mathbf 0$; the nonlinear trajectory equations otherwise involve the indeterminate ratio $0/0$. The stationary-point rule is the required computational convention whenever the local rates vanish; the numerical implementation must handle these zero-rate locations through this stationary convention.} Table~\ref{tab:gamma_calibration_summary} summarizes the final calibrated fields; the values describe the effective, event-specific propagation coefficient and are not direct measurements of fuel properties.
	
	\begin{table}[htbp]
		\centering
		\caption{Summary of the calibrated spatial propagation coefficients for the three daily experiments.}
		\label{tab:gamma_calibration_summary}
		\begin{tabular}{lccc}
			\toprule
			Experiment & Control locations & $\gamma$ range ($\mathrm{km\,h^{-1}}$) & Mean $\gamma$ ($\mathrm{km\,h^{-1}}$) \\
			\midrule
			September 3, 2020 & 48 & $0$--$0.3500$ & $0.1173$ \\
			September 4, 2020 & 35 & $0$--$0.3039$ & $0.1086$ \\
			September 5, 2020 & 35 & $0$--$0.2679$ & $0.1116$ \\
			\bottomrule
		\end{tabular}
	\end{table}
	
	\paragraph{Physical Interpretation and Backing-Fire Rate.}
	
	Equation~\eqref{eq:RH} retains the physical content of established
	surface-fire theory: the Rothermel model
	\cite{rothermel1972mathematical,andrews2018rothermel} relates propagation
	to fuel-bed characteristics and environmental conditions, but its direct
	pointwise application would require detailed fuel-bed data that are
	unavailable over the study domain. Equation~\eqref{eq:RH} therefore
	represents the wind and moisture effects explicitly and delegates the
	unresolved local differences to $\gamma(x,y)$; it is a spatially adjusted
	rate-of-spread relation in the physical context of Rothermel-type theory,
	not a direct implementation of the full model.
	
	In addition to the head-fire rate $R_H$, propagation of the complete fire
	perimeter requires a representation of the slower spread occurring in the
	backing direction. The local backing-fire rate is therefore defined as
	\begin{equation}
		R_B(x,y,t)=\kappa R_H(x,y,t),
		\label{eq:RB_RH_relation}
	\end{equation}
	with backing-to-head ratio $\kappa\in(0,1)$ treated as a case-specific
	model parameter. These two rates are the operational inputs of the
	trajectory model developed next. The ratio \(\kappa\) is prescribed separately for each experiment. The selected values are treated as case-specific modeling assumptions rather than independently estimated physical quantities. Their use does not imply that these ratios are uniquely determined by the available satellite observations.

	\subsection{Fuel and Directional Propagation Parameters}
	\label{subsec:fuel_propagation_parameters}
	
	Published fuel descriptions provide useful references for interpreting the
	vegetation and fuel characteristics of the Pantanal, but they cannot be
	mapped directly to a single fuel class for the study area. For example, the
	bulk-density information reported by \cite{Brown1981BulkDensity}
	provides useful qualitative references for fine-grass and grass--shrub fuel
	strata. Similarly, the grassland, savanna, and wetland fuel environments
	discussed in Section~2 provide relevant analogies for the different types of
	combustible vegetation occurring in the region. However, the ERA5-Land
	vegetation analysis indicates substantial spatial heterogeneity within the
	case-study domain. Consequently, no single literature-based fuel category is
	assigned to the entire wildfire perimeter. Instead, this heterogeneity is
	represented through the spatially varying propagation coefficient
	$\gamma(x,y)$ introduced in Eq.~\eqref{eq:RH}.
	
	The directional propagation model distinguishes between head-fire and
	backing-fire spread. At each point of the evolving perimeter, the head-fire
	\rev{direction is the wind-aligned major-axis direction of the local propagation wavelet, whereas the actual outward velocity at a perimeter point is selected by the geometric orthogonality condition in Section~3.4. The backing-fire direction is opposite to the wind-aligned head-fire axis. These axis directions need not coincide with the outward perimeter normal or the computed propagation trajectory.} The corresponding local rates of spread are denoted by
	$R_H(x,y,t)$ and $R_B(x,y,t)$, respectively. More specifically, over a time
	interval $[t_0,t_1]$, the head-fire rate of spread represents the distance
	traveled by the fire front from its location at time $t_0$ in the head-fire
	direction, divided by the elapsed time $t_1-t_0$. Similarly, the backing-fire
	rate of spread represents the corresponding propagation rate in the opposite
	direction. 
	
	For each perimeter point and time step, the local directional propagation
	law is characterized by three parameters, $a$, $b$, and $c$. Once the local
	head- and backing-fire rates have been determined, \rev{the parameters $a$, $b$, and $c$ are calculated as}
	\begin{equation}
		a=
		\frac{R_H+R_B}
		{2\left(1+0.25\left(1+U^*\right)\right)},
		\qquad
		b=\frac{R_H+R_B}{2},
		\qquad
		c=\frac{R_H-R_B}{2}.
		\label{eq:bc}
	\end{equation}
	{The dimensionless factor $0.25$ in the definition of $a$ is a prescribed shape parameter controlling the transverse-to-longitudinal spread relation; it is not derived from $R_H$ and $R_B$ alone. It should be treated as a phenomenological modeling choice, not as a universal empirical constant, until justified or calibrated for the relevant fuel conditions.}
	
	The parameter $a$ characterizes the transverse component of the local
	directional propagation law and incorporates the influence of wind intensity
	through the normalized wind speed $1+U^*$ defined in Eq.~\eqref{eq:RH};
	$b$ is the mean of the local head- and backing-fire rates, and $c$ measures
	their directional asymmetry. Since $1+U^*$ is dimensionless, $a$ has the units
	of $R_H$; for fixed rates, increasing $1+U^*$ decreases $a$ relative to $b$,
	elongating the local propagation pattern along the head-fire direction. The parameters are evaluated locally at each perimeter point and
	updated as conditions change. These parameters together with wind direction define the
	directional characteristics required by the trajectory equations below.

	\subsection{Hourly Fire-Perimeter Propagation}
	\label{subsec:perimeter_propagation}
	
	The satellite observations provide the information used to construct the
	initial fire perimeter, which is represented computationally by an ordered
	set of points. For the trajectory calculations, the geographical coordinates
	of these points are transformed into a local Cartesian coordinate system so
	that propagation distances and vector components can be expressed in
	consistent distance units. The initial fire perimeter is denoted by
	\begin{equation}
		\mathcal{P}^{(0)}
		=
		\left\{
		P_i^{(0)}=(x_i^{(0)},y_i^{(0)}):i=1,\ldots,n
		\right\}.
		\label{eq:initial_perimeter}
	\end{equation}
	The calculated perimeter points are subsequently transformed back to
	geographical coordinates for visualization and comparison with the satellite
	observations.
	
	The trajectory calculations are implemented in Python 
	to generate all maps, plots, and comparison images. For the numerical calculations,
	the initial polygon is discretized so that neighboring computational points
	are separated by no more than $0.1$ km, and the propagation interval is one
	hour. At each propagation step, the environmental variables and corresponding
	local propagation parameters are evaluated at the current perimeter points
	before the trajectory equations are solved.
	
	The objective is to determine the evolution of the perimeter at successive
	hourly intervals. For each point $P_i^{(k)}$ of the perimeter at hour $k$,
	a propagation vector
	\begin{equation}
		\mathbf{v}_i^{(k)}
		=
		\left(v_{1,i}^{(k)},v_{2,i}^{(k)}\right)
	\end{equation}
	is determined as a local propagation velocity (in $\mathrm{km\,h^{-1}}$). \rev{For a time step $\Delta t$ measured in hours, its endpoint defines the next perimeter point:}
	\begin{equation}
		P_i^{(k+1)}
		=
		P_i^{(k)}+\rev{\Delta t\,}\mathbf{v}_i^{(k)}.
		\label{eq:point_update}
	\end{equation}
	\rev{All experiments use $\Delta t=1\,\mathrm{h}$, so this dimensional clarification leaves the reported one-hour numerical updates unchanged.} Applying this update to all perimeter points generates
	$\mathcal{P}^{(k+1)}$, which then becomes the input perimeter for the
	subsequent time step.
	
	The present study considers only outward propagation of the observed fire
	perimeter, consistent with the objective of modeling the expansion of the
	fire domain; interior processes (residual combustion, heterogeneous fuel
	consumption) are outside the scope of the formulation.
	
	To determine the propagation vector associated with each perimeter point,
	let
	\begin{equation}
		\mathbf{u}_i^{(k)}
		=
		\left(u_{1,i}^{(k)},u_{2,i}^{(k)}\right)
	\end{equation}
	denote the local tangent to the current perimeter at $P_i^{(k)}$, and let
	$\Theta_i^{(k)}$ denote the local direction of wind. For
	compactness, the indices $i$ and $k$ are omitted in the equations below.
	
	For a candidate propagation vector $\mathbf{v}=(v_1,v_2)$, define
	\begin{equation}
		r=\sqrt{v_1^2+v_2^2},
		\label{eq:r}
	\end{equation}
	and its orientation relative to $\Theta$ through
	\begin{equation}
		C=
		\frac{v_1\cos\Theta+v_2\sin\Theta}{r},
		\qquad
		S_{\Theta}=
		\frac{v_2\cos\Theta-v_1\sin\Theta}{r}.
		\label{eq:relative_direction}
	\end{equation}
	The quantities entering the propagation equations are
	\begin{equation}
		D=a^2C^2+b^2S_{\Theta}^2,
		\label{eq:D}
	\end{equation}
	\begin{equation}
		Q=a^2C^2+(b^2-c^2)S_{\Theta}^2,
		\label{eq:Q}
	\end{equation}
	and
	\begin{equation}
		N=ca^2C+ab\sqrt{Q},
		\label{eq:N}
	\end{equation}
	where $a$, $b$, and $c$ are the local propagation parameters described in
	the preceding subsection.
	
	\rev{Following the Randers-based directional propagation formulation of \cite{dehkordi2022applications}, we introduce the directional function $H(\varphi)=N(\varphi)/D(\varphi)$, where $\varphi$ is the angle between the candidate propagation vector and the local wind direction, so that $\cos\varphi=C$ and $\sin\varphi=S_{\Theta}$. Imposing the \rev{one-hour-normalized velocity} relation $r=H(\varphi)$ and substituting Eqs.~\eqref{eq:D}--\eqref{eq:N} gives the first nonlinear condition,}
	\begin{equation}
		\frac{rD}{N}=1,
		\label{eq:trajectory_condition1}
	\end{equation}
	which determines the admissible propagation magnitude according to the
	local directional spread relation.
	
	A second condition is required to ensure that the propagated point is
	consistent with the local orientation of the current perimeter. Define
	\begin{equation}
		D'=2(b^2-a^2)S_{\Theta}C,
		\label{eq:Dprime}
	\end{equation}
	and
	\begin{equation}
		N'
		=
		-ca^2S_{\Theta}
		+
		ab\,
		\frac{(b^2-c^2-a^2)S_{\Theta}C}{\sqrt{Q}}.
		\label{eq:Nprime}
	\end{equation}
	The associated directional term is
	\begin{equation}
		H'=\frac{N'D-ND'}{D^2}.
		\label{eq:Hprime}
	\end{equation}
	\rev{The second nonlinear condition follows from Finsler orthogonality between the propagation direction and the tangent to the evolving perimeter, within the geometric framework of \cite{dehkordi2022applications, dehkordi2026satellite}. For the directional travel-time function $F(\mathbf{v})=r/H(\varphi)$, we impose $\partial F/\partial\mathbf{v}\cdot\mathbf{u}=0$. Differentiating with respect to the components of $\mathbf{v}$, using the directional derivative $H'$ in Eq.~\eqref{eq:Hprime}, and applying the first condition $r=H(\varphi)$ yields}
	\begin{equation}
		r(v_1u_1+v_2u_2)
		-
		H'(v_1u_2-v_2u_1)
		=
		0.
		\label{eq:trajectory_condition2}
	\end{equation}
	
	\rev{For $R_H=R_B=0$, the model defines $\mathbf v=\mathbf 0$ and the nonlinear equations must be bypassed; the stationary rule applies at zero-rate locations. For positive spread rates,} Equations~\eqref{eq:trajectory_condition1} and
	\eqref{eq:trajectory_condition2} are solved simultaneously in Python for
	$(v_1,v_2)$ at each discretized point of the current perimeter. Since only
	outward propagation is modeled, an admissible numerical solution must point
	toward the exterior of the current fire domain. If
	$\mathbf{n}_i^{(k)}$ denotes the local outward direction at
	$P_i^{(k)}$, the retained vector satisfies
	\begin{equation}
		\mathbf{v}_i^{(k)}
		\cdot
		\mathbf{n}_i^{(k)}
		>0.
		\label{eq:outward_condition}
	\end{equation}
	This condition excludes numerical solutions corresponding to propagation
	toward the interior of the current perimeter.
	
	\rev{For clarity, the directional law also admits a geometric characterization. In wind-aligned coordinates $(z_{\parallel},z_{\perp})$, the positive-rate indicatrix $r=H(\varphi)$ is the ellipse
		\begin{equation}
			\frac{(z_{\parallel}-c)^2}{b^2}+\frac{z_{\perp}^2}{a^2}=1,
			\label{eq:wind_ellipse}
		\end{equation}
		whose origin lies strictly inside it because $0<c<b$ when $0<R_B<R_H$. Its strict convexity implies a unique supporting point for each prescribed unit outward normal $\mathbf n$. Writing $\mathbf e=(\cos\Theta,\sin\Theta)$, $\mathbf e_{\perp}=(-\sin\Theta,\cos\Theta)$, $n_{\parallel}=\mathbf n\cdot\mathbf e$, and $n_{\perp}=\mathbf n\cdot\mathbf e_{\perp}$, this supporting velocity is
		\begin{equation}
			\mathbf v_{\mathrm{sup}}=c\mathbf e+\frac{b^2n_{\parallel}\mathbf e+a^2n_{\perp}\mathbf e_{\perp}}{\sqrt{b^2n_{\parallel}^2+a^2n_{\perp}^2}}.
			\label{eq:support_velocity}
		\end{equation}
		This expression characterizes the geometrically selected outward branch and provides an analytical consistency check for solutions of Eqs.~\eqref{eq:trajectory_condition1}--\eqref{eq:trajectory_condition2}. It does not assert that the numerical implementation uses this closed form. The supporting point is unique for each prescribed outward unit normal because the ellipse is strictly convex. Moreover, its outward projection is strictly positive: $\mathbf v_{\mathrm{sup}}\cdot\mathbf n=c n_{\parallel}+\sqrt{b^2n_{\parallel}^2+a^2n_{\perp}^2}>0$, since $b>c\geq0$ and $a>0$. Thus, for positive rates and a well-defined outward normal, the geometric formulation establishes existence and uniqueness of the outward supporting velocity.} Here, in the numerical result, (i) the general nonlinear solver is implemented so the computational framework can accommodate general anisotropic indicatrices $H(\varphi)$ where closed-form supporting points do not exist, and (ii) for the translated-ellipse case, Eq. (\ref{eq:support_velocity}) guarantees theoretical existence and uniqueness of the outward root selected by the solver.
	
	Once the outward vector has been determined,
	Eq.~\eqref{eq:point_update} provides the propagated point. Repeating the
	procedure along the complete perimeter produces the next hourly fire
	perimeter. The environmental variables and corresponding local propagation
	parameters are then updated before the subsequent propagation step. Starting
	from the initial satellite-derived perimeter, the recursive procedure
	therefore generates
	\begin{equation}
		\mathcal{P}^{(0)}
		\longrightarrow
		\mathcal{P}^{(1)}
		\longrightarrow
		\mathcal{P}^{(2)}
		\longrightarrow
		\cdots
		\longrightarrow
		\mathcal{P}^{(T)},
		\label{eq:perimeter_sequence}
	\end{equation}
	where $\mathcal{P}^{(k)}$ denotes the calculated fire perimeter after $k$
	hourly propagation steps.

	\section{Results}
	
	This section evaluates the proposed flexible data-assisted propagation
	framework on the Pantanal case of
	Section~2: three daily experiments combine satellite-derived perimeters
	with the vegetation and wind data, and the calculated propagation is
	compared with the corresponding satellite observations.
	\label{sec:results}
	
	\subsection{Pantanal Wildfire Event and Experimental Design}
	\label{subsec:case_study_results}
	
	The numerical experiments consider the wildfire observed in the Brazilian
	Pantanal from September~3 to 5, 2020, through three consecutive daily
	experiments. Each experiment follows the estimate--freeze--propagate
	protocol. Starting from the satellite-derived initial perimeter, the
	available satellite-observed fire pattern guides the sector-wise estimation
	of $\gamma(x,y)$ based on the propagation calculated over the first 12 model
	hours. The resulting field is then held fixed, with no further adjustment
	during hours 13--24, while the hourly ERA5 wind forcing continues to vary.
	{Thus, the latter interval is a fixed-parameter computational continuation, not an independent temporal validation of subsequent propagation.} Section~\ref{subsec:pixel_metrics}
	quantifies the spatial agreement at the pixel level for both the complete
	daily simulation and the {fixed-parameter late-window development}.
	Table~\ref{tab:data_roles} summarizes the data sources and their respective
	roles in the proposed framework.
	
	\begin{table}[htbp]
		\centering
		\caption{Role of each data source in the present calibration and propagation experiments.}
		\label{tab:data_roles}
		\begin{tabular}{p{0.22\textwidth}p{0.30\textwidth}p{0.38\textwidth}}
			\toprule
			Data source & Role in the model & Interpretation in this study \\
			\midrule
			Initial daily VIIRS observation & Construction of the initial perimeter & Initial condition for each daily experiment \\
			VIIRS pattern used for calibration & Observation-guided sector-wise estimation of $\gamma(x,y)$
			over the first 12 model hours & Calibration information for each daily experiment \\
			Subsequent VIIRS observation & Comparison with propagation obtained using fixed $\gamma(x,y)$ & {Fixed-parameter model continuation; not independent temporal validation} \\
			Hourly ERA5 wind &  $1+U^*$ and $\Theta$ in Eq.~\eqref{eq:RH} & External atmospheric forcing during hourly propagation \\
			ERA5-Land vegetation &  Characterization of spatial vegetation heterogeneity & Motivation for a nonuniform coefficient field; not direct fuel calibration \\
			Prescribed moisture $M$ & Moisture attenuation in Eq.~\eqref{eq:RH} & Event-level assumption examined through structural sensitivity \\
			\bottomrule
		\end{tabular}
	\end{table}

	\subsection{Model Parameters and Spatial Propagation Coefficient}
	\label{subsec:parameters_results}
	
	The three numerical experiments used $\alpha=1.00$, $\beta=0.08$, and
	$U_{\mathrm{ref}}=1~\mathrm{m\,s^{-1}}$. The choices of $\alpha$ and
	$\beta$ provide physically interpretable baseline values for the numerical
	experiments. In particular, $\alpha=1.00$ represents a linear dependence of
	the head-fire rate of spread on the normalized wind speed, whereas
	$\beta=0.08$ introduces a moderate exponential reduction in the rate of
	spread as fuel moisture increases. These values are adopted as reference
	values rather than as uniquely calibrated parameters for the Pantanal.
	Their influence on the calculated fire propagation is examined separately
	through the sensitivity analysis. The backing-to-head ratio was case-specific:
	$\kappa=0.1875$, $0.125$, and $0.25$ for model inputs 3, 4, and 5,
	respectively. These values correspond to
	backing-fire rates of $18.75\%$, $12.5\%$, and $25\%$ of the respective
	head-fire rates. The use of a backing-to-head ratio is consistent with
	established directional fire-spread formulations, in which the backing-fire
	rate is expressed as a fraction of the head-fire rate
	\cite{andrews2018rothermel}. In the present study, the three values are
	prescribed case-specific parameters rather than independently calibrated
	backing-fire rates. Fuel moisture was
	prescribed as $M=10\%$ in all three cases. 
	This relatively low moisture level was selected to represent the unusually
	dry conditions prevailing
	in the Pantanal during 2020, particularly in September, when the region
	experienced severe drought and pronounced water deficits \cite{pelissari2023dynamics}. Across the three experiments, the hourly wind speed ranged
	from $0.19$ to $3.04\,\mathrm{m\,s^{-1}}$
	(mean $1.98\,\mathrm{m\,s^{-1}}$) on September 3,
	from $1.50$ to $2.56\,\mathrm{m\,s^{-1}}$
	(mean $2.11\,\mathrm{m\,s^{-1}}$) on September 4,
	and from $1.52$ to $3.44\,\mathrm{m\,s^{-1}}$
	(mean $2.49\,\mathrm{m\,s^{-1}}$) on September 5
	(see Figs.~2b and~2c).
	
	\begin{figure}[htbp]
		\centering
		\begin{subfigure}[t]{0.38\textwidth}
			\centering
			\includegraphics[width=\linewidth]{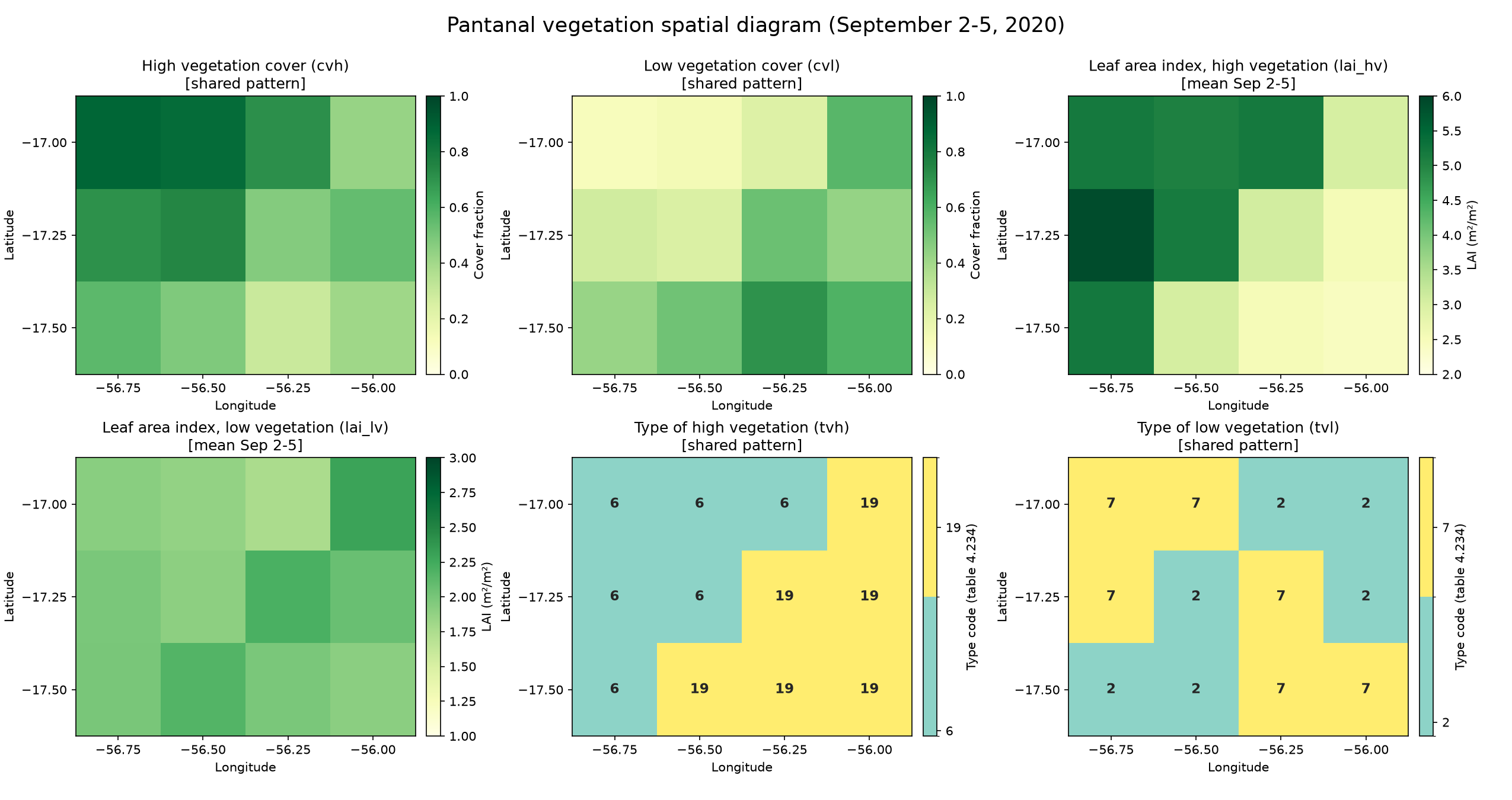}
			\caption{Vegetation spatial distribution.}
			\label{fig:vegetation_spatial}
		\end{subfigure}
		\hfill
		\begin{subfigure}[t]{0.3\textwidth}
			\centering
			\includegraphics[width=\linewidth]{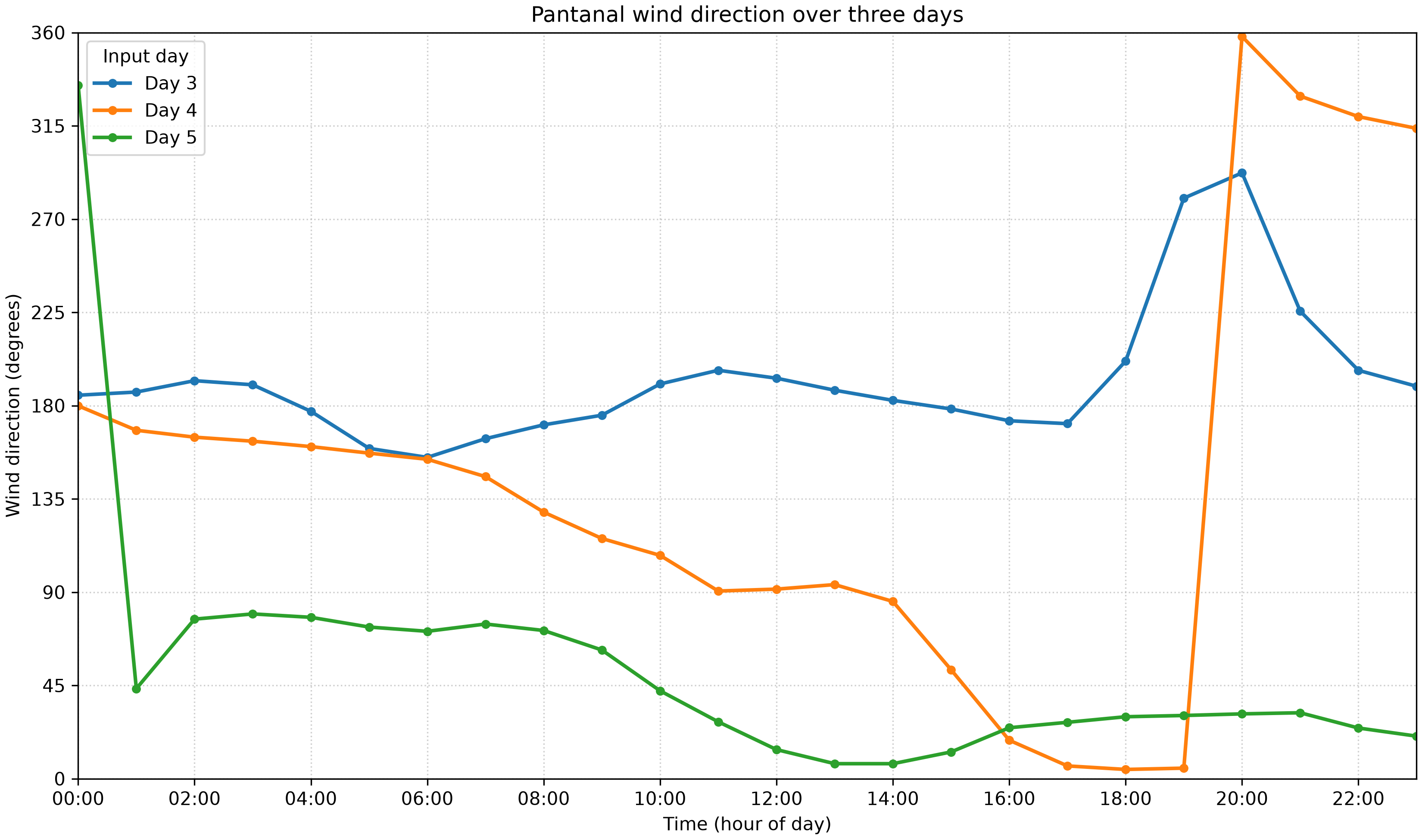}
			\caption{Wind direction.}
			\label{fig:wind_direction}
		\end{subfigure}
		\hfill
		\begin{subfigure}[t]{0.3\textwidth}
			\centering
			\includegraphics[width=\linewidth]{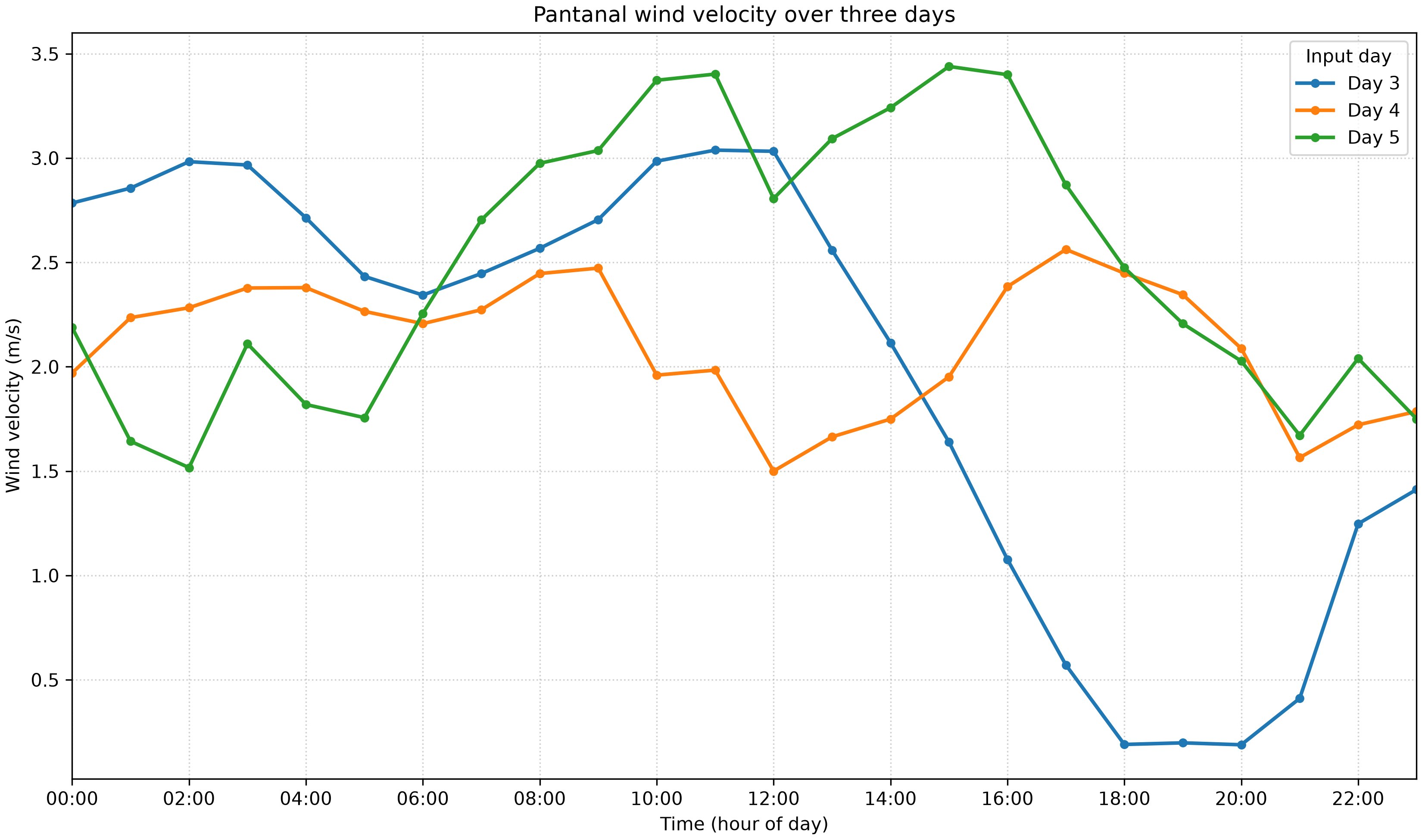}
			\caption{Wind velocity.}
			\label{fig:wind_velocity}
		\end{subfigure}
		\caption{Environmental information used to characterize the Pantanal case
			study: (a) spatial vegetation conditions over the analyzed dates; (b) wind
			direction during the three-day period; and (c) corresponding wind velocity.}
		\label{fig:environmental_results}
	\end{figure}

	Spatial variability in fire propagation was represented through the local
	baseline coefficient $\gamma(x,y)$.
	
\begin{figure}[htbp]
		\centering
		\begin{subfigure}[t]{0.18\textwidth}
			\centering
			\includegraphics[width=.8\linewidth]{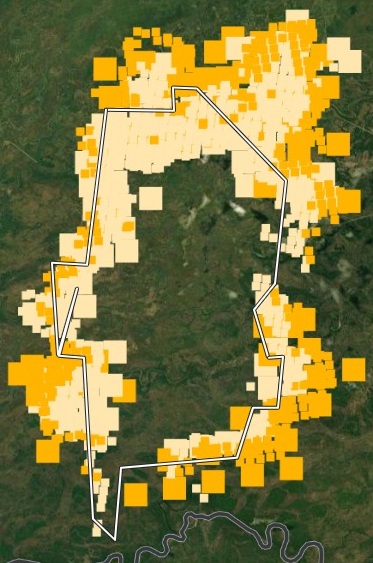}
			\caption{Observed pattern.}
			\label{fig:sep03_calib_obs}
		\end{subfigure}\hfill
		\begin{subfigure}[t]{0.30\textwidth}
			\centering
			\includegraphics[width=\linewidth]{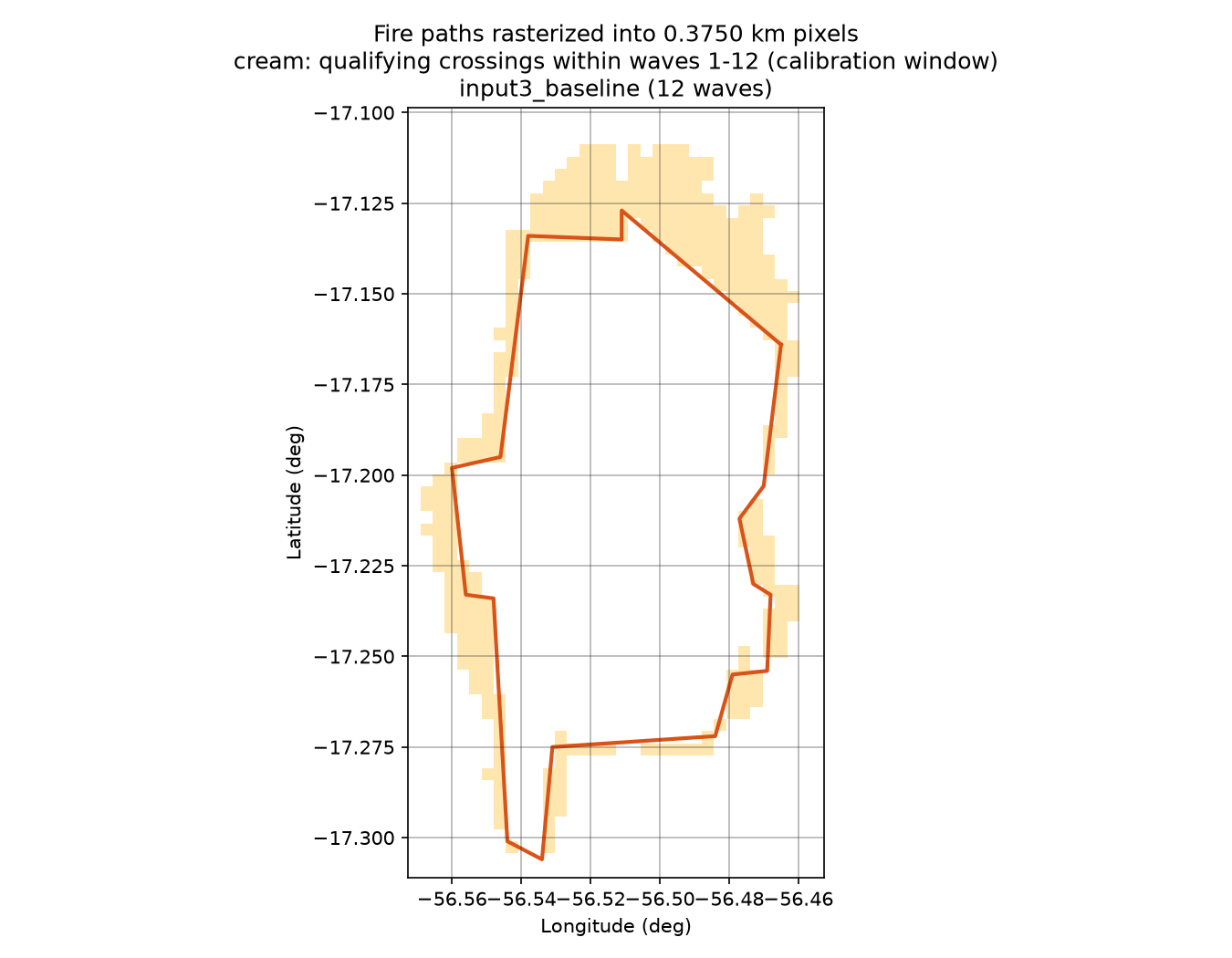}
			\caption{12-wave model calibration.}
			\label{fig:sep03_gamma_window}
		\end{subfigure}\hfill
		\begin{subfigure}[t]{0.30\textwidth}
			\centering
			\includegraphics[width=\linewidth]{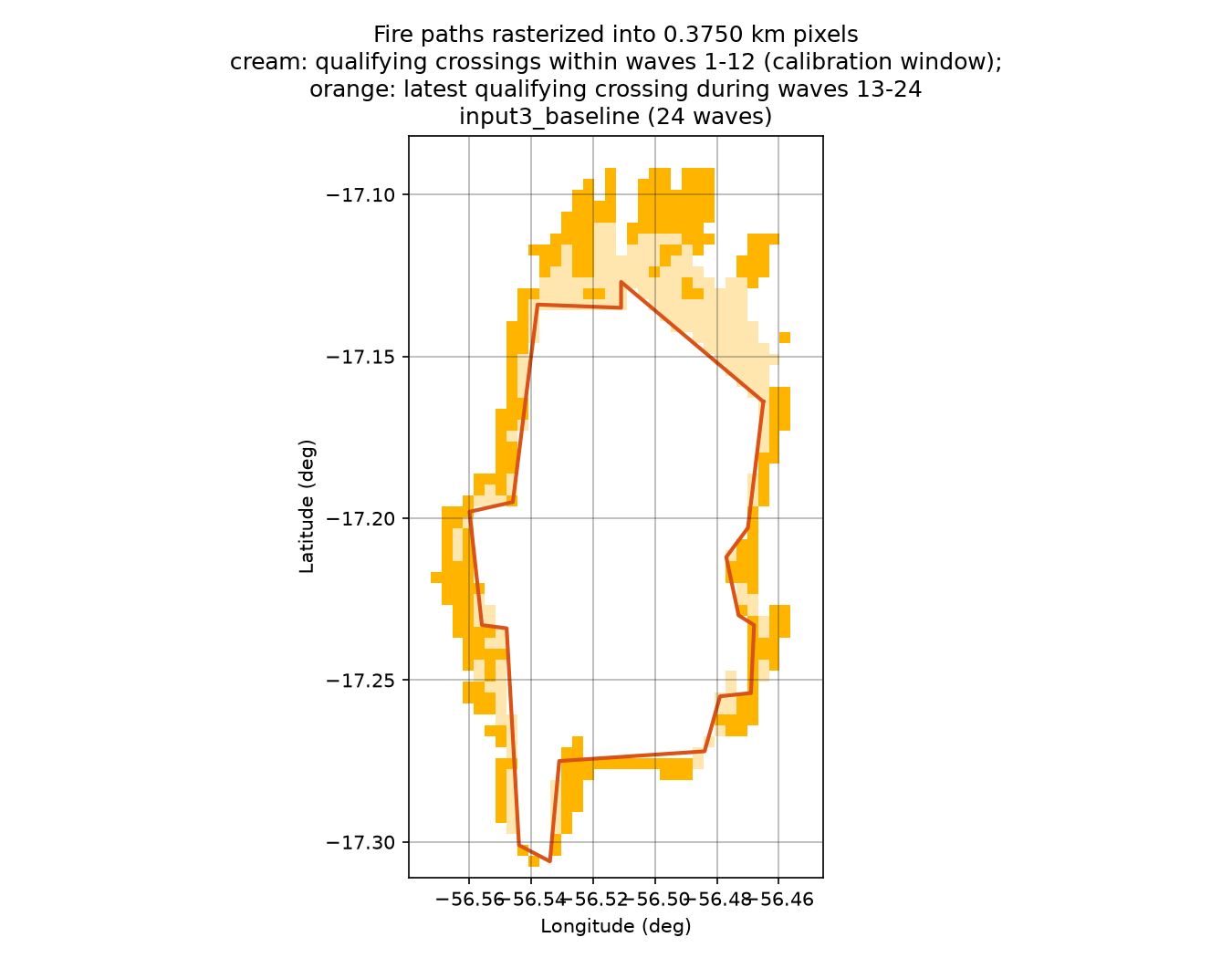}
			\caption{Complete 24-h calculation.}
			\label{fig:sep03_gamma_full}
		\end{subfigure}
		
		\begin{subfigure}[t]{0.21\textwidth}
			\centering
			\includegraphics[width=.8\linewidth]{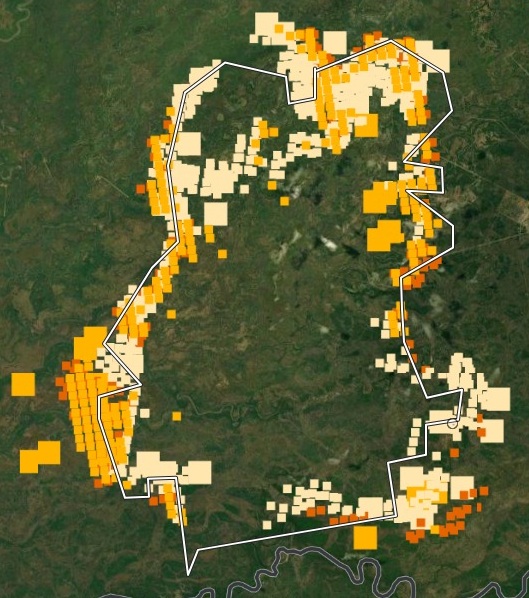}
			\caption{Observed pattern.}
			\label{fig:sep04_calib_obs}
		\end{subfigure}\hfill
		\begin{subfigure}[t]{0.30\textwidth}
			\centering
			\includegraphics[width=\linewidth]{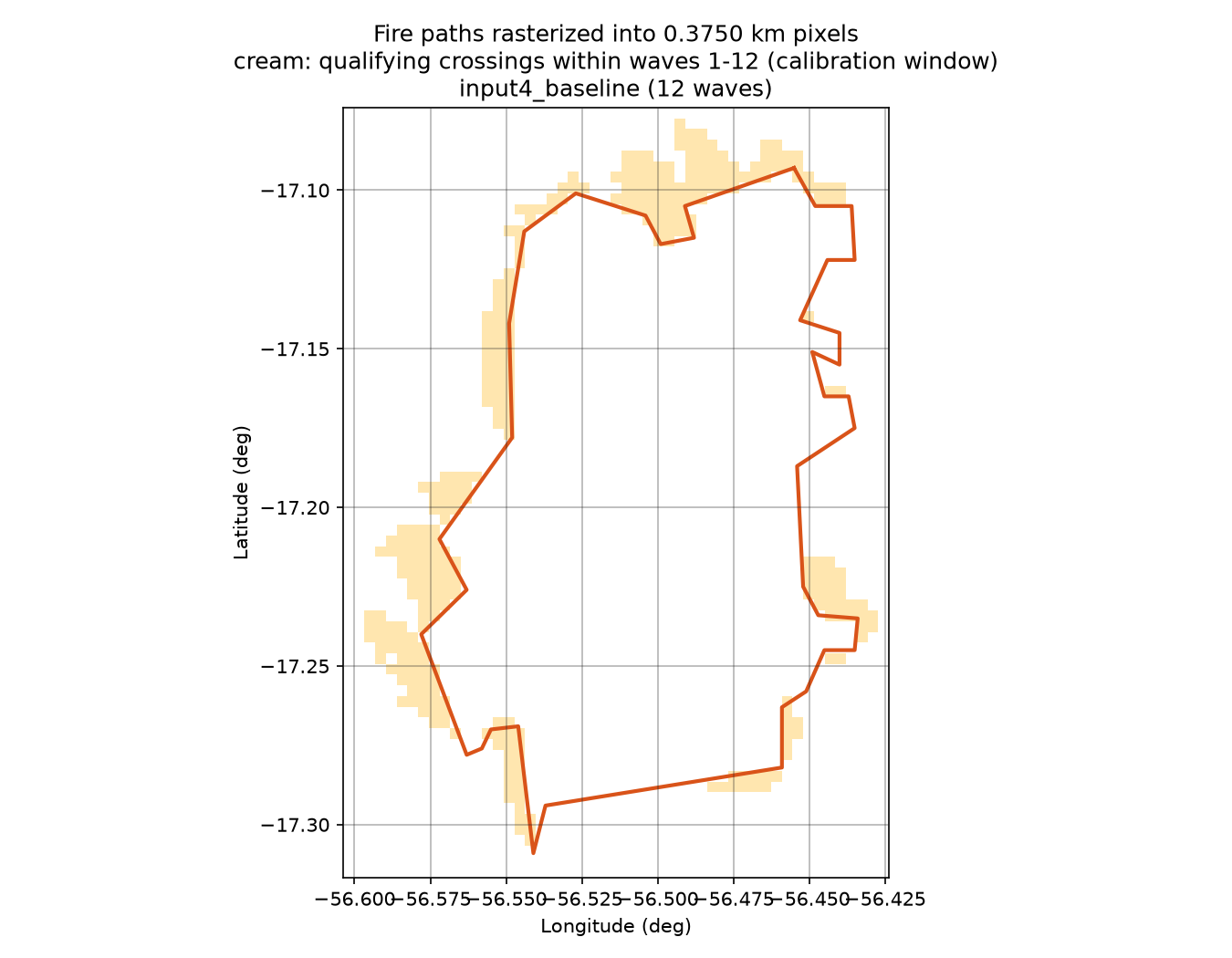}
			\caption{12-wave calibration window.}
			\label{fig:sep04_gamma_window}
		\end{subfigure}\hfill
		\begin{subfigure}[t]{0.30\textwidth}
			\centering
			\includegraphics[width=\linewidth]{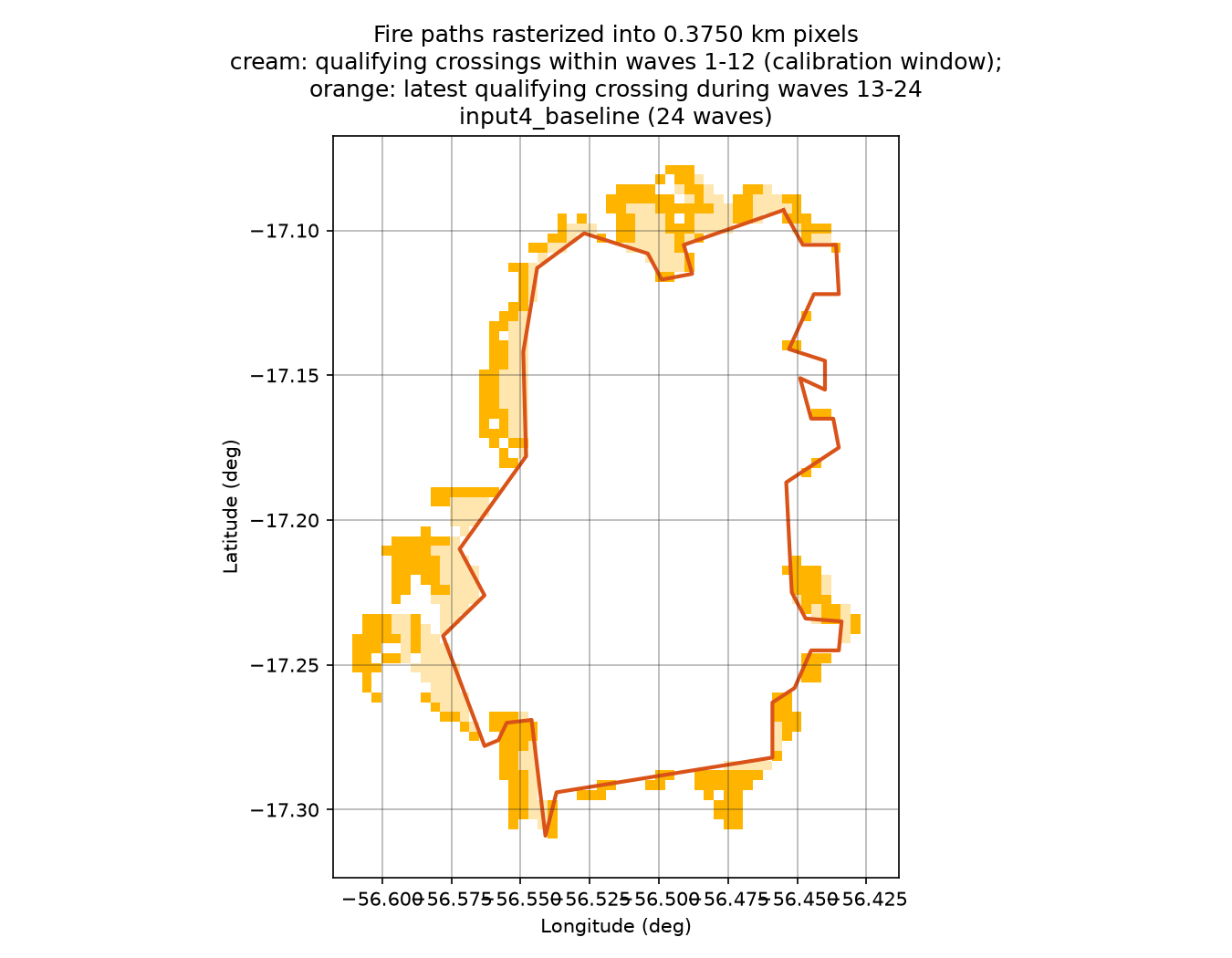}
			\caption{Complete 24-h calculation.}
			\label{fig:sep04_gamma_full}
		\end{subfigure}
		
		\begin{subfigure}[t]{0.22\textwidth}
			\centering
			\includegraphics[width=.8\linewidth]{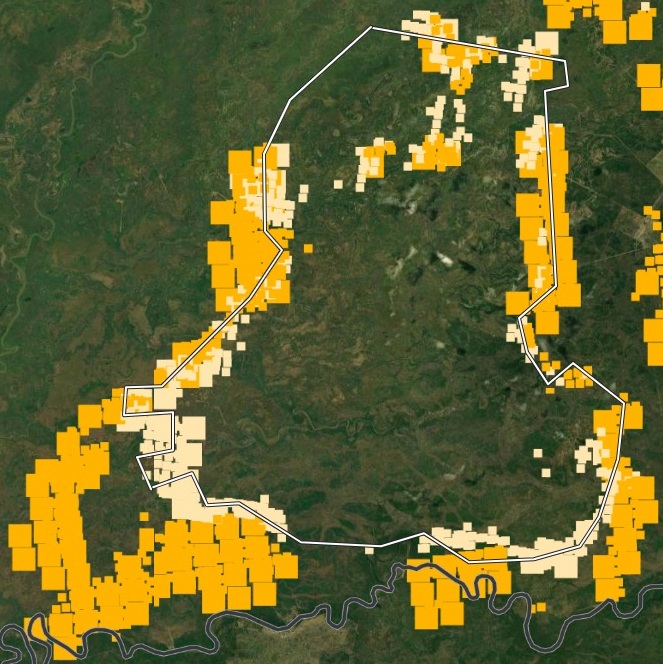}
			\caption{Observed pattern.}
			\label{fig:sep05_calib_obs}
		\end{subfigure}\hfill
		\begin{subfigure}[t]{0.30\textwidth}
			\centering
			\includegraphics[width=\linewidth]{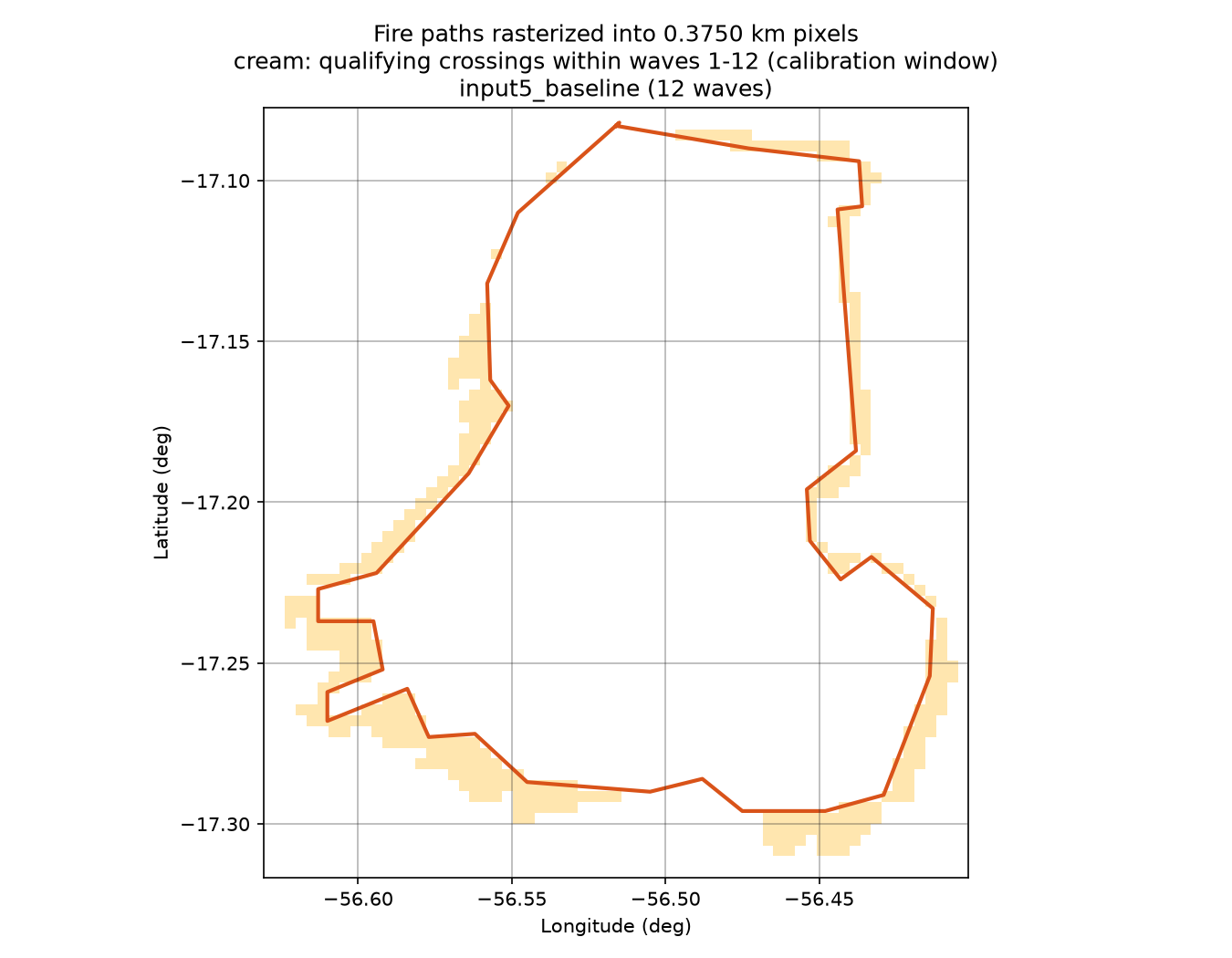}
			\caption{12-wave calibration window.}
			\label{fig:sep05_gamma_window}
		\end{subfigure}\hfill
		\begin{subfigure}[t]{0.30\textwidth}
			\centering
			\includegraphics[width=\linewidth]{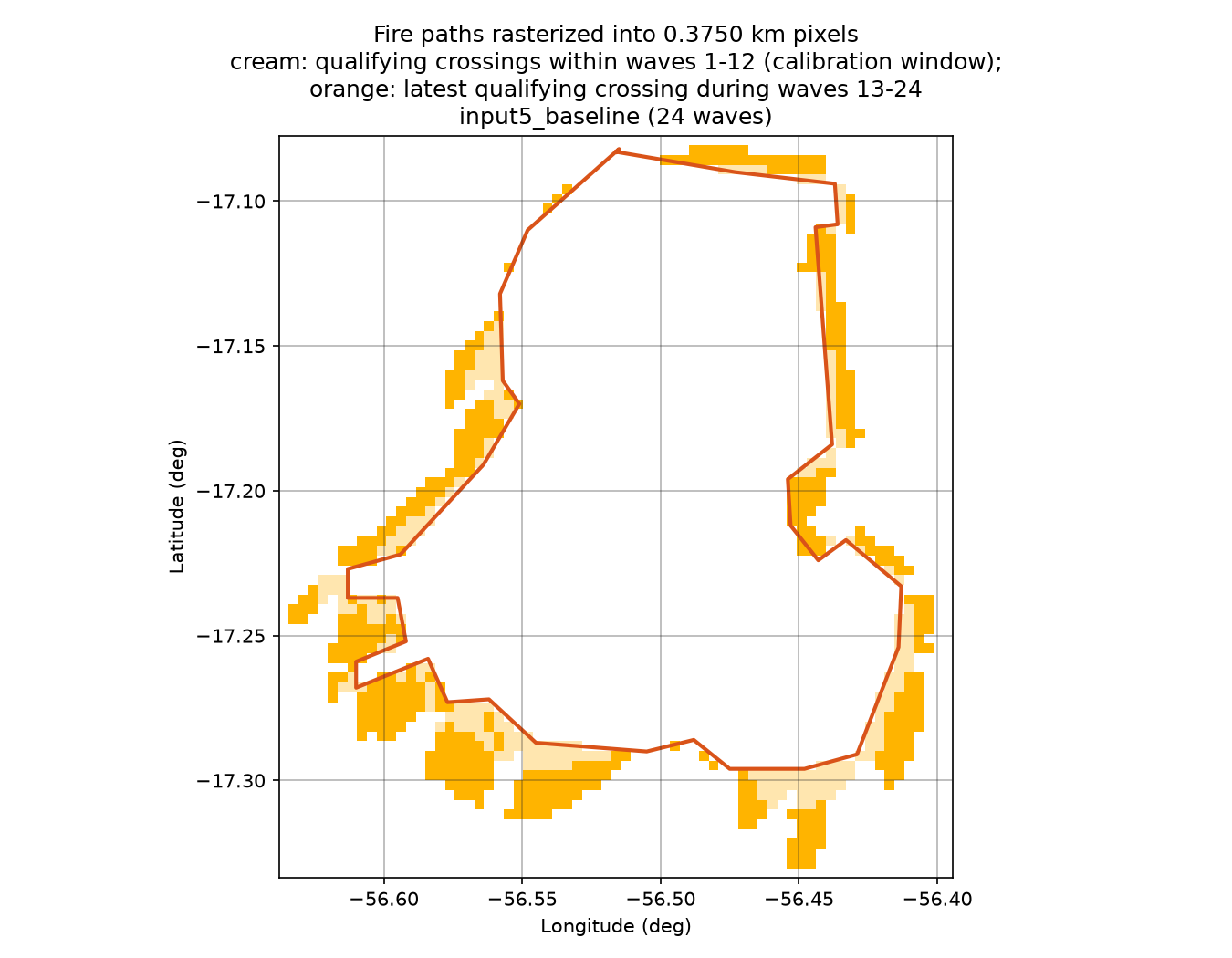}
			\caption{Complete 24-h calculation.}
			\label{fig:sep05_gamma_full}
		\end{subfigure}
		
		\caption{Determination of the local propagation coefficient $\gamma(x,y)$
			for the September~3, 4, and 5 experiments (rows). Left column: available
			VIIRS active-fire pattern displayed over the satellite background. Middle
			column: propagation calculated over the 12-h model calibration window
			(cream: pixels associated with waves 1--12). Right column: complete 24-h
			calculation with the calibrated $\gamma(x,y)$ field held fixed (orange:
			pixels whose latest qualifying crossing occurred during waves 13--24).
			The middle-column pattern is compared with the available satellite
			observations to guide the selection of sector-wise $\gamma$ values, whereas
			the subsequent propagation is generated without further adjustment of
			$\gamma(x,y)$ and provides {a fixed-parameter continuation rather than independent temporal validation}.
			The dark-orange line denotes the initial satellite-derived perimeter.}
		\label{fig:gamma_calibration}
	\end{figure}

	Figure~\ref{fig:gamma_calibration} illustrates the procedure in the
	rasterized pixel representation,
	which places the calculated front in the same georeferenced form as the
	satellite observations. For each daily experiment, the left column shows
	the observed active-fire pattern (satellite background)  and the middle column shows the 12-wave
	calculation restricted to the calibration window (cream).
	Because each daily VIIRS overpass provides only a single spatial
	snapshot, no time-resolved variation of the coefficient can be
	identified from the data; $\gamma(x,y)$ is therefore treated as a
	time-invariant effective field anchored on this first 12-h window.
	The right column of
	Fig.~\ref{fig:gamma_calibration} shows the complete 24-h calculation:
	the cream pixels reproduce the calibration-window envelope of the
	middle column, while the orange pixels (the wave color of waves
	13--24) mark the area whose latest qualifying crossing
	occurred during the remaining 12~h. This orange advance is produced by
	the frozen field under the observed hourly wind forcing, plays no role
	in the estimation, and is strongly sector-dependent in all three days;
	{it represents a fixed-parameter continuation, not an independent prediction of observations withheld from calibration.} The two-stage display thus documents, on the same pixel grid, both the 12-h calibration-window envelope used to estimate $\gamma(x,y)$   and the subsequent 12-h fixed-parameter continuation..
	
	The selected values should therefore be interpreted as effective local
	propagation coefficients for the present case study rather than as direct
	measurements of fuel properties or parameters associated uniquely with
	particular vegetation classes. The use of spatially varying values is
	supported independently by the substantial heterogeneity in vegetation type,
	fractional cover, and LAI identified from the ERA5-Land data
	(Fig.~\ref{fig:vegetation_spatial}), whereas the numerical magnitudes of
	$\gamma(x,y)$ were estimated by comparing the available satellite-observed
	fire pattern with the propagation calculated over the first 12 model hours
	of each daily experiment. Accordingly, the vegetation data justify
	the spatial structure of the coefficient field, while the satellite-observed
	fire evolution provides the information used to calibrate its numerical
	values. The sector values
	were interpolated over the geographic domain and the interpolated field was
	sampled at the current perimeter points during each numerical step.

\begin{figure}[htbp]
		\centering
		\begin{subfigure}[t]{0.34\textwidth}
			\centering
			\includegraphics[width=\linewidth]{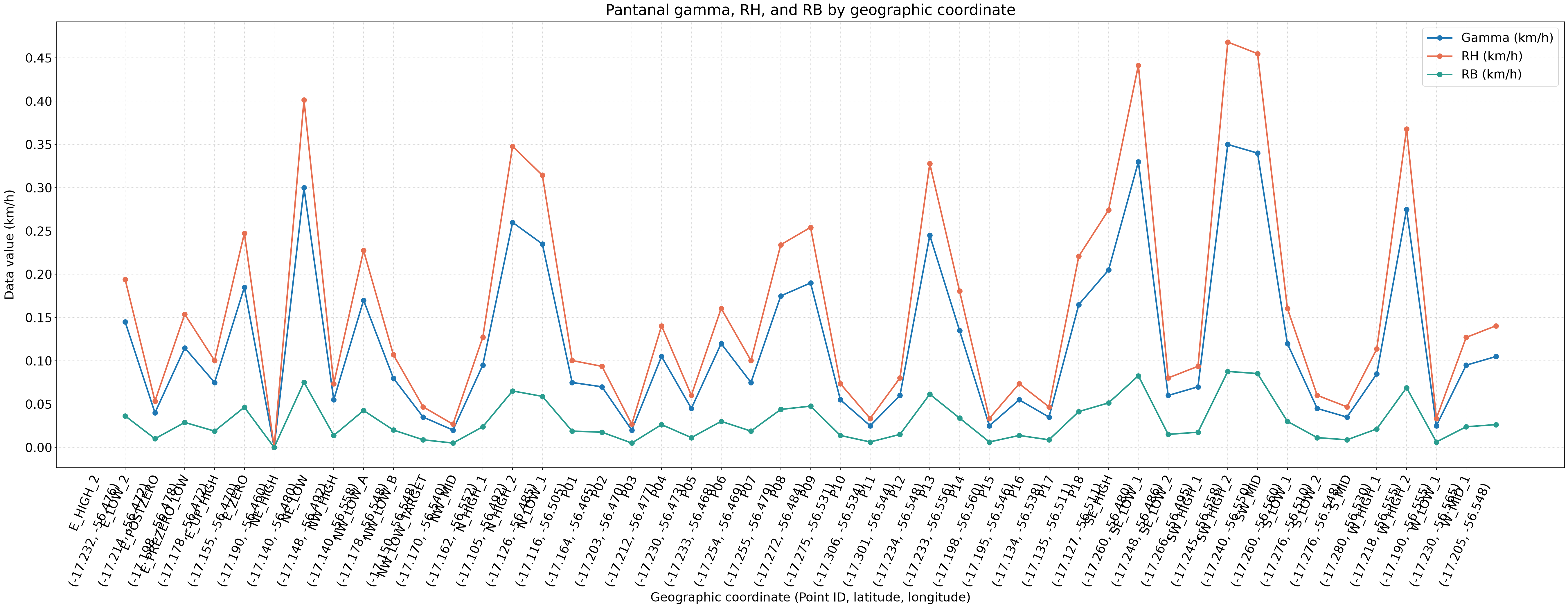}
			\caption{September~3.}
			\label{fig:gamma_rh_rb_input3}
		\end{subfigure}
		\hfill
		\begin{subfigure}[t]{0.29\textwidth}
			\centering
			\includegraphics[width=\linewidth]{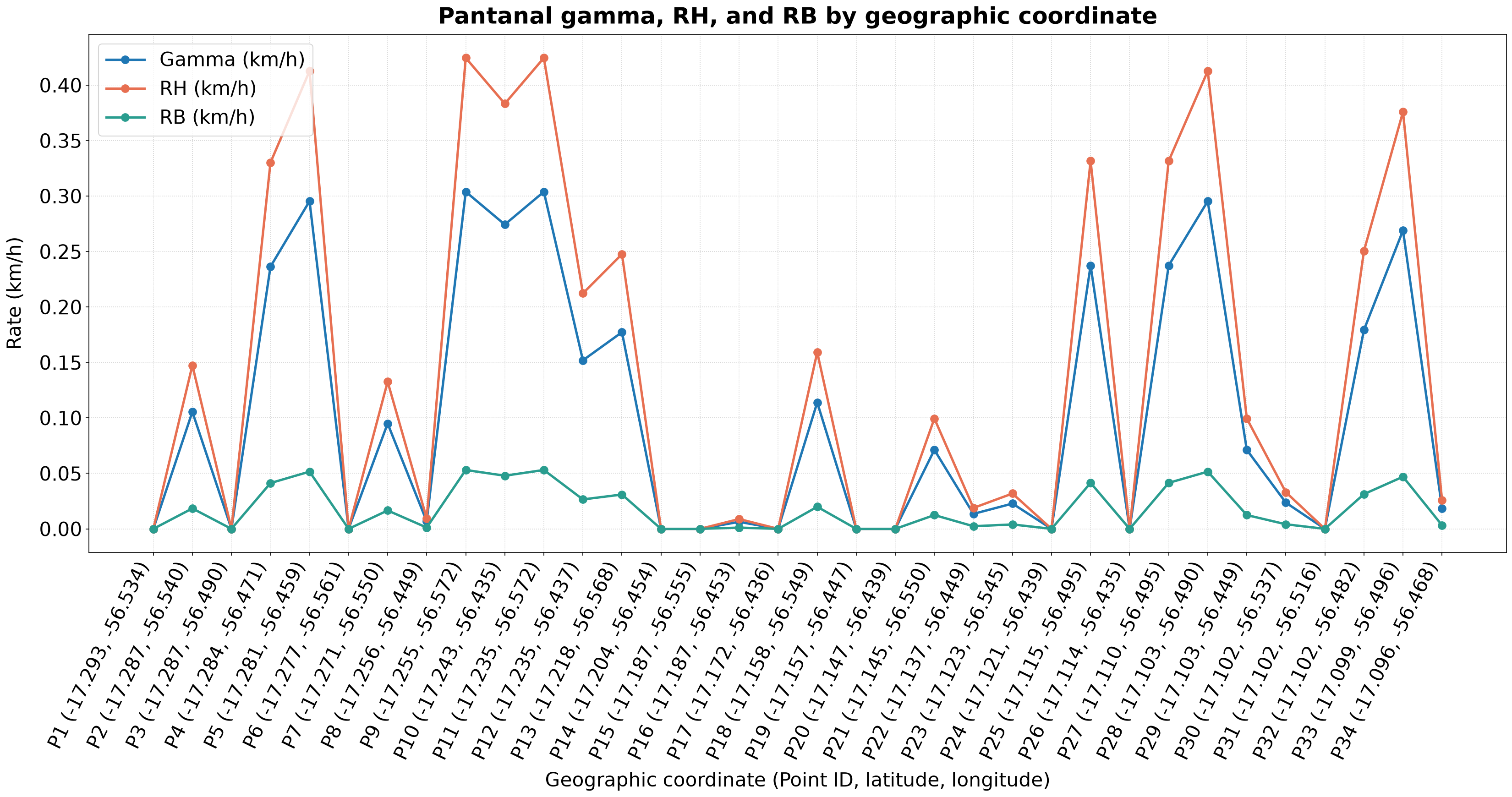}
			\caption{September~4.}
			\label{fig:gamma_rh_rb_input4}
		\end{subfigure}
		\hfill
		\begin{subfigure}[t]{0.34\textwidth}
			\centering
			\includegraphics[width=\linewidth]{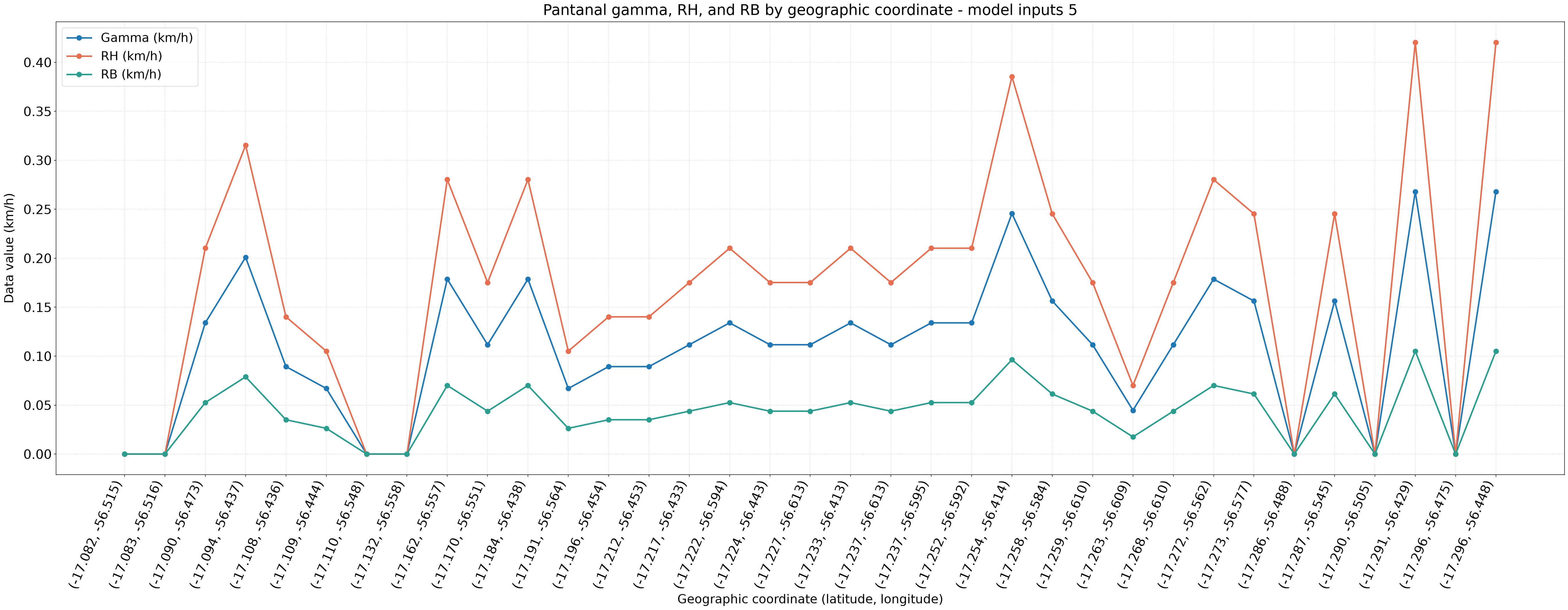}
			\caption{September~5.}
			\label{fig:gamma_rh_rb_input5}
		\end{subfigure}
		\caption{Spatial distributions of the observation-calibrated baseline
			propagation coefficient $\gamma(x,y)$ and of the local head- and
			backing-fire rates obtained by combining this frozen coefficient field
			with the environmental conditions in
			Fig.~\ref{fig:environmental_results} for the three model-input cases.
			These rates are the spatially varying quantities used in the
			trajectory-based propagation calculations.}
		\label{fig:gamma_rh_rb_results}
	\end{figure}

	In Fig.~\ref{fig:gamma_rh_rb_results}, three panels present the geographic distributions of $\gamma$, $R_H$,
	and $R_B$ for the three model-input cases. The spatial variation of
	$\gamma$ originates in the observation-guided calibration over the first 12 model hours,
	with the ERA5-Land vegetation analysis providing the independent
	justification for allowing spatial variation, whereas the time-dependent
	wind magnitude and direction enter the local propagation calculation
	through $1+U^*$ and $\Theta$. The resulting head- and backing-fire rates are the operational
	spread rates used by the trajectory formulation, with the backing rate
	remaining lower than the head-fire rate consistently with
	Eq.~\eqref{eq:RB_RH_relation}.
	
	\paragraph{Pixel Representation of the Calculated Trajectories}

	For comparison with the satellite observations, the calculated trajectories
	are rasterized onto a grid of $375$~m pixels, corresponding to the native
	resolution of the VIIRS active-fire product
	\cite{schroeder2014viirs}. This procedure
	converts the continuous trajectory output into a georeferenced representation
	comparable with the discrete VIIRS detections. A pixel is retained as simulated fire activity when it is crossed by at least three distinct calculated trajectories; \rev{``distinct'' refers to separate calculated trajectories, not three visits by the same trajectory. The threshold is a prescribed rasterization filter, not a VIIRS detection criterion. } This threshold suppresses isolated crossings and emphasizes locations of persistent calculated propagation. The resulting
	pixel maps are used for visual and sector-wise comparison and provide the
	basis for the quantitative spatial metrics described in
	Section~\ref{subsec:pixel_metrics}.
	
	\subsection{Fire Propagation on September 3, 2020}
	\label{subsec:sep03}
	
\begin{figure}[htbp]
		\centering
		
		\begin{subfigure}[t]{0.14\textwidth}
			\centering
			\includegraphics[width=.8\linewidth]{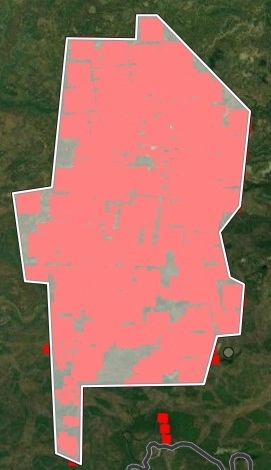}
			\caption{Initial satellite-derived perimeter.}
			\label{fig:sep03_initial}
		\end{subfigure}
		\hfill
		\begin{subfigure}[t]{0.32\textwidth}
			\centering
			\includegraphics[width=\linewidth]{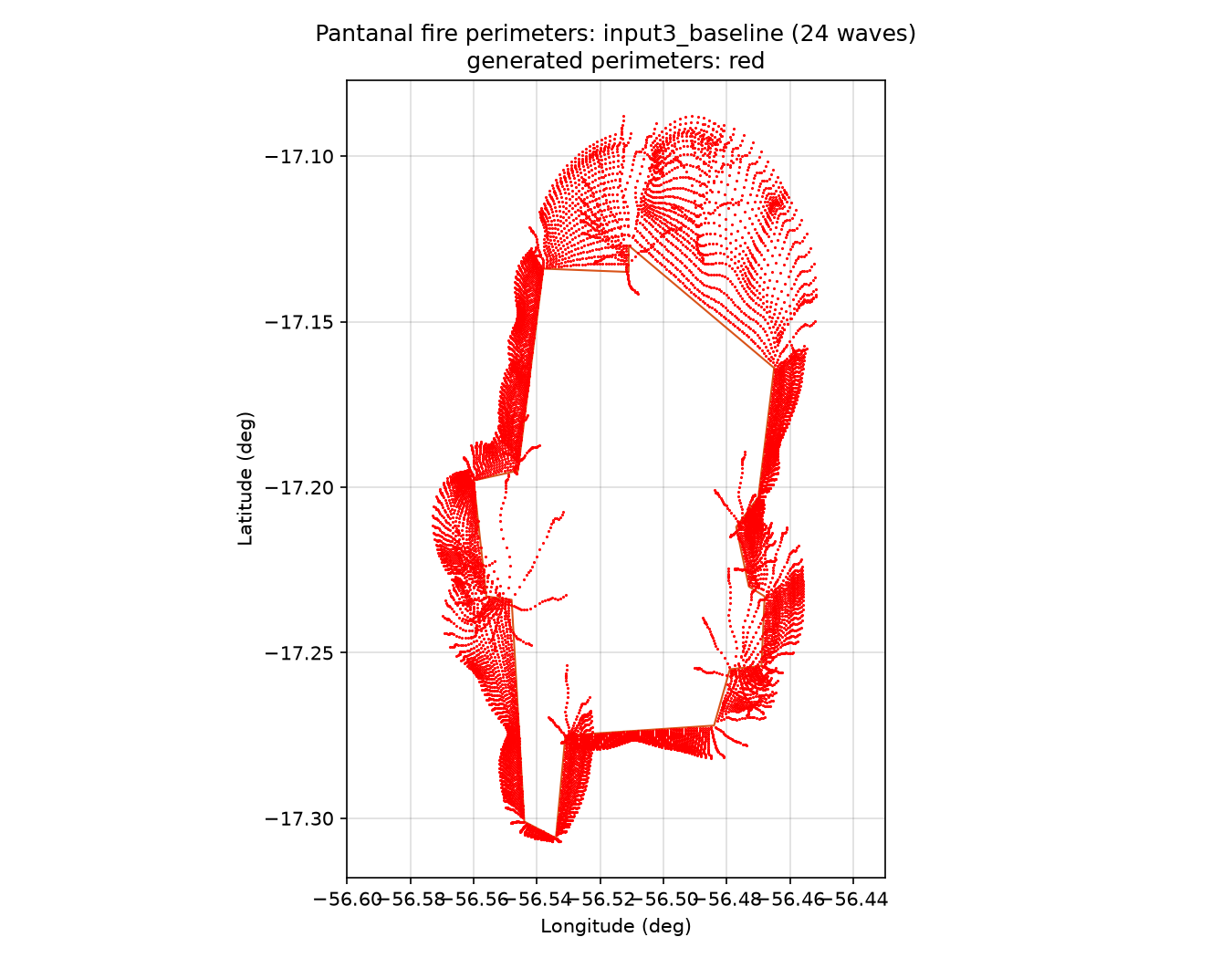}
			\caption{Calculated fire propagation.}
			\label{fig:sep03_simulation}
		\end{subfigure}
		\hfill
		\begin{subfigure}[t]{0.32\textwidth}
			\centering
			\includegraphics[width=\linewidth]{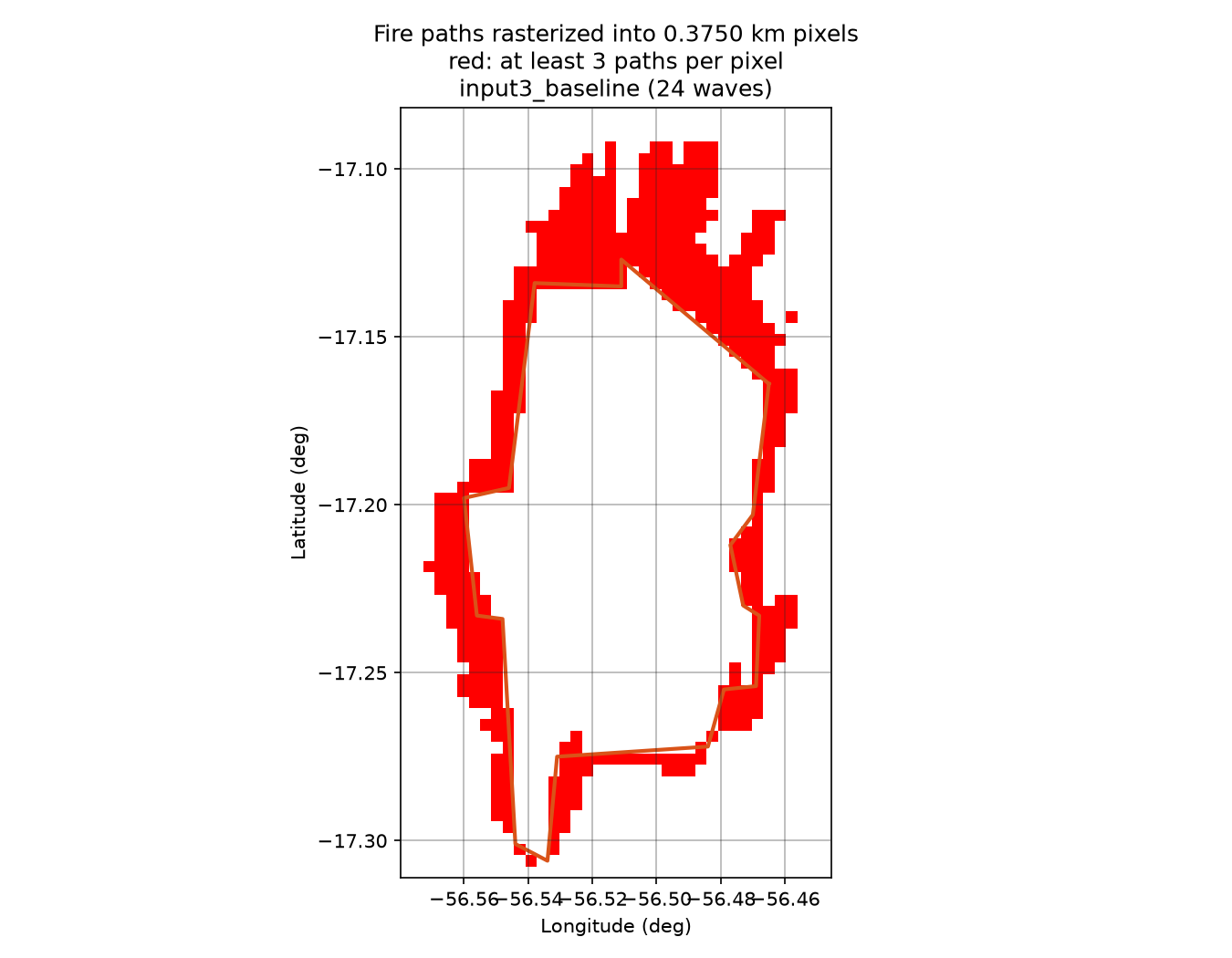}
			\caption{Calculated propagation rasterized into 375-m pixels.}
			\label{fig:sep03_simulation_pix}
		\end{subfigure}
		\hfill
		\begin{subfigure}[t]{0.18\textwidth}
			\centering
			\includegraphics[width=\linewidth]{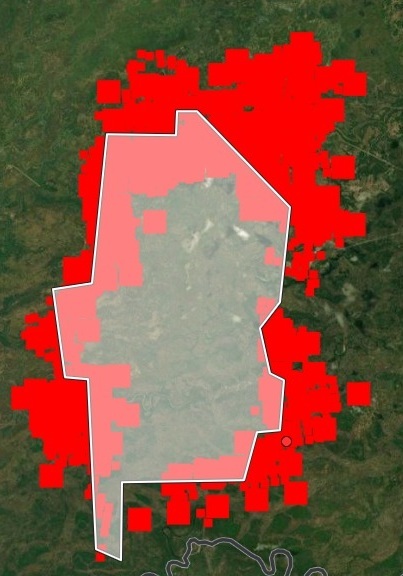}
			\caption{Subsequent satellite observation.}
			\label{fig:sep03_observed}
		\end{subfigure}
		
		\caption{Wildfire propagation experiment for September 3, 2020:
			(a) satellite observation used to construct the initial fire
			perimeter; (b) hourly fire propagation calculated by the proposed
			trajectory model; (c) the same calculated trajectories rasterized
			into 375-m pixels, displayed at the native resolution of the VIIRS
			active-fire product; and (d) subsequent satellite observation used
			for comparison.}
		\label{fig:sep03}
	\end{figure}

	Figure~\ref{fig:sep03} presents the first daily experiment: the
	satellite-derived boundary in Fig.~\ref{fig:sep03_initial} defines the
	initial perimeter, from which the hourly propagation of
	Fig.~\ref{fig:sep03_simulation} is calculated with the daily
	$\gamma(x,y)$ field calibrated over the first 12 model hours and then fixed. The
	same calculation is summarized at the observational scale in
	Fig.~\ref{fig:sep03_simulation_pix}.
	
	The calculated trajectories show strongly nonuniform propagation ---
	substantial displacement in some sectors, limited advance in others ---
	mirroring the heterogeneous pattern of the subsequent satellite
	detections. In the rasterized panel, the sectors of strongest
	calculated development appear as contiguous blocks of marked pixels,
	coinciding with the areas where the subsequent detections concentrate.

	\subsection{Fire Propagation on September 4, 2020}
	\label{subsec:sep04}
	
\begin{figure}[htbp]
		\centering
		
		\begin{subfigure}[t]{0.17\textwidth}
			\centering
			\includegraphics[width=\linewidth]{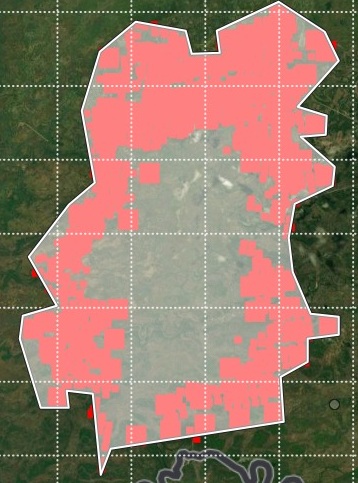}
			\caption{Initial satellite-derived perimeter.}
			\label{fig:sep04_initial}
		\end{subfigure}
		\hfill
		\begin{subfigure}[t]{0.29\textwidth}
			\centering
			\includegraphics[width=\linewidth]{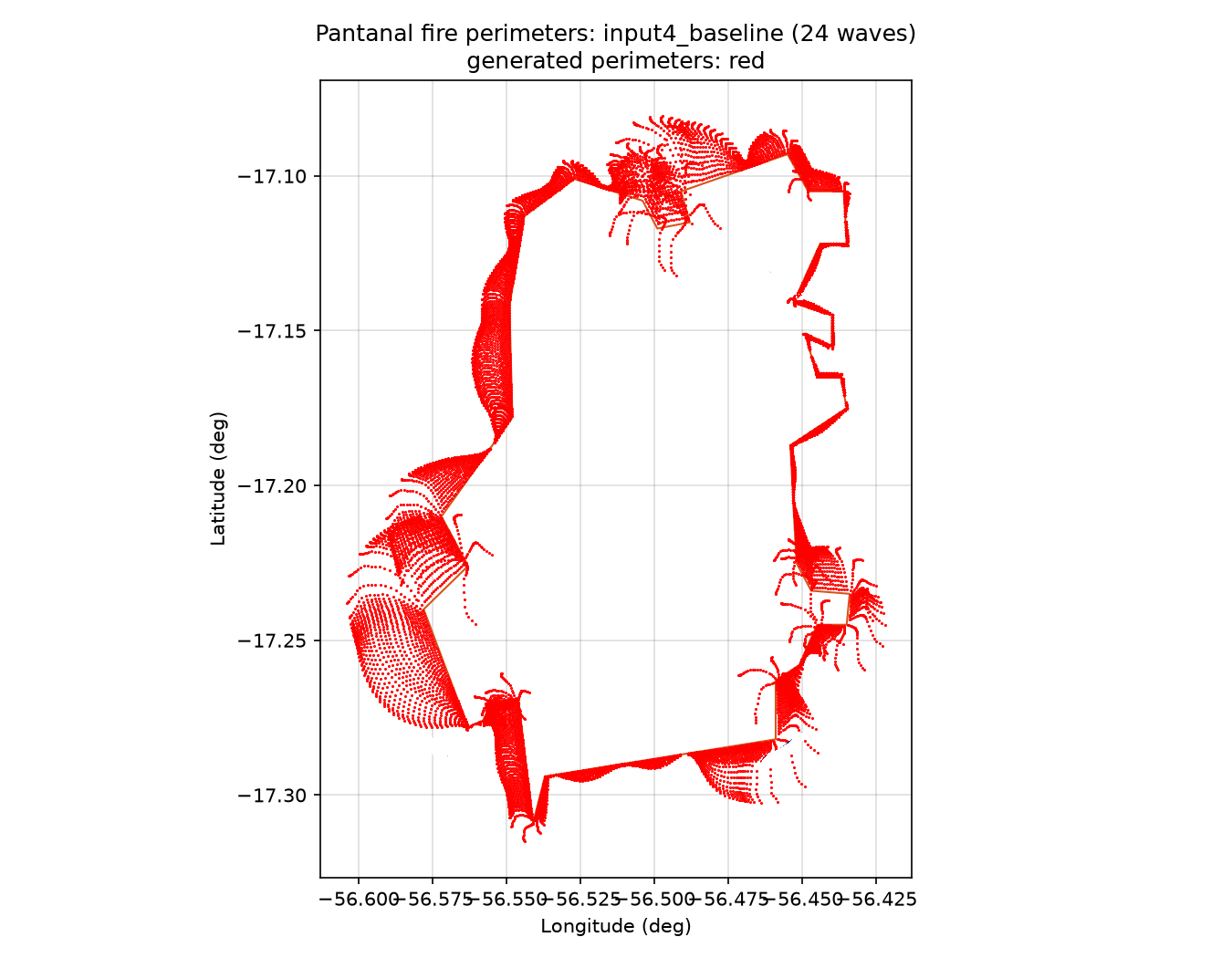}
			\caption{Calculated fire propagation.}
			\label{fig:sep04_simulation}
		\end{subfigure}
		\hfill
		\begin{subfigure}[t]{0.29\textwidth}
			\centering
			\includegraphics[width=\linewidth]{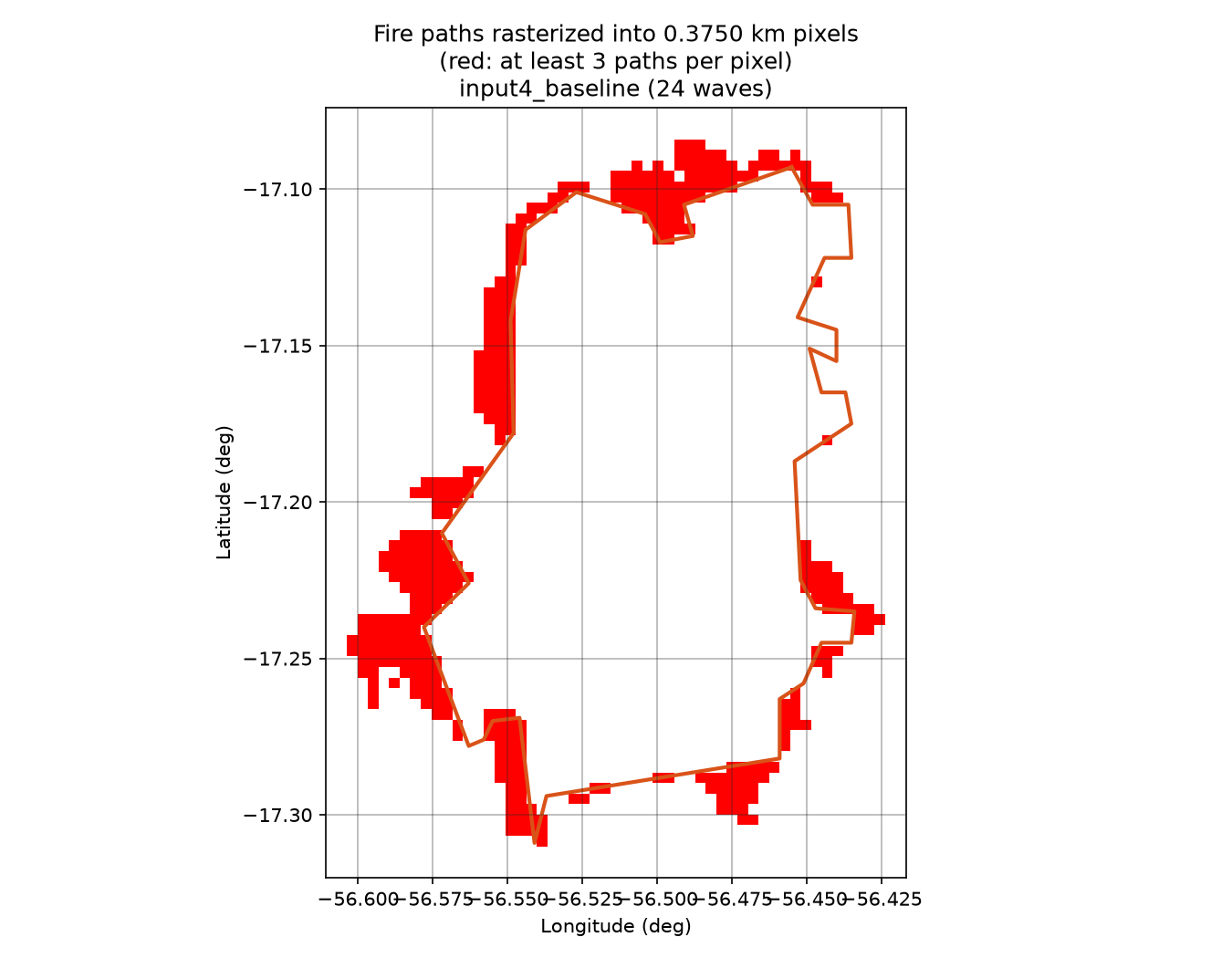}
			\caption{Calculated propagation rasterized into 375-m pixels.}
			\label{fig:sep04_simulation_pix}
		\end{subfigure}
		\hfill
		\begin{subfigure}[t]{0.2\textwidth}
			\centering
			\includegraphics[width=\linewidth]{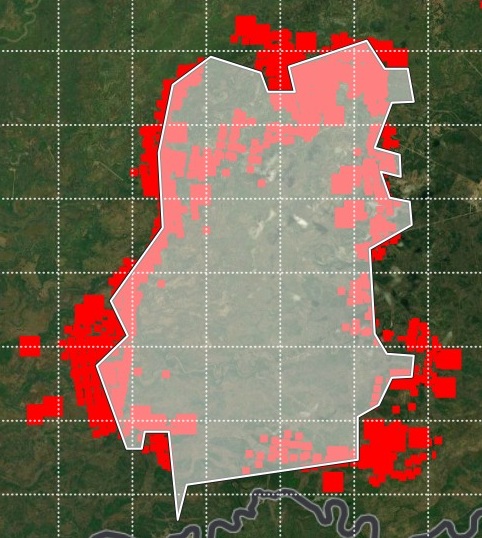}
			\caption{Subsequent satellite observation.}
			\label{fig:sep04_observed}
		\end{subfigure}
		
		\caption{Wildfire propagation experiment for September 4, 2020:
			(a) satellite observation used to construct the updated initial fire
			perimeter; (b) hourly fire propagation calculated by the proposed
			trajectory model; (c) the same calculated trajectories rasterized
			into 375-m pixels, displayed at the native resolution of the VIIRS
			active-fire product; and (d) subsequent satellite observation used
			for comparison.}
		\label{fig:sep04}
	\end{figure}

	For the second experiment, the updated satellite-derived boundary on
	\rev{September 4 was used to reinitialize the model} (Fig.~\ref{fig:sep04}).
	This daily reinitialization prevents accumulated geometric discrepancies
	from one interval from influencing the next. As in the first
	experiment, the calculation is shown both as continuous trajectories and in
	the rasterized, satellite-scale form
	(Figs.~\ref{fig:sep04_simulation}
	and~\ref{fig:sep04_simulation_pix}).

	The calculated front again develops preferentially in particular sectors
	rather than expanding uniformly around the initial perimeter. This behavior
	is also apparent in the subsequent satellite observation, where active-fire
	detections are concentrated in specific portions of the domain. The
	comparison therefore illustrates how representing local differences in
	propagation behavior allows the fitted model to reproduce this nonuniformity.
	The pixel panel of Fig.~\ref{fig:sep04_simulation_pix} highlights the same
	sectors and provides a direct visual link between the calculated
	trajectories and the discrete satellite detections.

	\subsection{Fire Propagation on September 5, 2020}
	\label{subsec:sep05}
	
\begin{figure}[htbp]
		\centering
		
		\begin{subfigure}[t]{0.19\textwidth}
			\centering
			\includegraphics[width=\linewidth]{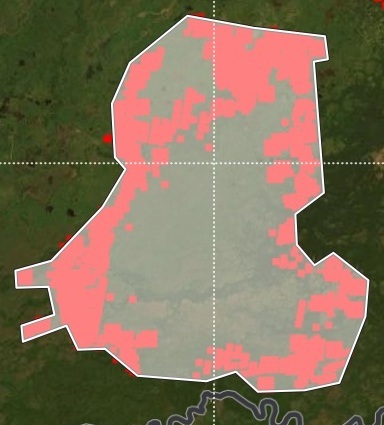}
			\caption{Initial satellite-derived perimeter.}
			\label{fig:sep05_initial}
		\end{subfigure}
		\hfill
		\begin{subfigure}[t]{0.28\textwidth}
			\centering
			\includegraphics[width=\linewidth]{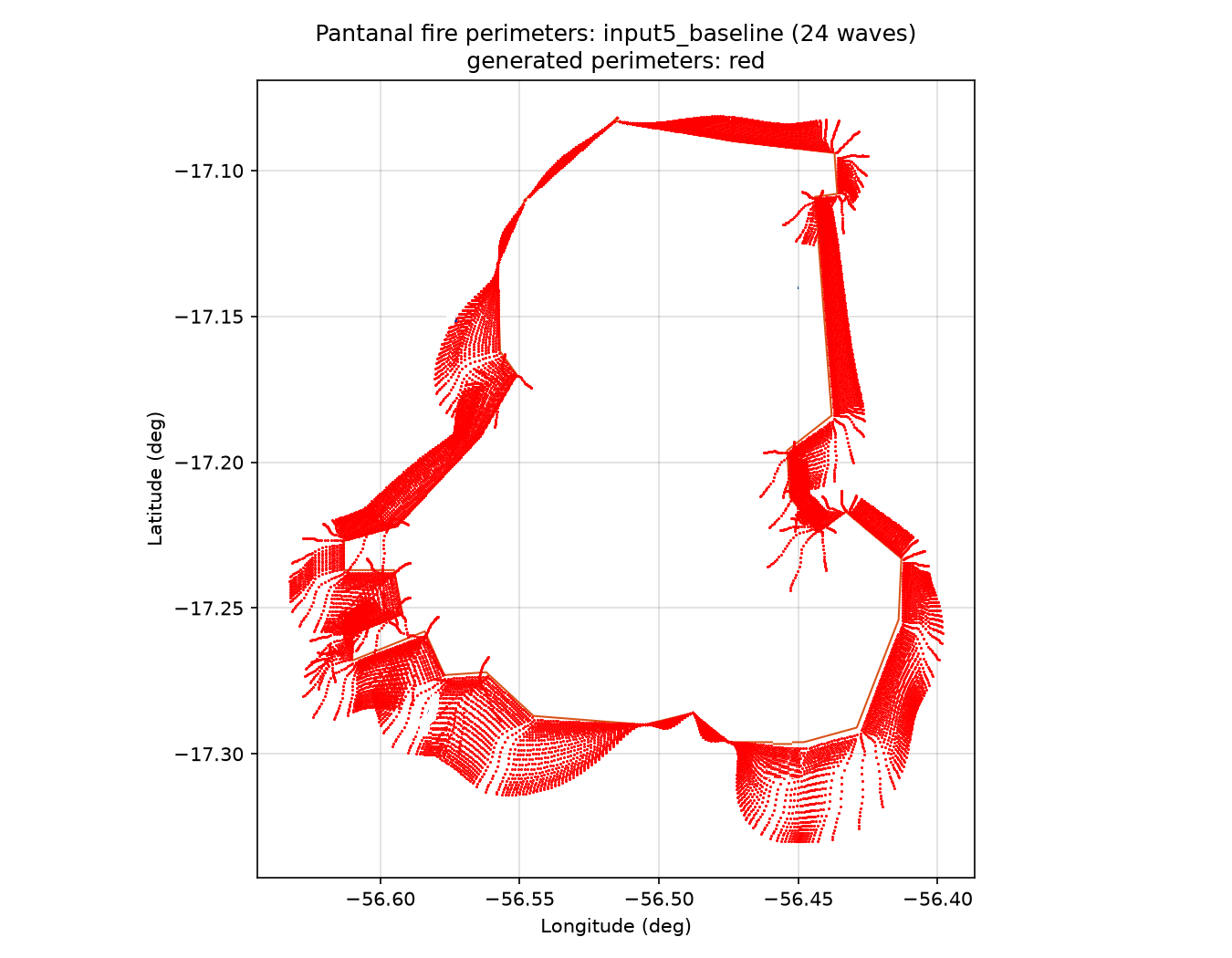}
			\caption{Calculated fire propagation.}
			\label{fig:sep05_simulation}
		\end{subfigure}
		\hfill
		\begin{subfigure}[t]{0.28\textwidth}
			\centering
			\includegraphics[width=\linewidth]{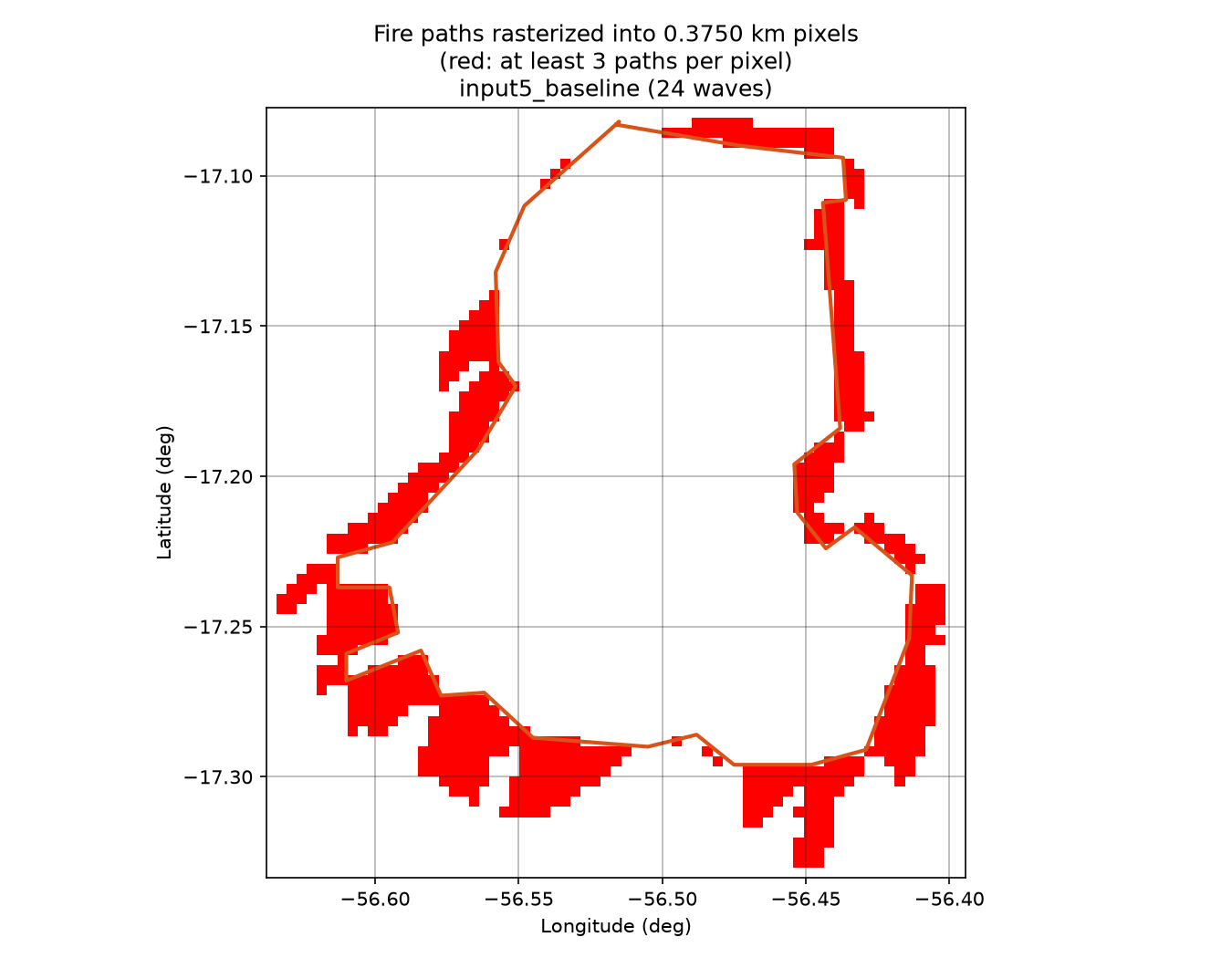}
			\caption{Calculated propagation rasterized into 375-m pixels.}
			\label{fig:sep05_simulation_pix}
		\end{subfigure}
		\hfill
		\begin{subfigure}[t]{0.21\textwidth}
			\centering
			\includegraphics[width=\linewidth]{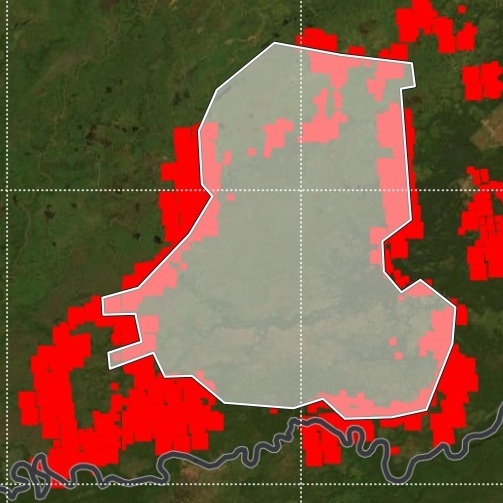}
			\caption{Subsequent satellite observation.}
			\label{fig:sep05_observed}
		\end{subfigure}
		
		\caption{Wildfire propagation experiment for September 5, 2020:
			(a) satellite observation used to construct the initial fire
			perimeter; (b) hourly fire propagation calculated by the proposed
			trajectory model; (c) the same calculated trajectories rasterized
			into 375-m pixels, displayed at the native resolution of the VIIRS
			active-fire product; and (d) subsequent satellite observation used
			for comparison.}
		\label{fig:sep05}
	\end{figure}

	\rev{The third experiment uses a newly initialized satellite-derived perimeter for the September 5 case. The exact acquisition timestamp of the initial satellite observation is unavailable in the retained archive, as discussed in Section~\ref{subsec:data_provenance}.} As shown in Fig.~\ref{fig:sep05_initial}, its
	geometry differs considerably from those of the preceding events,
	providing an additional demonstration of the trajectory formulation for a more
	irregular initial boundary. The calculated propagation is shown in
	Fig.~\ref{fig:sep05_simulation}, and its rasterized, satellite-scale
	rendering is shown in Fig.~\ref{fig:sep05_simulation_pix}.
	
	Despite the more irregular initial geometry, the model again generates
	spatially differentiated propagation, with strong advance in some sectors
	and comparatively small displacement in others. The subsequent satellite
	observation exhibits a similarly heterogeneous distribution of active-fire
	detections, supporting the ability of the trajectory formulation to
	represent nonuniform propagation along irregular fire boundaries. The
	rasterized panel of Fig.~\ref{fig:sep05_simulation_pix} confirms this
	pattern at the observational scale, with marked pixels concentrated in the
	same sectors in which the subsequent detections cluster.

	\subsection{Coordinate-Referenced Sectoral Comparison}
	\label{subsec:sector_comparison}
	
	The panels in Figs.~\ref{fig:sep03}--\ref{fig:sep05} are georeferenced by
	latitude and longitude, but the VIIRS detections are discrete active-fire
	pixels rather than closed, simultaneous perimeter polygons. Consequently, overlap indices such as intersection over union or pointwise
	Hausdorff distance cannot be obtained reliably from the displayed panels
	alone without an explicit image-registration and rasterization procedure.
	We therefore first use a directional sector comparison, which is consistent
	with the spatial resolution of the available observations; the registered
	pixel-level analysis is presented separately in
	Section~\ref{subsec:pixel_metrics}. Strong, moderate, and
	limited development refer to the relative outward extent within each daily
	panel, not to universal spread-rate thresholds. The rasterized model panels
	(Figs.~\ref{fig:sep03_simulation_pix}--\ref{fig:sep05_simulation_pix})
	present the results at the observational pixel scale, while
	Table~\ref{tab:sector_comparison} summarizes the qualitative sector-level
	comparison.
	
	\begin{table}[htbp]
		\centering
		\small
		\caption{Sectoral comparison of the calculated trajectories and subsequent VIIRS active-fire pattern. The daily $\gamma(x,y)$ field is calibrated using the propagation over the first 12 model hours and then held fixed during hours 13--24; the comparison therefore combines calibration-period agreement with {a fixed-parameter continuation (not independent validation) during the second 12~h}.}
		\label{tab:sector_comparison}
		\begin{tabular}{p{0.10\textwidth}p{0.31\textwidth}p{0.31\textwidth}p{0.17\textwidth}}
			\toprule
			Date & Subsequent VIIRS pattern & Calculated trajectory pattern & Spatial agreement \\
			\midrule
			September 3 & Extensive development beyond the initial boundary, particularly toward northern and lateral sectors, with spatial gaps & Pronounced northern and western advancement, with localized development elsewhere and limited motion in fitted low-$\gamma$ sectors & Broad directional agreement; local differences remain \\
			September 4 & Concentrated activity outside the updated boundary in several western, northern, eastern, and southern sectors & Strongest advancement in western and northern portions, with localized propagation toward the southeast & Partial-to-broad agreement in the dominant nonuniform pattern \\
			September 5 & Irregular outward activity concentrated mainly around western--southwestern and northern--northeastern portions & Pronounced western--southwestern advancement together with localized northern and eastern trajectories & Broad sectoral agreement; pixel-scale correspondence is not claimed \\
			\bottomrule
		\end{tabular}
	\end{table}
	
	\subsection{Pixel-Level Quantitative Agreement with the Observed Pattern}
	\label{subsec:pixel_metrics}
	
	The sector-level comparison in Table~\ref{tab:sector_comparison} provides a
	qualitative assessment of the propagation pattern, whereas the pixel
	representation enables a complementary quantitative comparison. The
	calculated and observed patterns were mapped onto the same $375$~m grid.
	For each day, the observed active-fire pattern was obtained by registering
	the initial-perimeter overlay of the satellite comparison map to the
	\texttt{Vertices} polygon of the corresponding input workbook. An isotropic
	iterative-closest-point procedure was used for this registration, resulting
	in errors of $38$--$80$~m, well below the $375$~m pixel size. The VIIRS
	detection pixels were then mapped onto the model grid to define the observed
	set $O$, while the subset $O_{\mathrm{new}}$ contains detections outside the
	initial perimeter. On the model side, $\mathbb{M}$ denotes the marked pixels
	of the complete 24-h simulation, and $\mathbb{M}_{\mathrm{late}}$ denotes
	those associated with waves 13--24, generated after the calibration period
	with the spatial field $\gamma(x,y)$ held fixed.
	
	\rev{Here the hit rate is $|\mathbb M\cap O|/|O|$, the precision is $|\mathbb M\cap O|/|\mathbb M|$, and the false-alarm ratio is $1-\text{precision}$; in late-window rows, $\mathbb M$ and $O$ are replaced by $\mathbb M_{\mathrm{late}}$ and $O_{\mathrm{new}}$.} Table~\ref{tab:pixel_metrics} reports the intersection-over-union (IoU), Dice
	coefficient, hit rate, and false-alarm ratio for the complete simulations and
	the corresponding late-window comparison. The full-run Dice coefficients are
	$0.58$, $0.57$, and $0.47$ for September~3, 4, and 5, respectively, while the
	model reproduces approximately $40$--$58\%$ of the satellite-detected
	active-fire pixels. These results indicate a meaningful spatial
	correspondence under a relatively strict comparison criterion: the
	three-trajectory retention rule defined above is applied, no neighborhood
	tolerance is introduced, and the comparison is performed at the native VIIRS
	resolution.
	
	The numerical values should also be interpreted in light of the observational
	limitations. VIIRS does not continuously observe the wildfire and provides
	active-fire detections only at the times of satellite overpasses. Consequently,
	fire propagation occurring between successive observations is not directly
	recorded. Detection may also be affected by cloud or smoke obscuration,
	vegetation canopy, sensor resolution, and the timing of fire activity relative
	to an overpass. Moreover, active-fire pixels identify locations of detected
	thermal activity rather than a continuous and precise fire perimeter.
	Therefore, an exact pixel-by-pixel reconstruction of the actual propagation
	throughout the 24-h interval cannot be established from these observations
	alone. The false-alarm ratios of $0.24$--$0.44$ should consequently not be
	interpreted exclusively as erroneous simulated propagation, since some
	calculated locations may correspond to fire development that was not captured
	during the available satellite observations.
	
	{The late-window results describe the continuation of the fitted model, but cannot independently establish forecasting skill.} Without further adjustment of
	$\gamma(x,y)$ after the 12-h model calibration window, the late-window simulations overlap
	approximately $23$--$28\%$ of the observed new development, with precisions
	\rev{of approximately $0.44$--$0.69$ against $O_{\mathrm{new}}$}. Considering the discrete nature
	of the satellite observations and the fact that $O_{\mathrm{new}}$ may also
	contain residual detections within previously affected areas, {these results indicate descriptive overlap between the modeled continuation and the available thermal detections, without establishing the accuracy of a predicted fire perimeter.} The pixel-level metrics should
	therefore be viewed as a quantitative complement to the trajectory and
	sector-wise comparisons rather than as an exact measure of the continuous
	physical fire-front evolution.
	
	\begin{table}[htbp]
		\centering
		\small
		\caption{Quantitative pixel-level agreement between the calculated
			propagation and the observed active-fire pattern at the $375$~m
			grid. Top: complete 24-h calculation ($\mathbb{M}$) against all observed
			detections ($O$). Bottom: {fixed-parameter late-window}
			increment ($\mathbb{M}_{\mathrm{late}}$, waves 13--24) against the
			observed new development outside the initial perimeter
			($O_{\mathrm{new}}$). Registration error is $38$--$80$~m.}
		\label{tab:pixel_metrics}
		\begin{tabular}{llrrrrrr}
			\toprule
			Day & Set & $|\mathbb{M}|$ & $|O|$ & IoU & Dice & Hit rate & False-alarm \\
			\midrule
			Sep 3 & full run & 491 & 806 & 0.407 & 0.578 & 0.465 & 0.236 \\
			Sep 4 & full run & 505 & 493 & 0.400 & 0.571 & 0.578 & 0.436 \\
			Sep 5 & full run & 805 & 1144 & 0.305 & 0.467 & 0.398 & 0.435 \\
			\midrule
			Sep 3 & late window & 286 & 789 & 0.223 & 0.365 & 0.248 & \rev{0.315} \\
			Sep 4 & late window & 301 & 484 & 0.204 & 0.339 & 0.275 & \rev{0.558} \\
			Sep 5 & late window & 507 & 1124 & 0.184 & 0.310 & 0.225 & \rev{0.501} \\
			\bottomrule
		\end{tabular}
	\end{table}

	Overall, the three experiments {illustrate} that the proposed formulation can
	represent strongly nonuniform {modeled} fire propagation across different irregular
	initial perimeters. Although the available satellite-observed fire pattern is used to calibrate
	$\gamma(x,y)$ over the first 12 model hours, the coefficient field is
	subsequently held fixed,
	and hours 13--24 are simulated without further adjustment. The resulting
	trajectory, pixel-level, and sector-wise comparisons indicate that relevant
	spatial features of the subsequent fire development are retained beyond {the fitting steps; because the observed comparison pattern also informed calibration, this should not be interpreted as independent prediction.} Remaining discrepancies can be attributed, at least in
	part, to the discrete nature of the satellite observations, unresolved local
	fuel conditions, prescribed moisture, environmental-data resolution, and
	processes not explicitly represented by the model.

	
	\subsection{One-at-a-Time Sensitivity of the Calculated Propagation}
	\label{subsec:pixel_sensitivity}
	
	A one-at-a-time sensitivity analysis was performed for the five principal
	model and forcing factors: $\alpha$, $\beta$, $\gamma$, wind speed $U$,
	and wind direction $\Theta$. For each daily experiment, one factor was
	halved or doubled relative to its reference value while all other inputs
	were kept unchanged. Together with the reference cases, this resulted in
	45 factor--level--day entries when the reference configuration is counted separately for each factor. Of these, 36 entries concern the four interpretable magnitude-related factors, including 12 repeated baseline entries; if the baseline is reused within each day, they correspond to 27 distinct configurations. A one-at-a-time sensitivity analysis was performed for the
	four principal magnitude-related model and forcing factors:
	$\alpha$, $\beta$, $\gamma$, and wind speed $U$.
	For each daily experiment, one factor was halved or doubled
	relative to its reference value while all other inputs were
	kept unchanged. Reusing the baseline configuration within
	each day yields 27 distinct simulation runs across the three
	days. Directional wind sensitivity, which requires additive
	angular rotations of the form
	$\Theta_{\delta} = \operatorname{mod}(\Theta + \delta, 360^\circ)$
	(e.g., $\delta = \pm 15^\circ, \pm 30^\circ$)
	rather than multiplicative angle scaling, is left for future
	work. Table~5 summarizes the parameter variations considered.
	
	\begin{table}[!htbp]
		\centering
		\small
		\caption{Parameter variations considered in the sensitivity analysis.
			For each experiment, one quantity was varied while all other model
			inputs were kept at their reference values.}
		\label{tab:sensitivity}
		
		\begin{tabularx}{\textwidth}{
				>{\centering\arraybackslash}p{0.10\textwidth}
				>{\centering\arraybackslash}p{0.18\textwidth}
				>{\centering\arraybackslash}p{0.22\textwidth}
				X}
			\toprule
			\textbf{Quantity} &
			\textbf{Reference case} &
			\textbf{Values considered} &
			\textbf{Expected model response} \\
			\midrule
			
			$\alpha$
			& $1.00$
			& $0.5\alpha,\;\alpha,\;2\alpha$
			& Modifies the dependence of $R_H$ on normalized wind speed through
			$(1+U^*)^\alpha$. \\
			
			\addlinespace
			
			$\beta$
			& $0.08$
			& $0.5\beta,\;\beta,\;2\beta$
			& Increasing $\beta$ strengthens the exponential reduction of $R_H$
			associated with fuel moisture. \\
			
			\addlinespace
			
			$\gamma$
			& Calibrated spatial field
			& $0.5\gamma,\;\gamma,\;2\gamma$
			& Scales $R_H$ and $R_B$ proportionally while preserving the spatial
			structure of the calibrated field. \\
			
			\addlinespace
			
			$U$
			& Observed ERA5 wind speed
			& $0.5U,\;U,\;2U$
			& Changes the magnitude of the wind forcing and the resulting rate of
			spread. \\
			
%
			
			\bottomrule
		\end{tabularx}
	\end{table}
	\FloatBarrier
	
	

	\begin{figure}[p]
		\centering
		\captionsetup{font=small}
		\noindent \rev{\textbf{$\alpha$}}\par
		\includegraphics[width=0.30\textwidth]{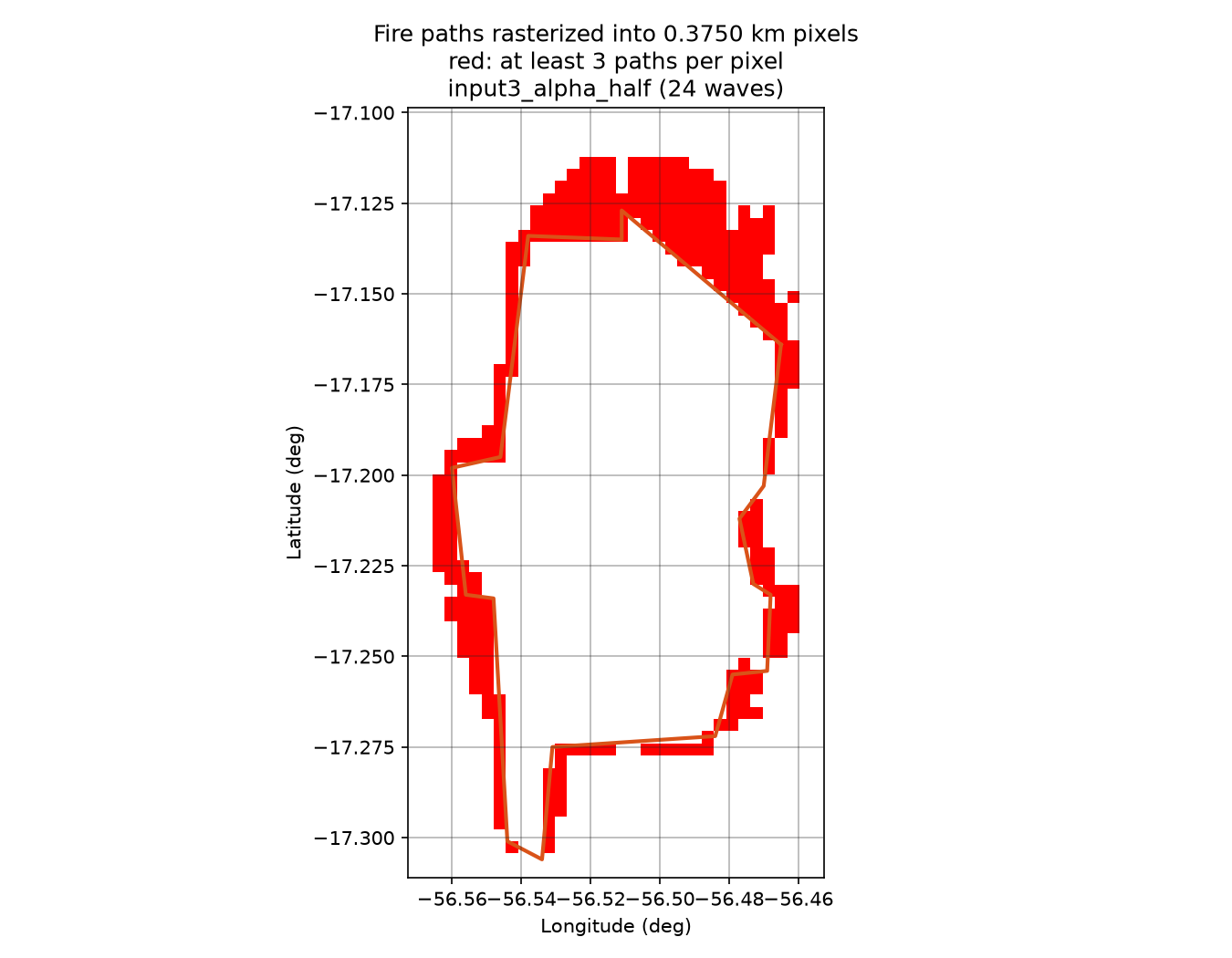}\hfill
		\includegraphics[width=0.30\textwidth]{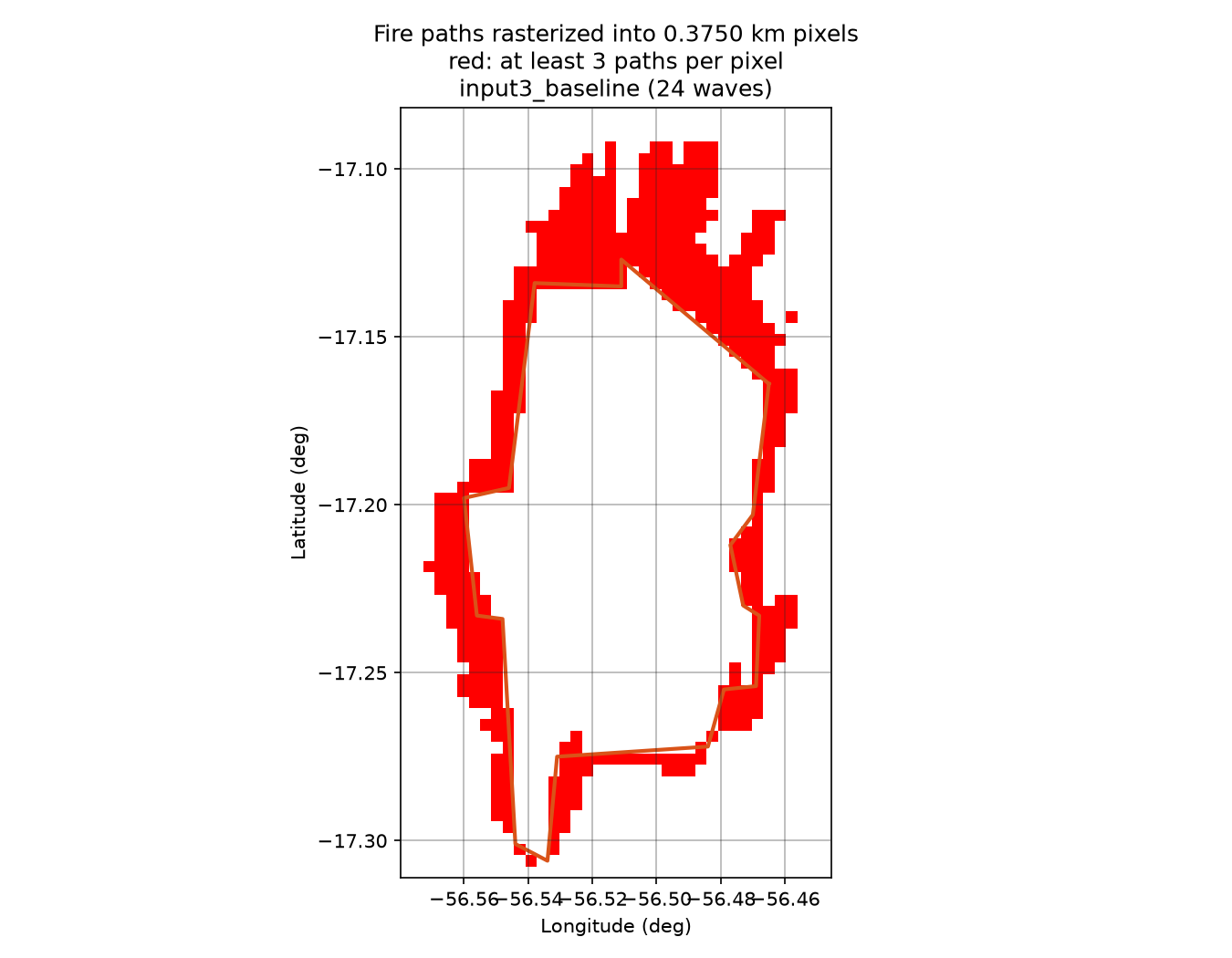}\hfill
		\includegraphics[width=0.30\textwidth]{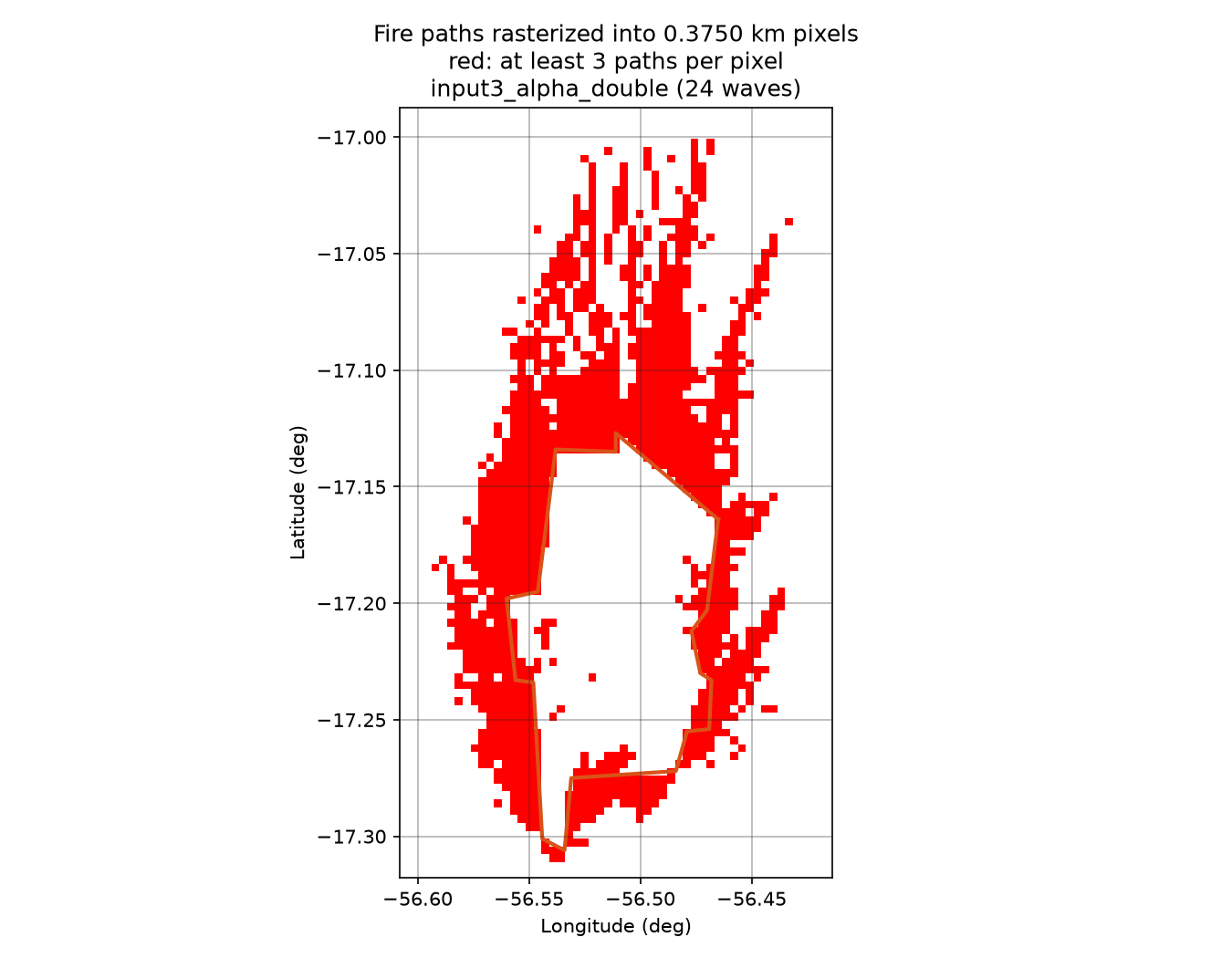}\par
		\noindent \rev{\textbf{$\beta$}}\par
		\includegraphics[width=0.30\textwidth]{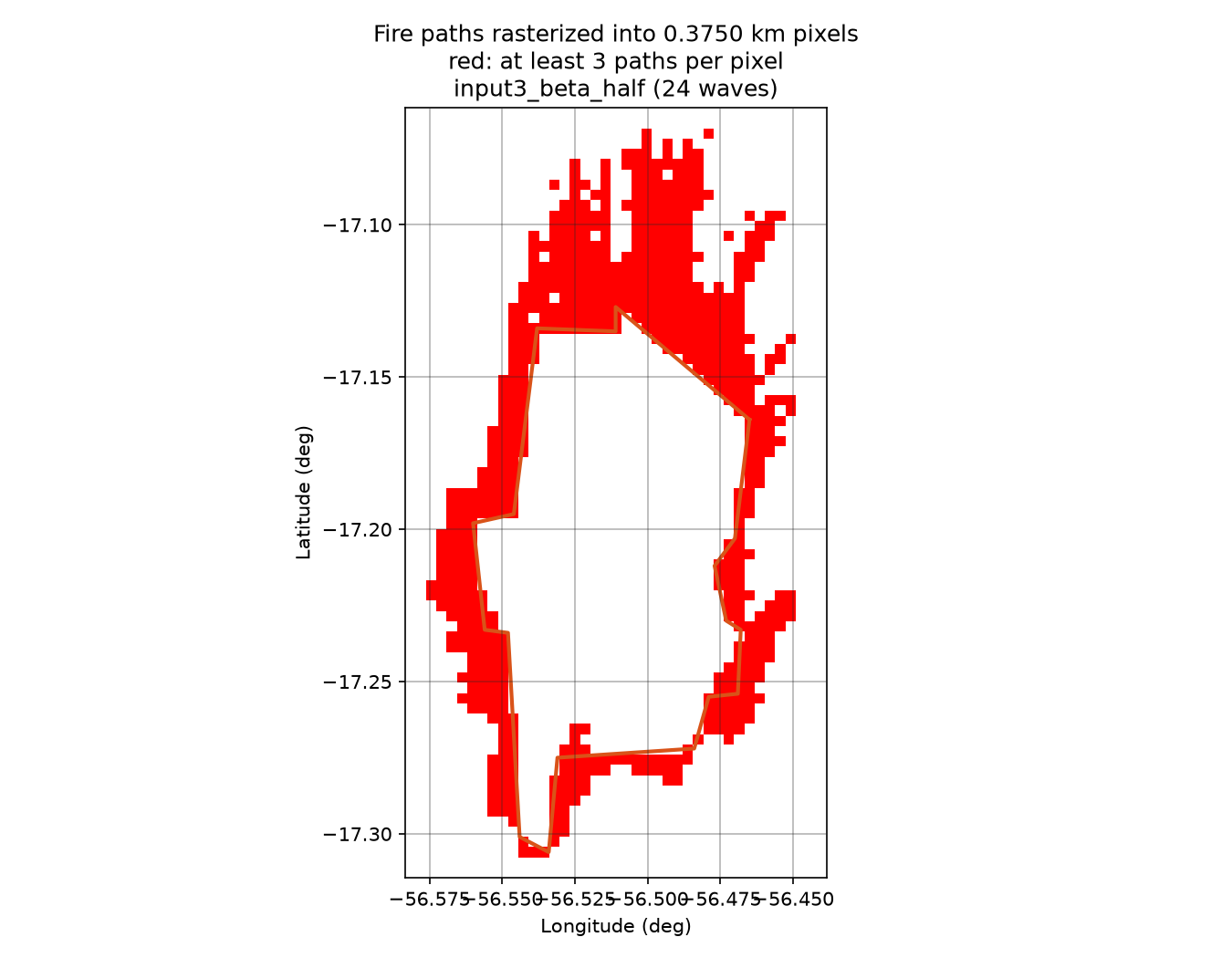}\hfill
		\includegraphics[width=0.30\textwidth]{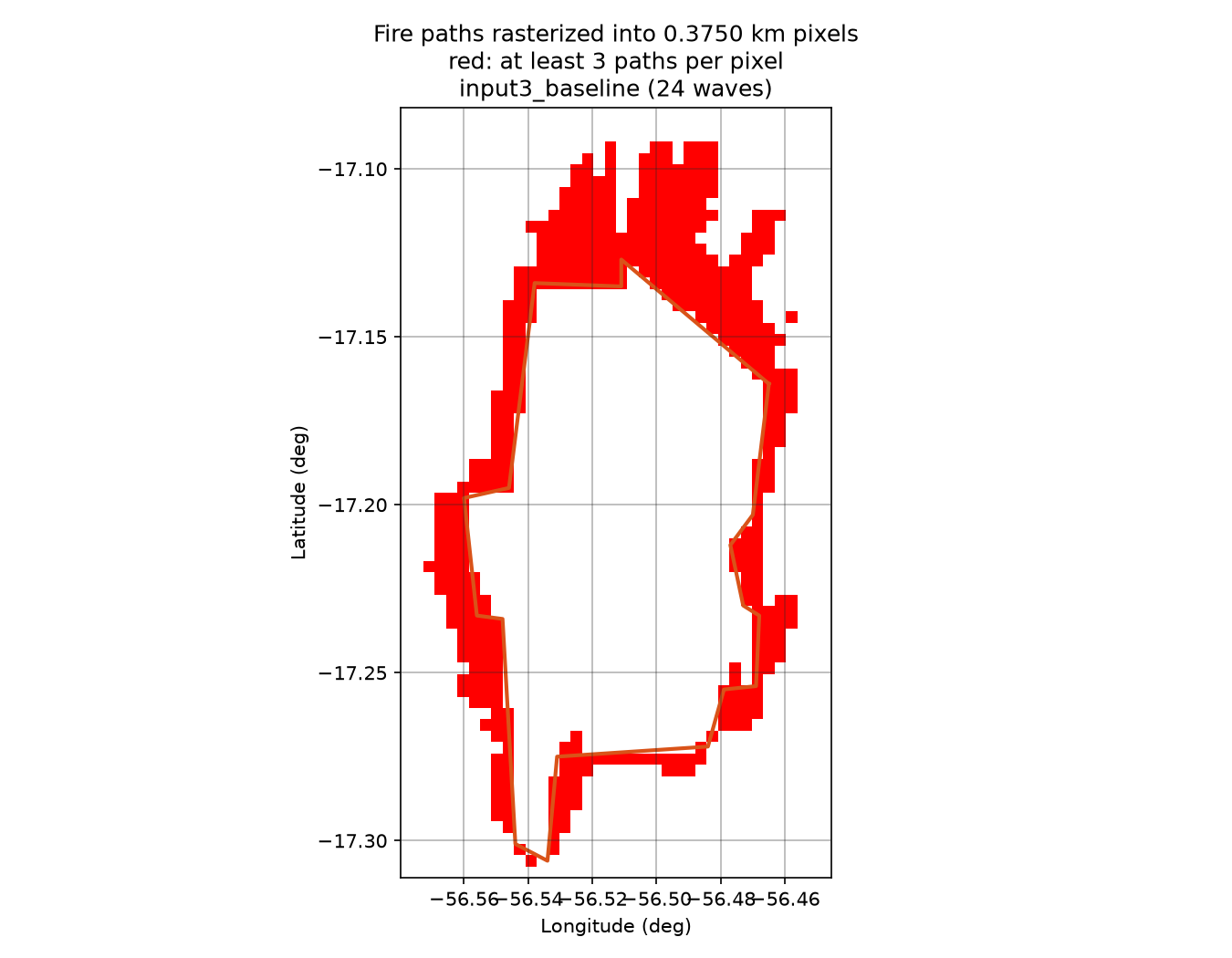}\hfill
		\includegraphics[width=0.30\textwidth]{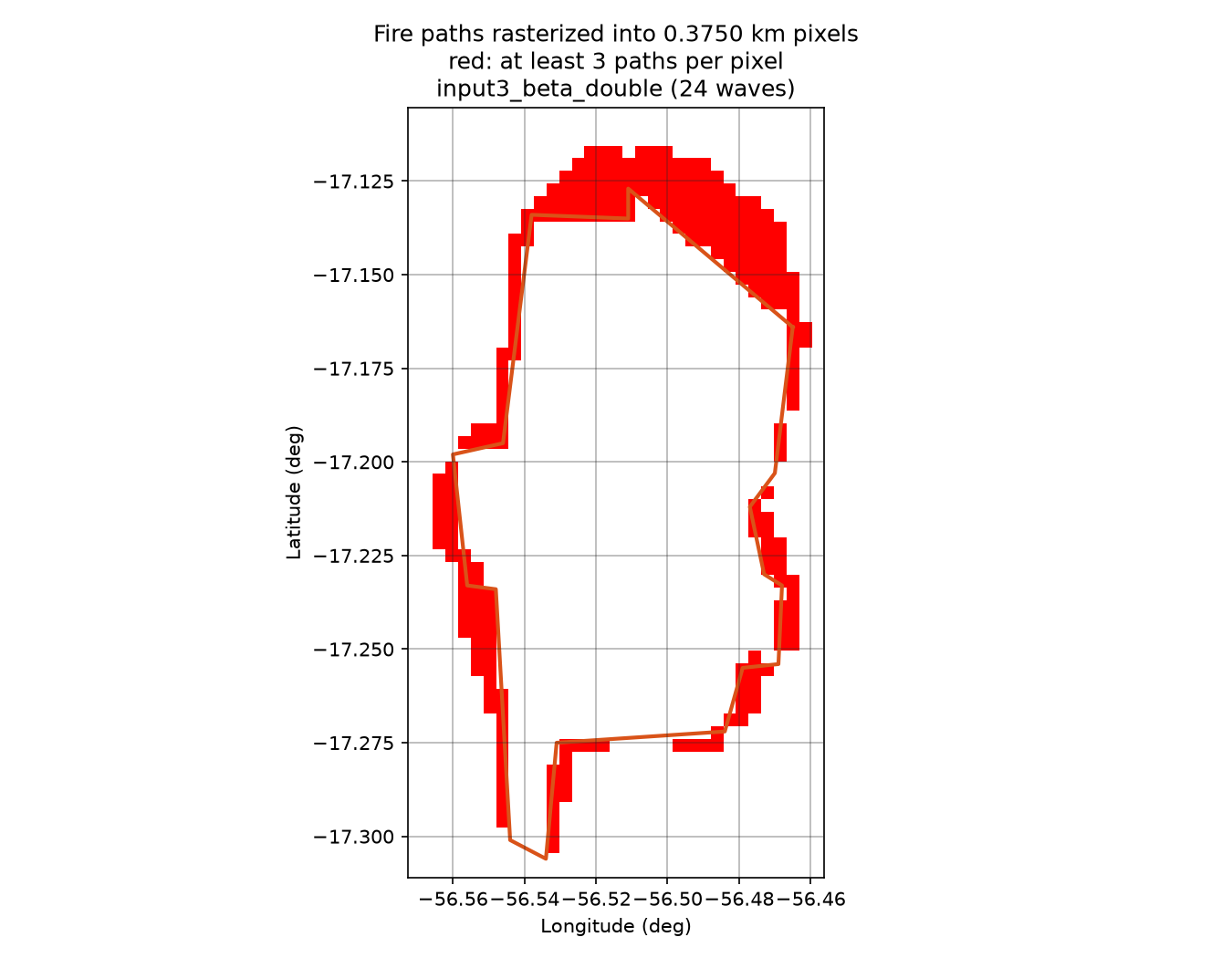}\par
		\noindent \rev{\textbf{$\gamma$}}\par
		\includegraphics[width=0.30\textwidth]{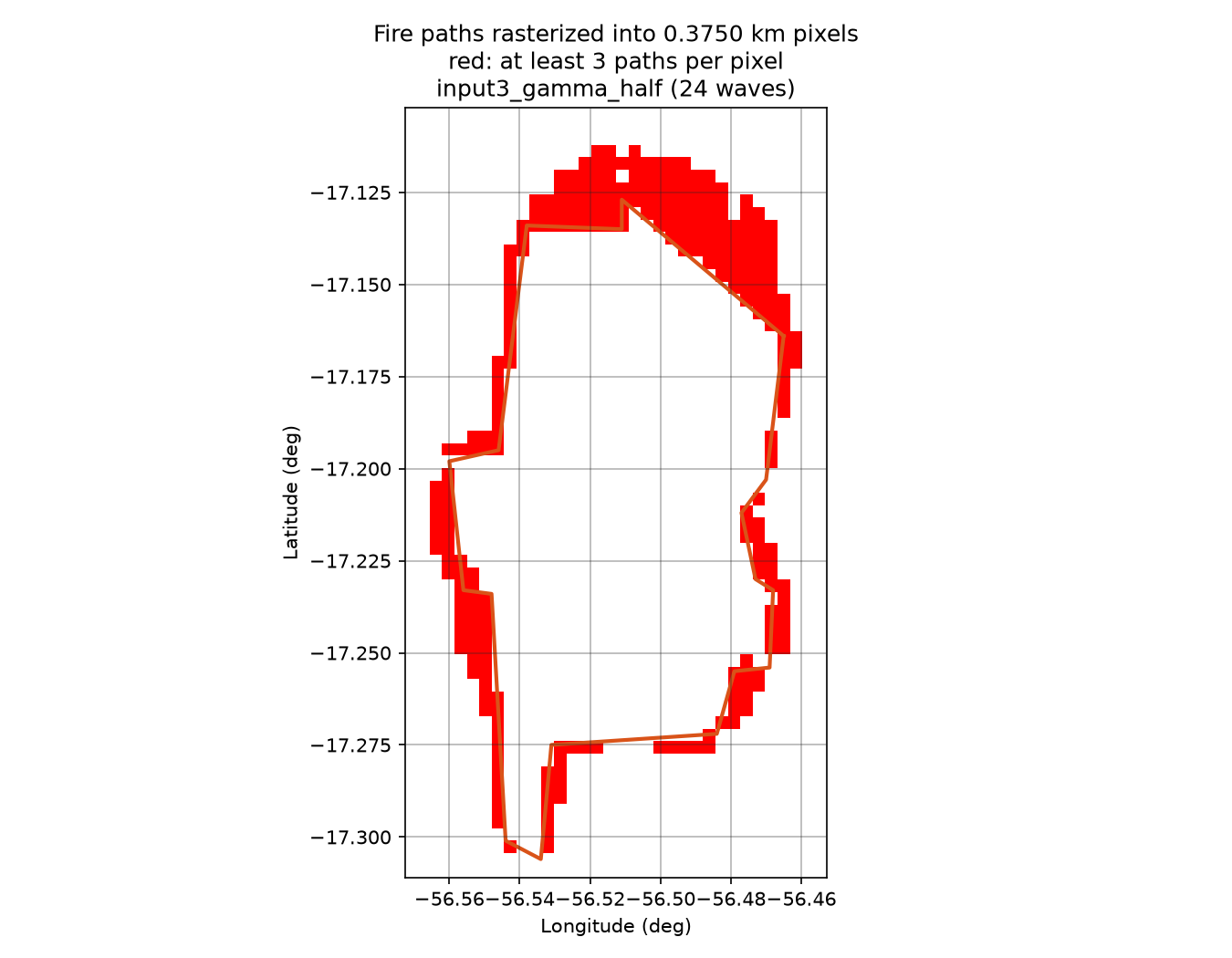}\hfill
		\includegraphics[width=0.30\textwidth]{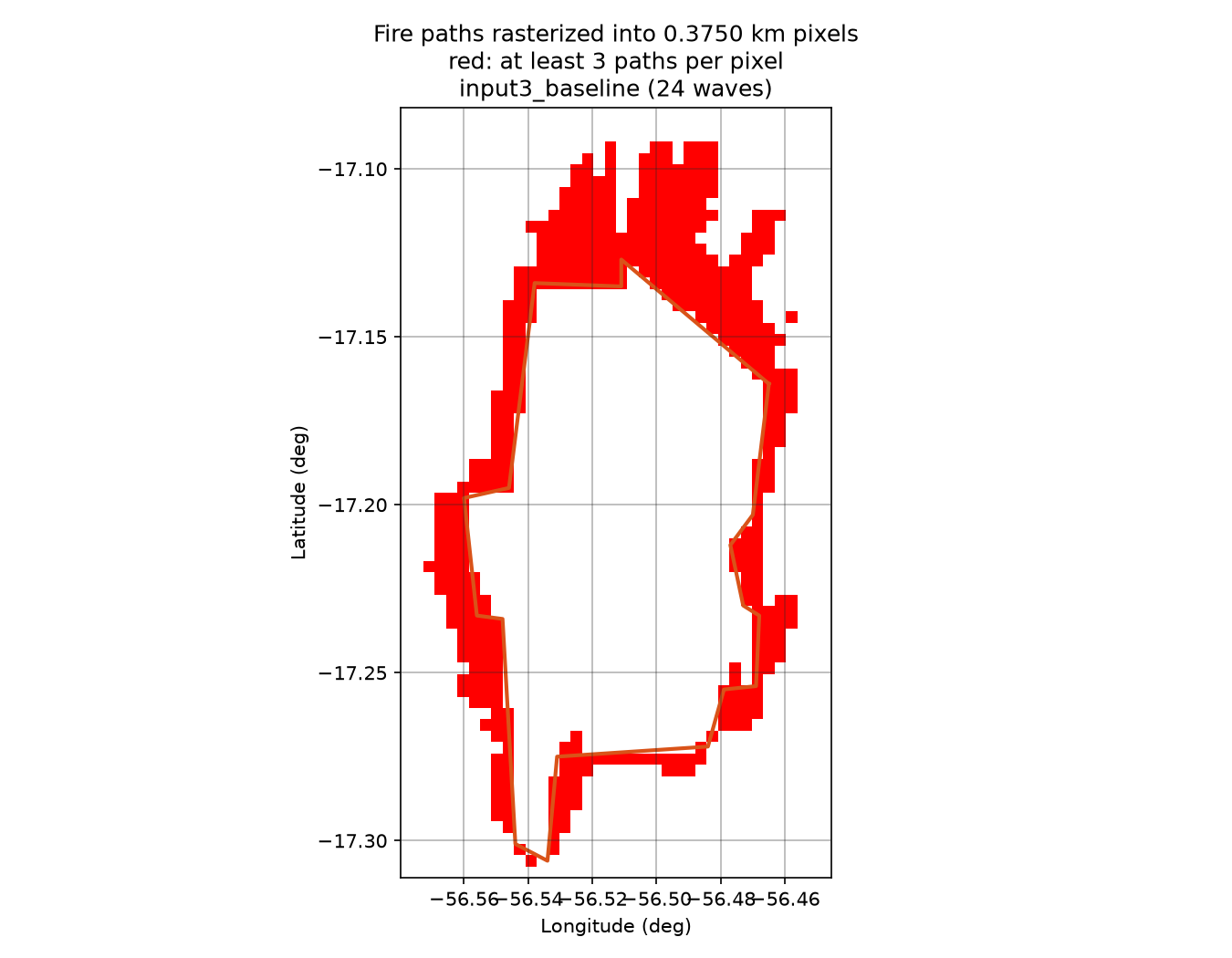}\hfill
		\includegraphics[width=0.30\textwidth]{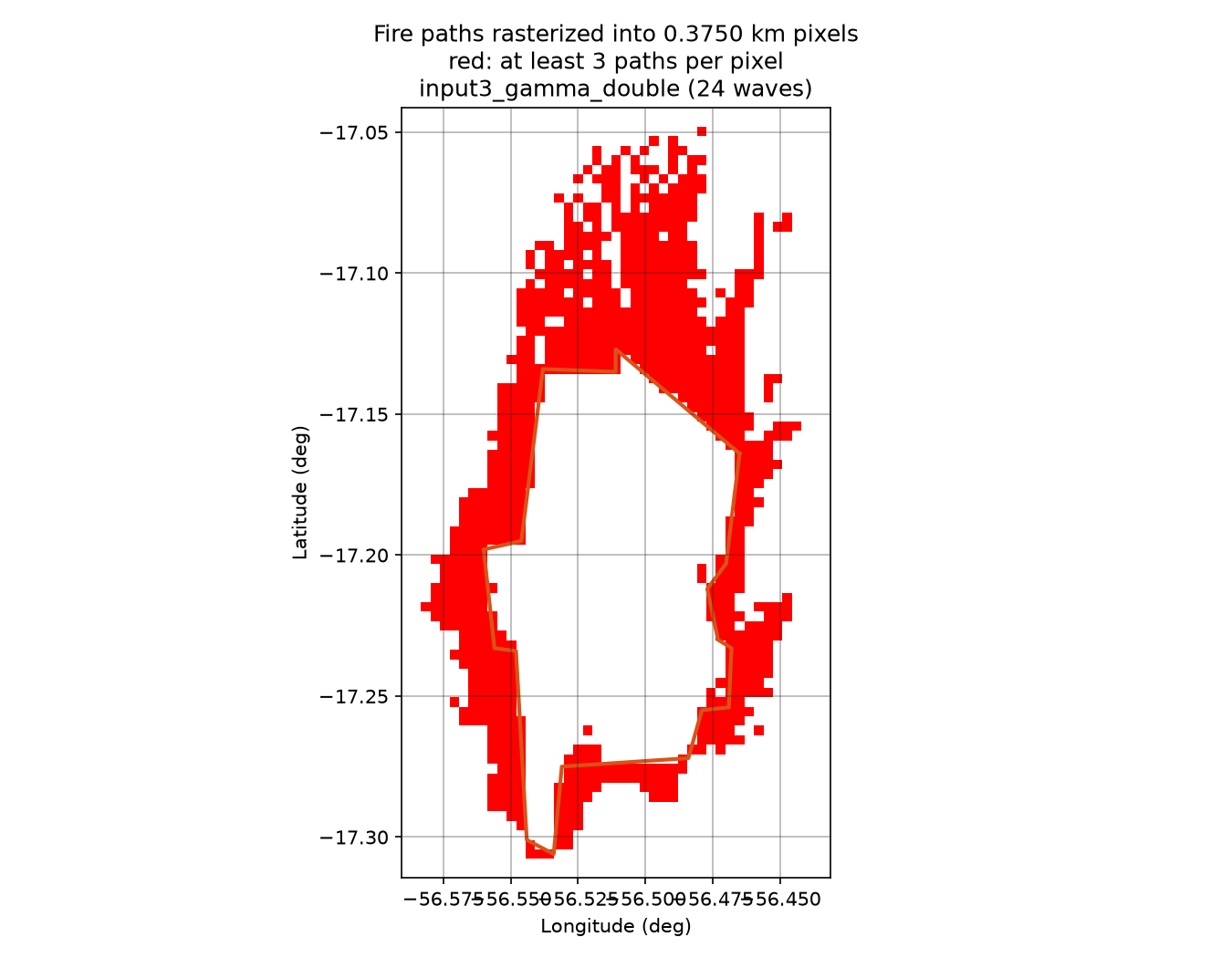}\par
		\noindent \rev{\textbf{Wind speed $U$}}\par
		\includegraphics[width=0.30\textwidth]{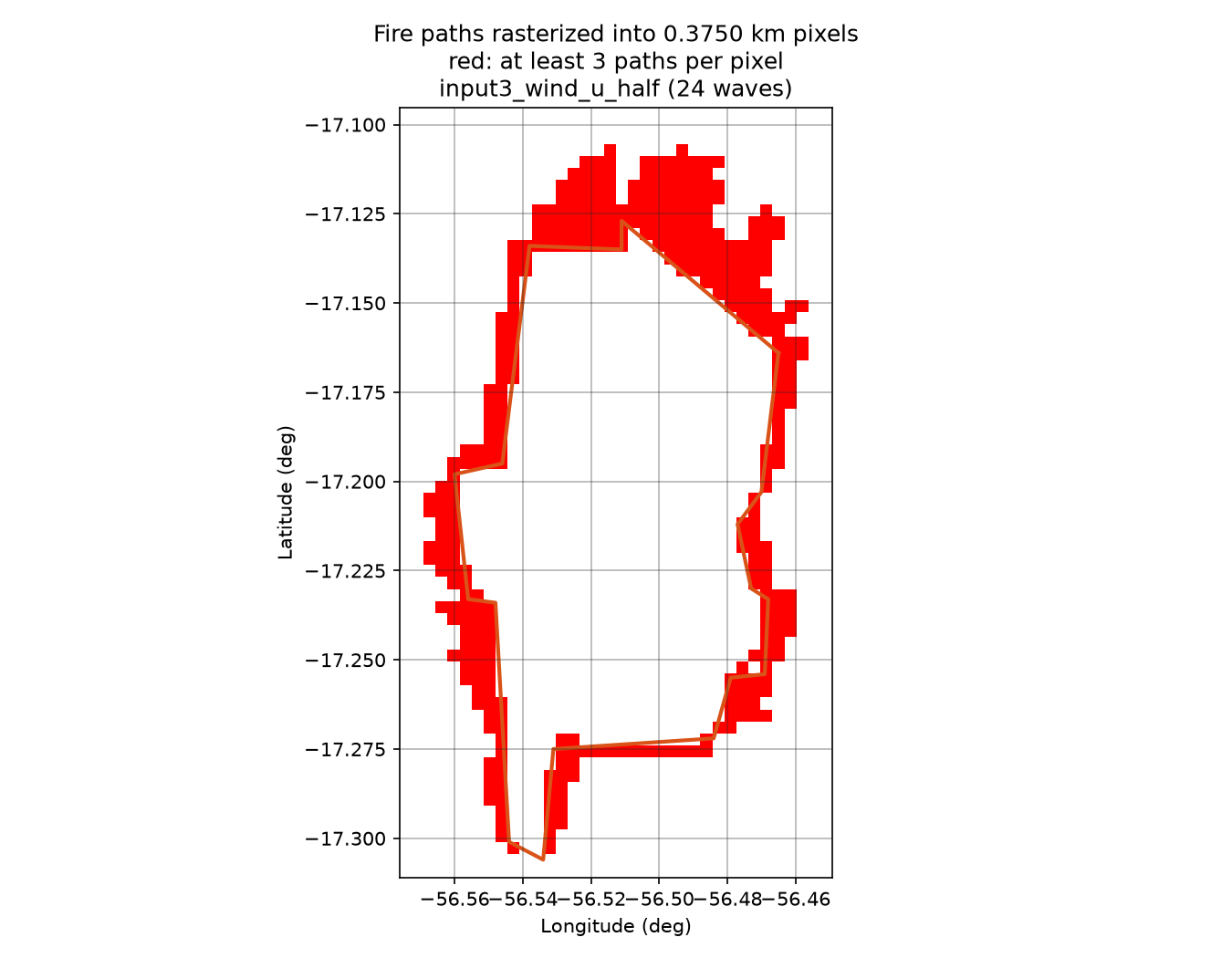}\hfill
		\includegraphics[width=0.30\textwidth]{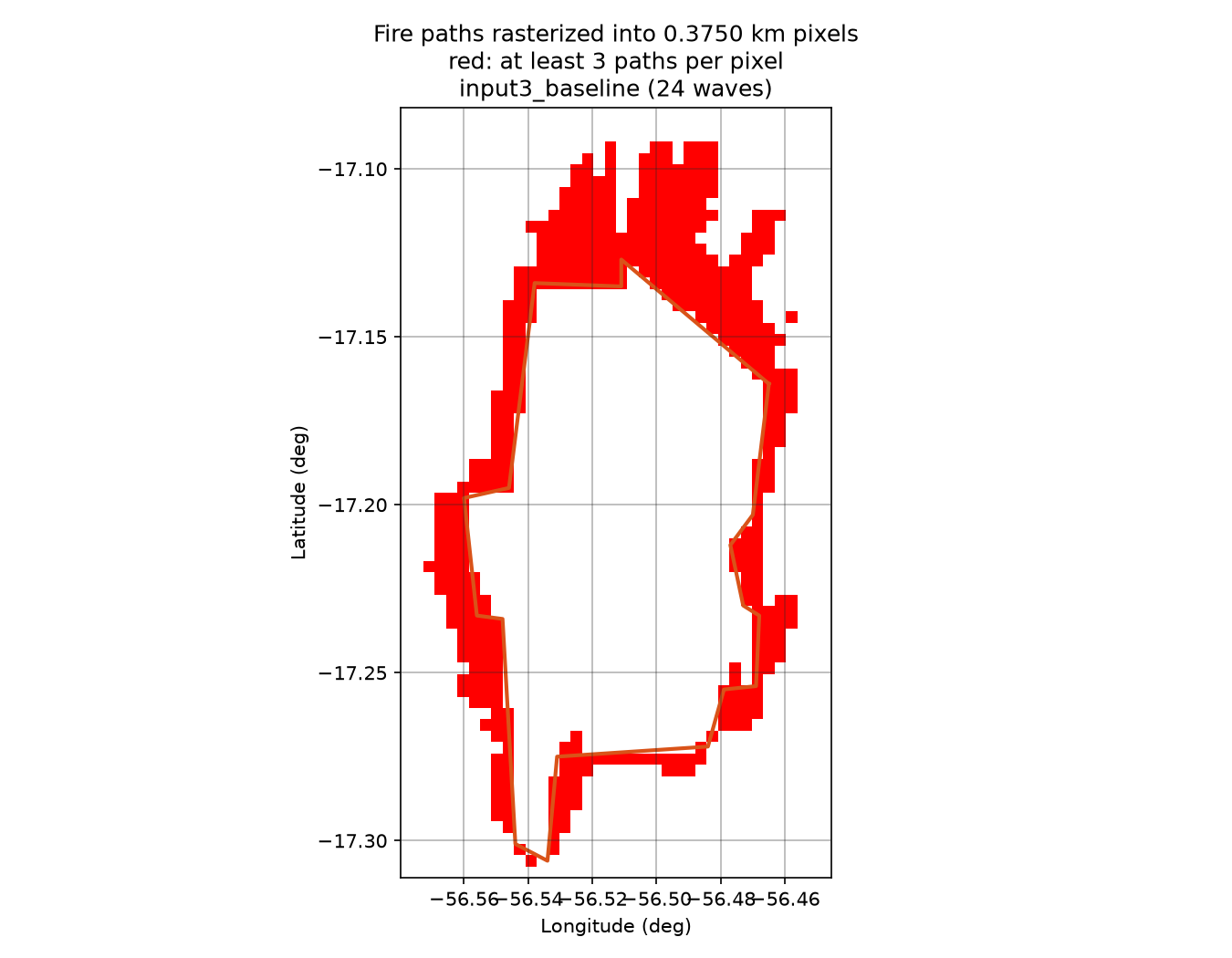}\hfill
		\includegraphics[width=0.30\textwidth]{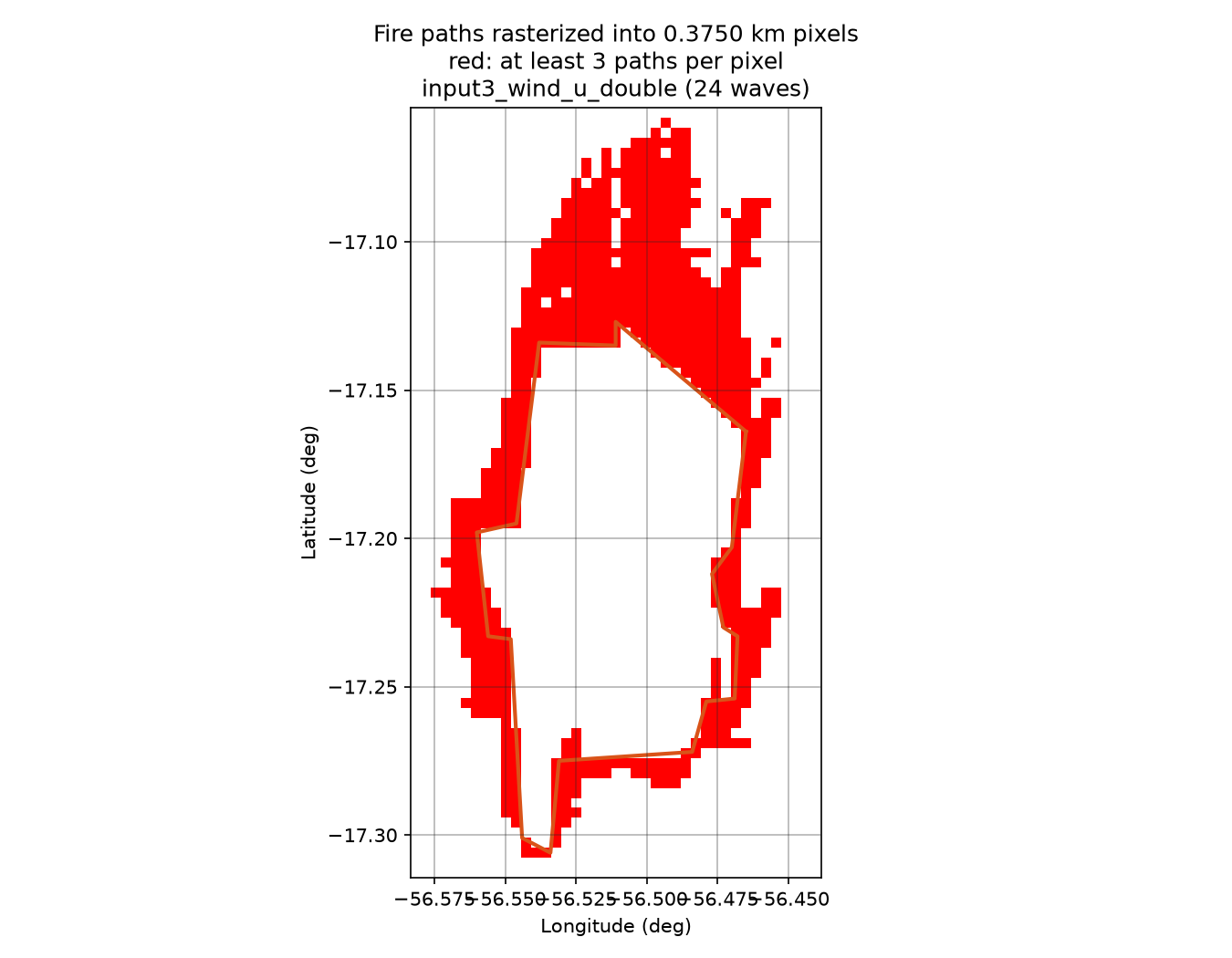}\par
		\caption{\rev{One-at-a-time sensitivity of calculated propagation for September 3, 2020. Rows show $\alpha$, $\beta$, $\gamma$, and wind speed $U$; columns show half, reference, and double values. The plotted numerical outputs are unchanged. Red pixels denote cells crossed by at least three calculated trajectories; the orange curve denotes the initial perimeter.}}
		\label{fig:pixel_sens}
	\end{figure}
	
	\begin{figure}[p]
		\centering
		\captionsetup{font=small}
		\noindent \rev{\textbf{$\alpha$}}\par
		\includegraphics[width=0.30\textwidth]{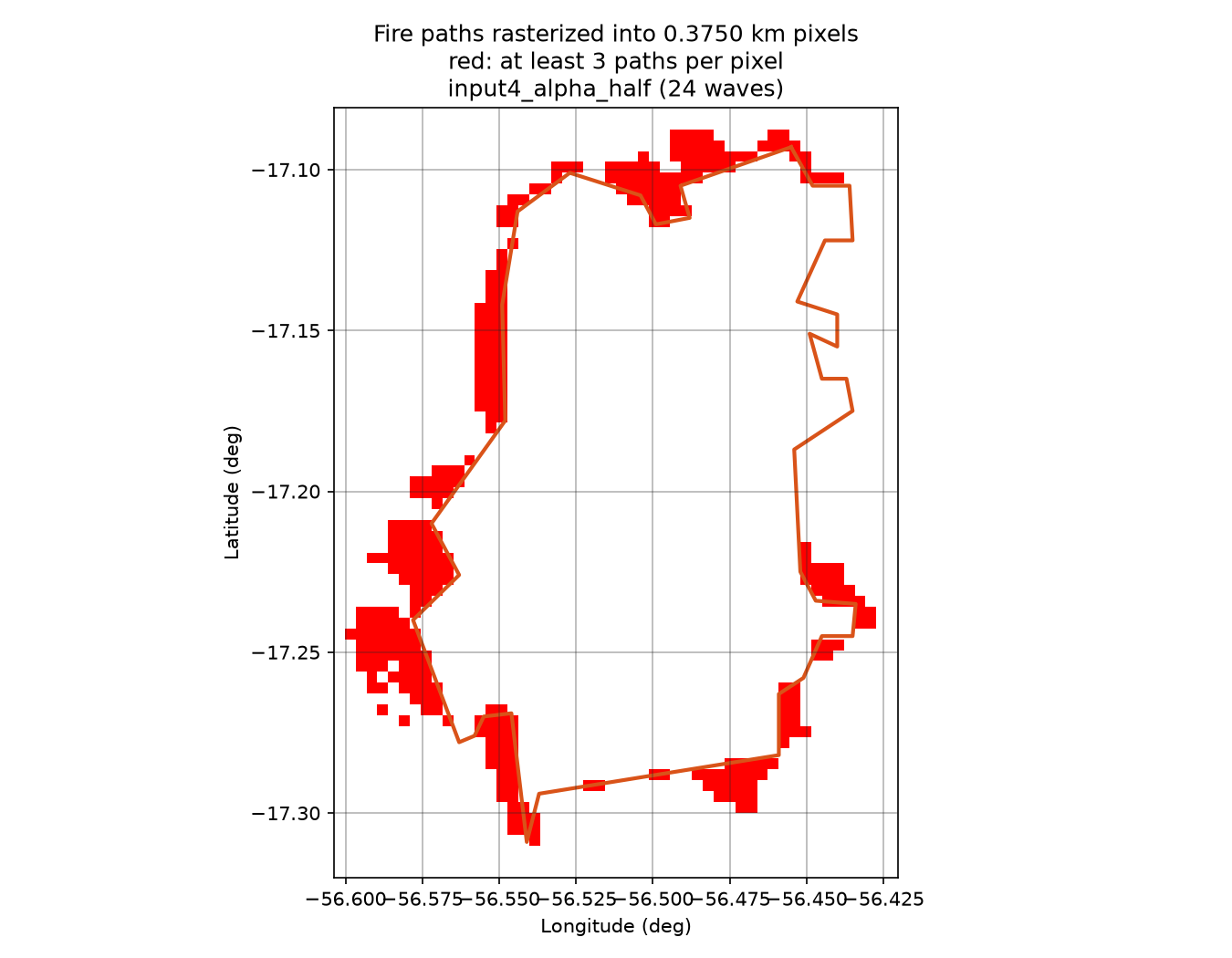}\hfill
		\includegraphics[width=0.30\textwidth]{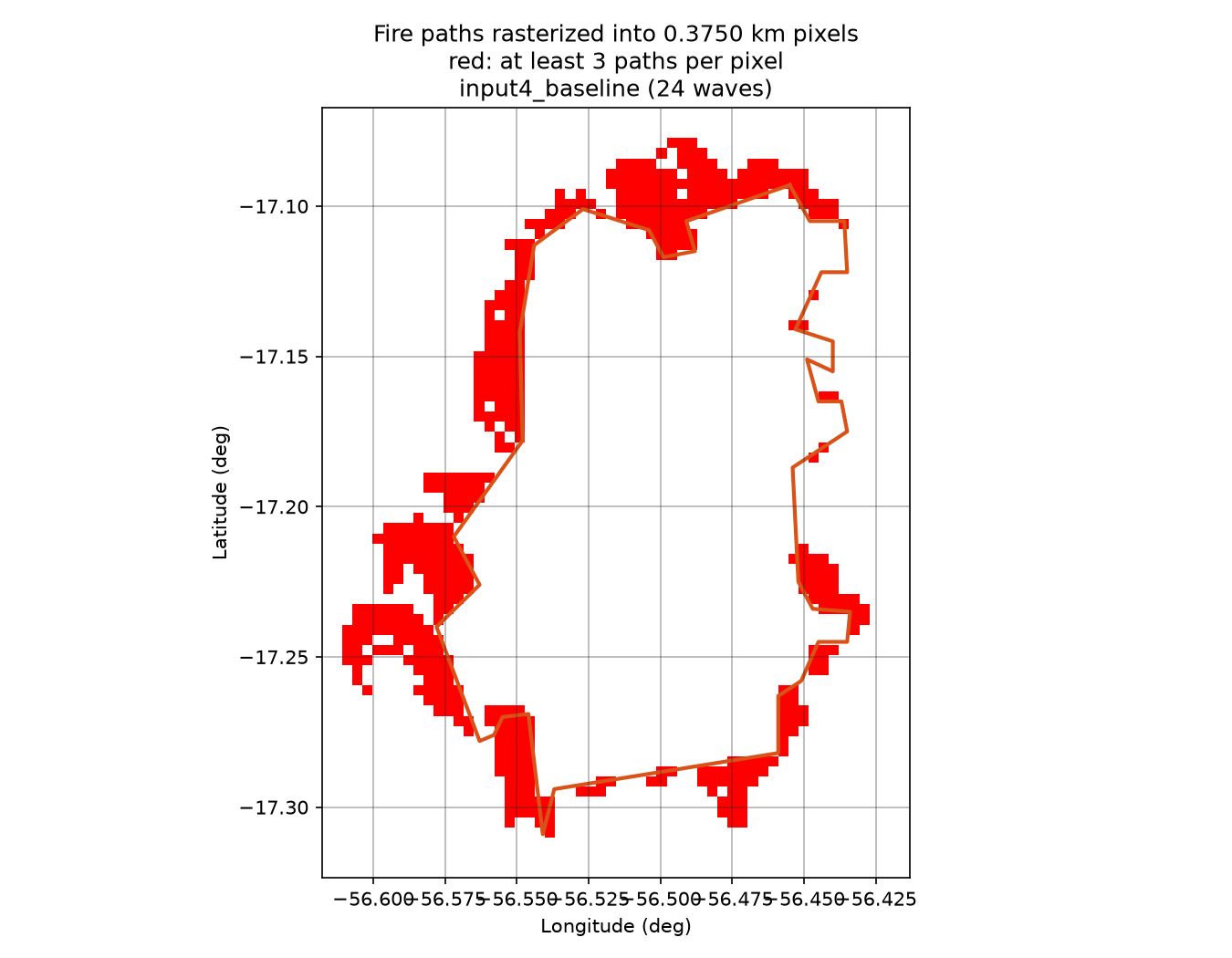}\hfill
		\includegraphics[width=0.30\textwidth]{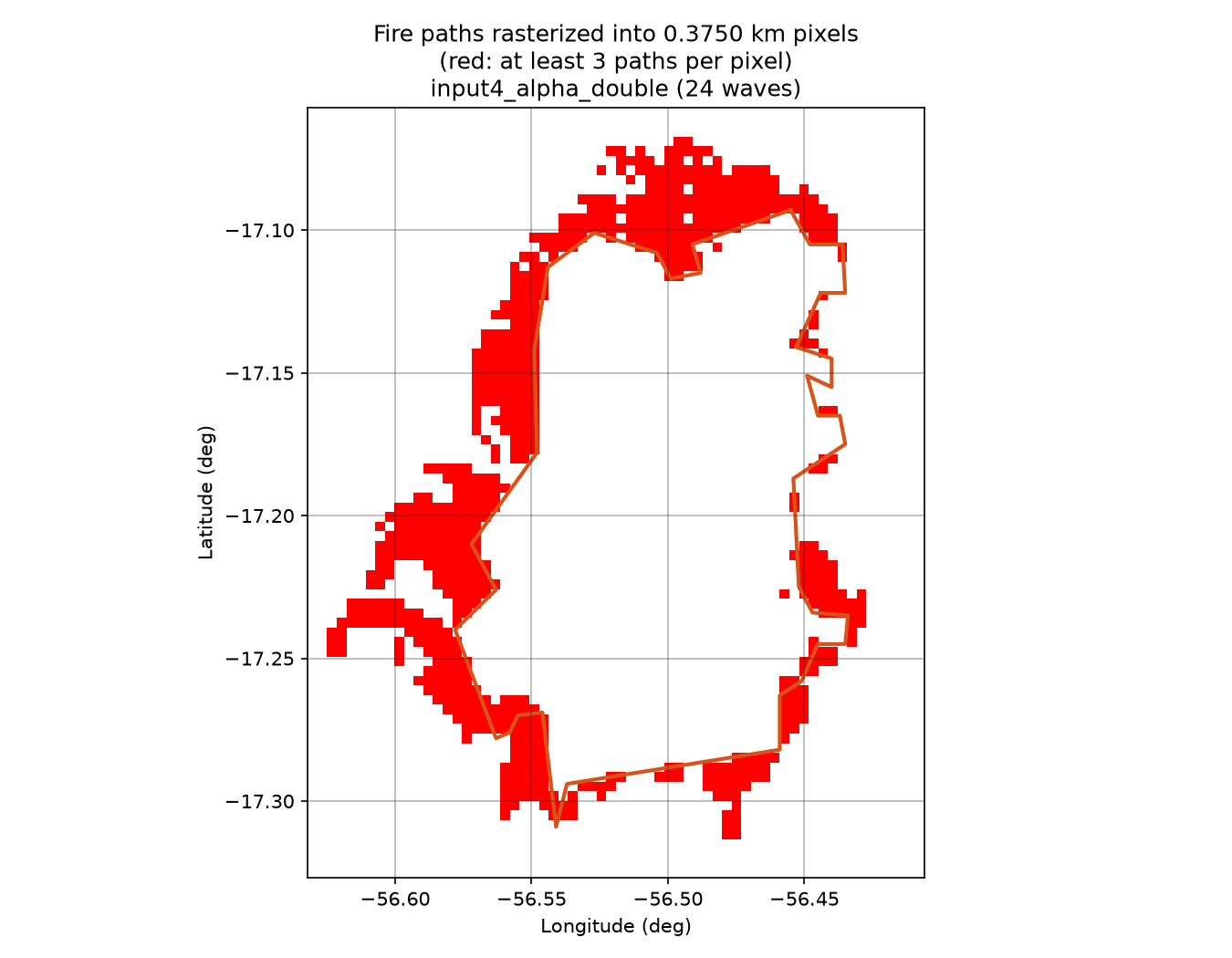}\par
		\noindent \rev{\textbf{$\beta$}}\par
		\includegraphics[width=0.30\textwidth]{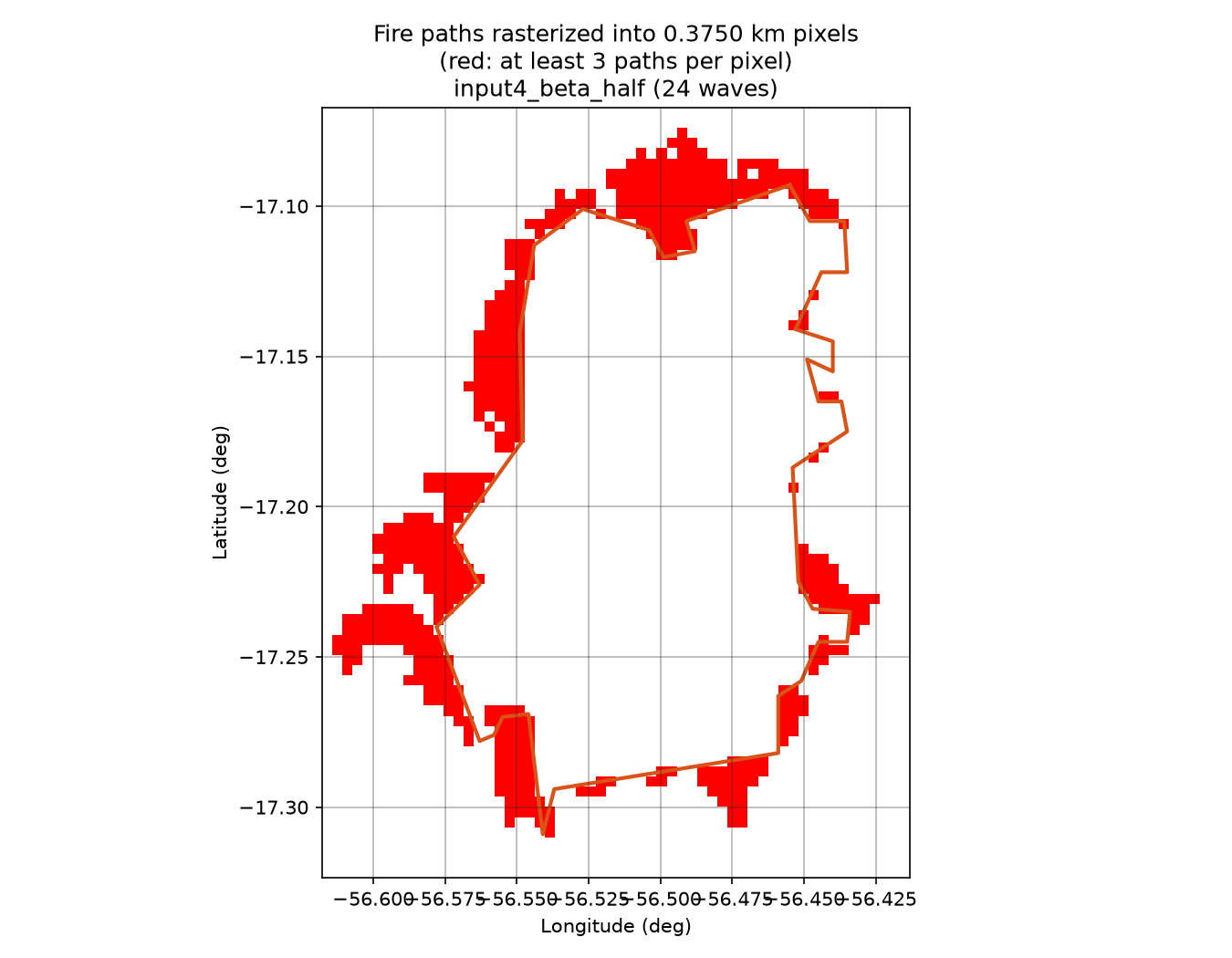}\hfill
		\includegraphics[width=0.30\textwidth]{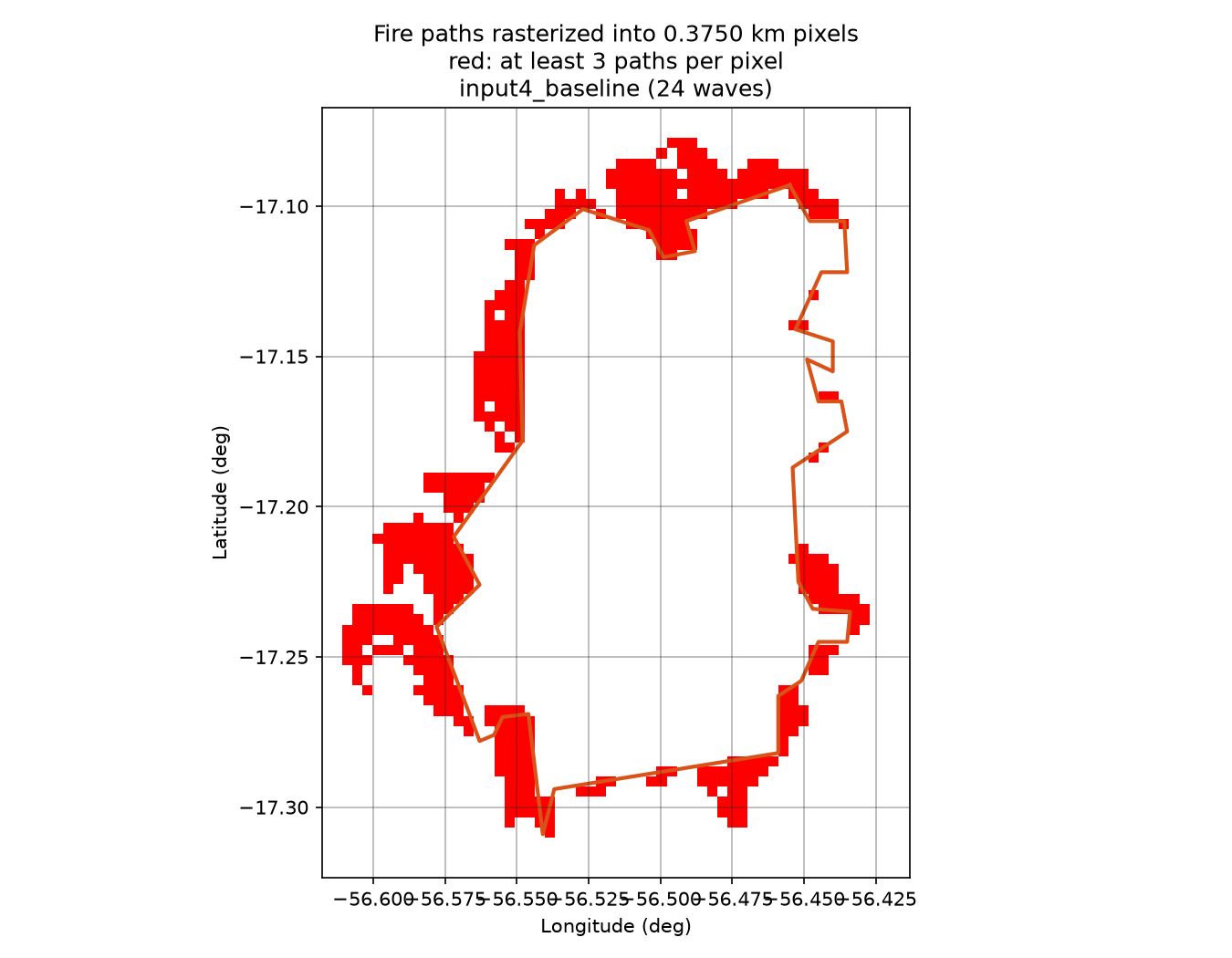}\hfill
		\includegraphics[width=0.30\textwidth]{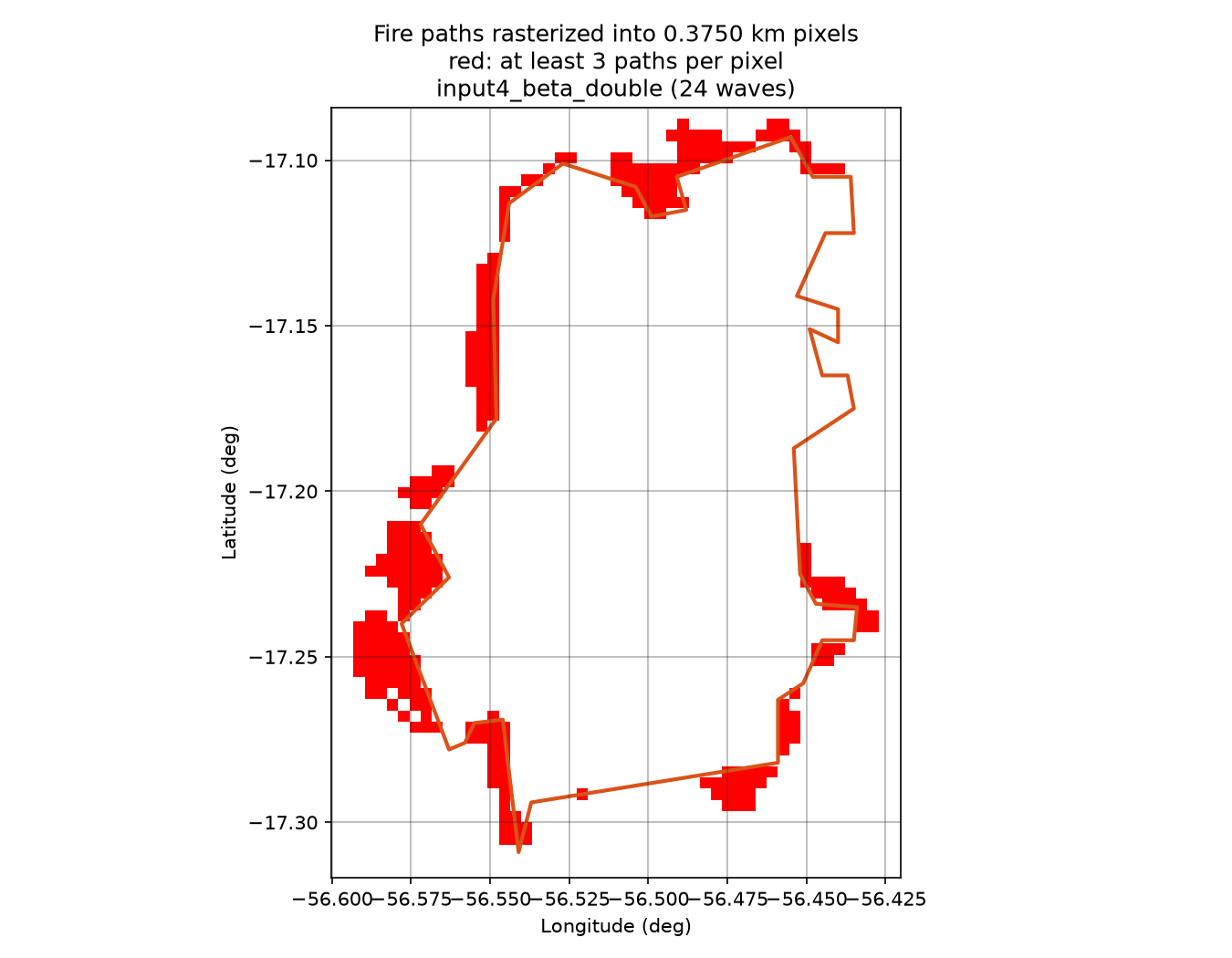}\par
		\noindent \rev{\textbf{$\gamma$}}\par
		\includegraphics[width=0.30\textwidth]{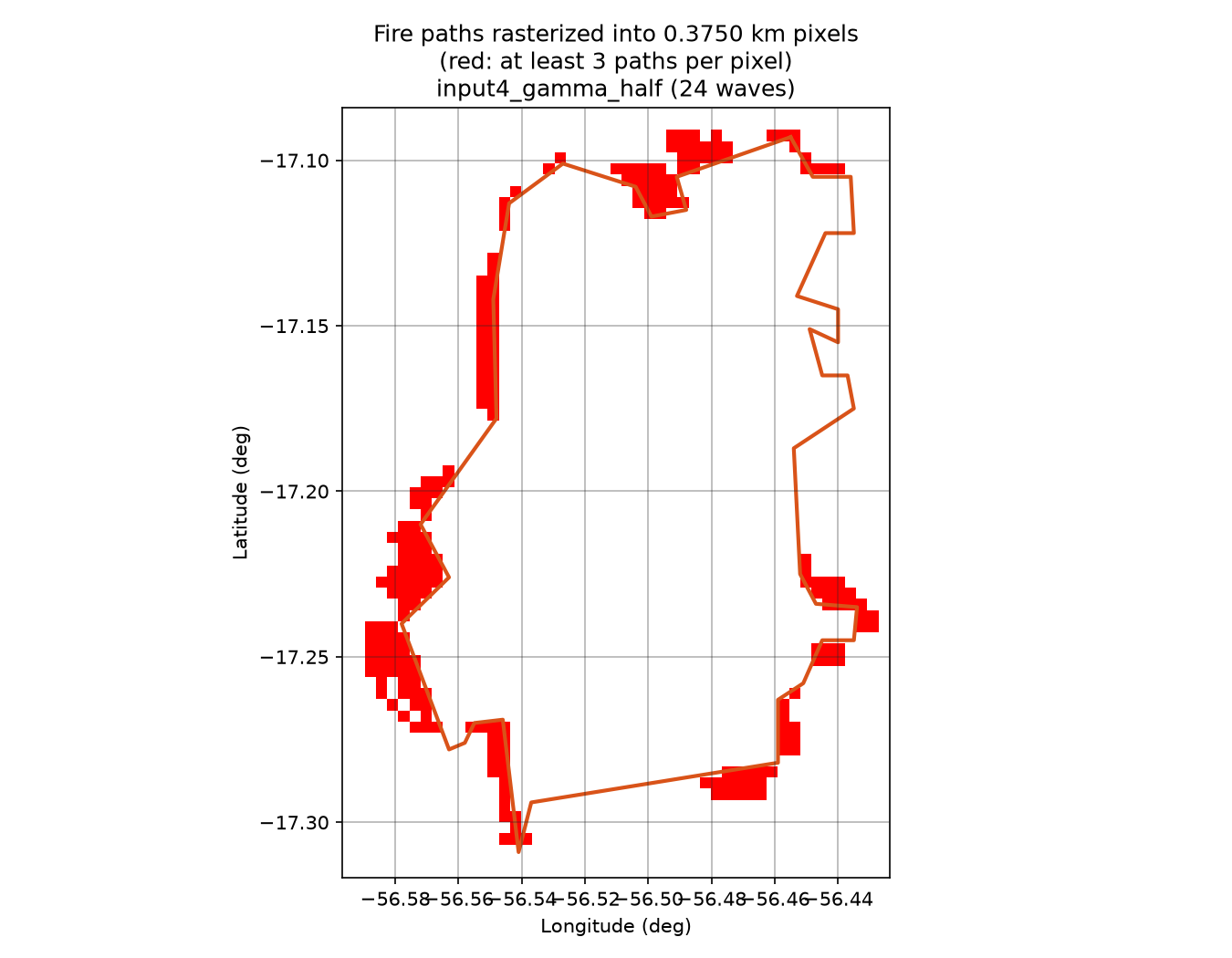}\hfill
		\includegraphics[width=0.30\textwidth]{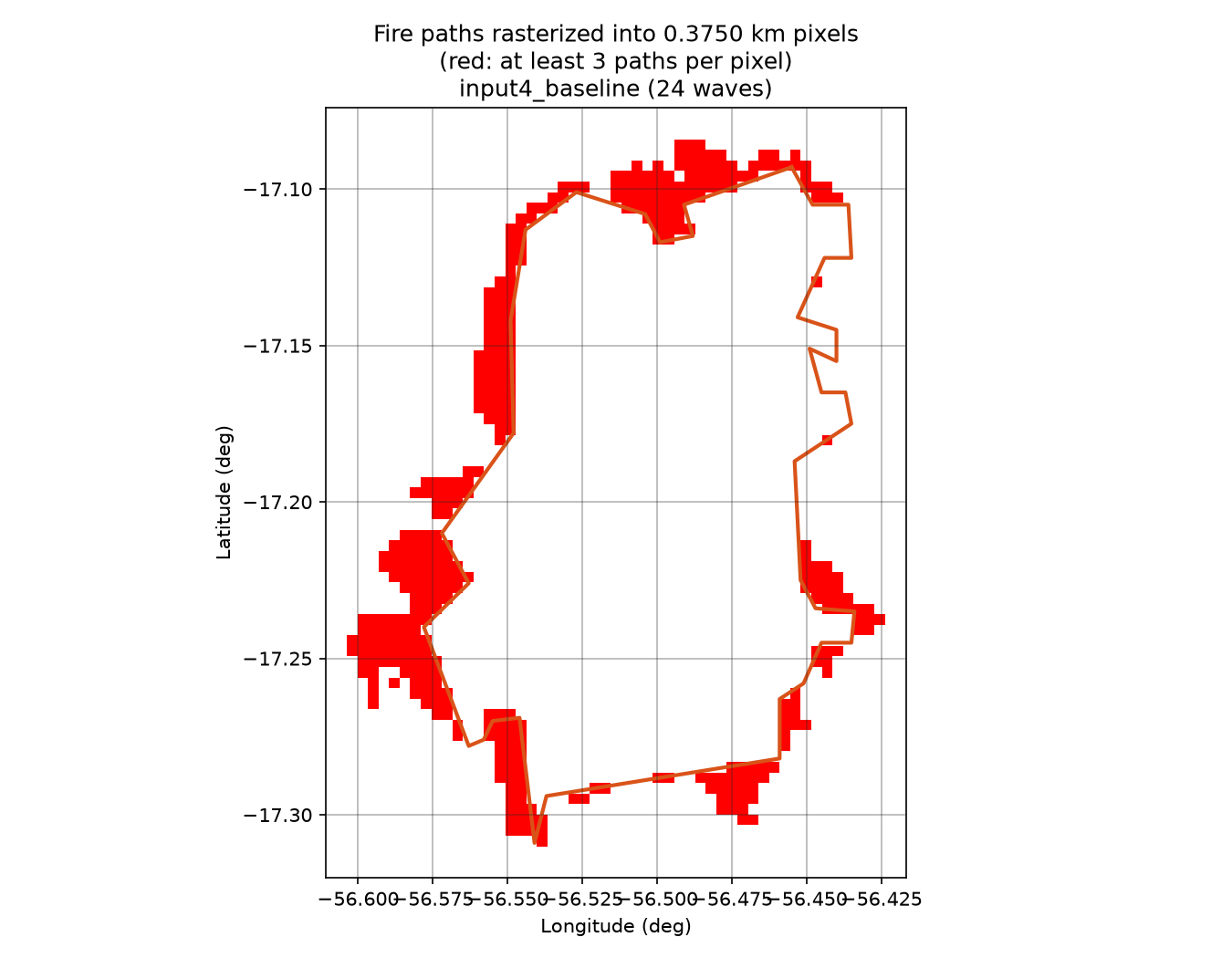}\hfill
		\includegraphics[width=0.30\textwidth]{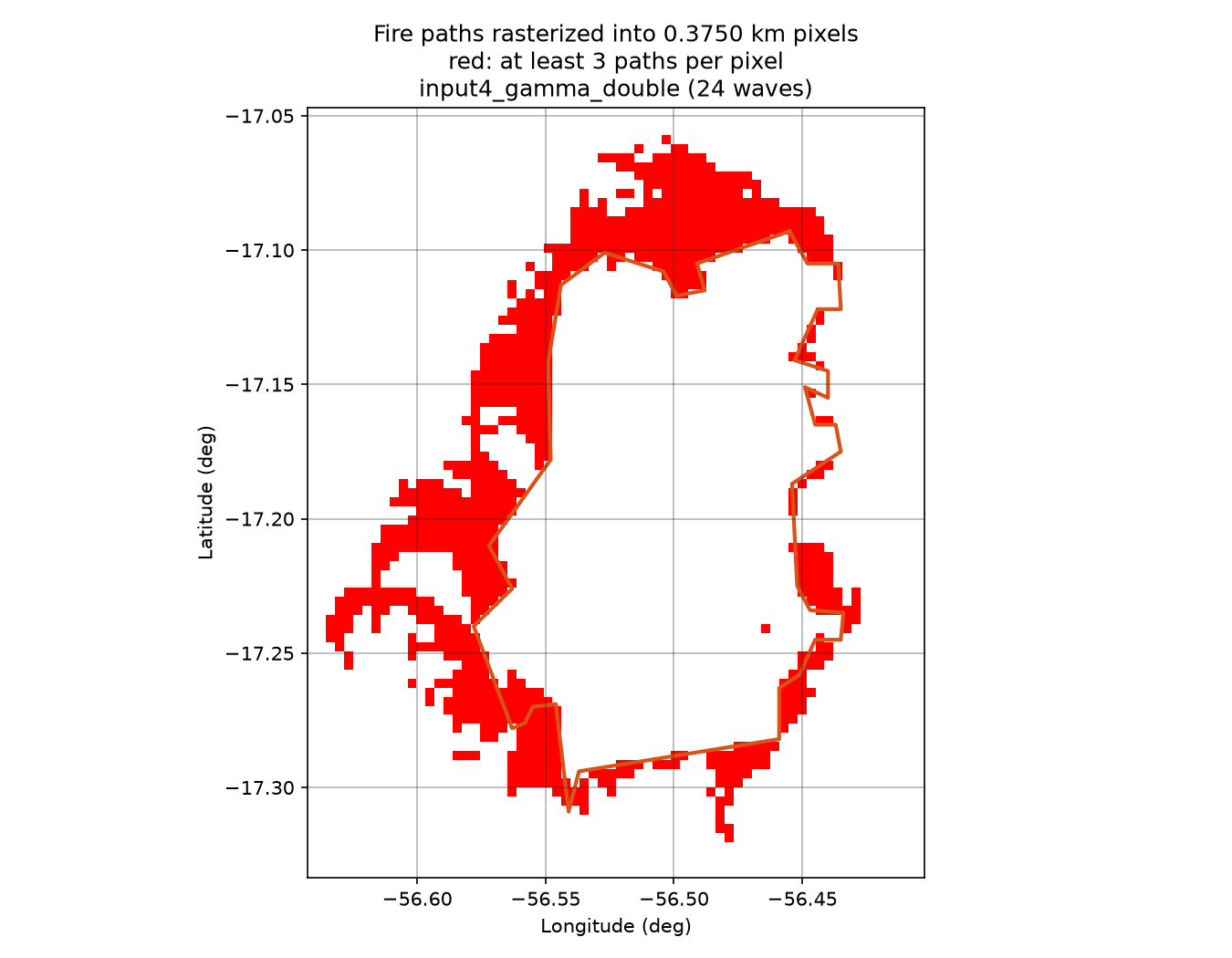}\par
		\noindent \rev{\textbf{Wind speed $U$}}\par
		\includegraphics[width=0.30\textwidth]{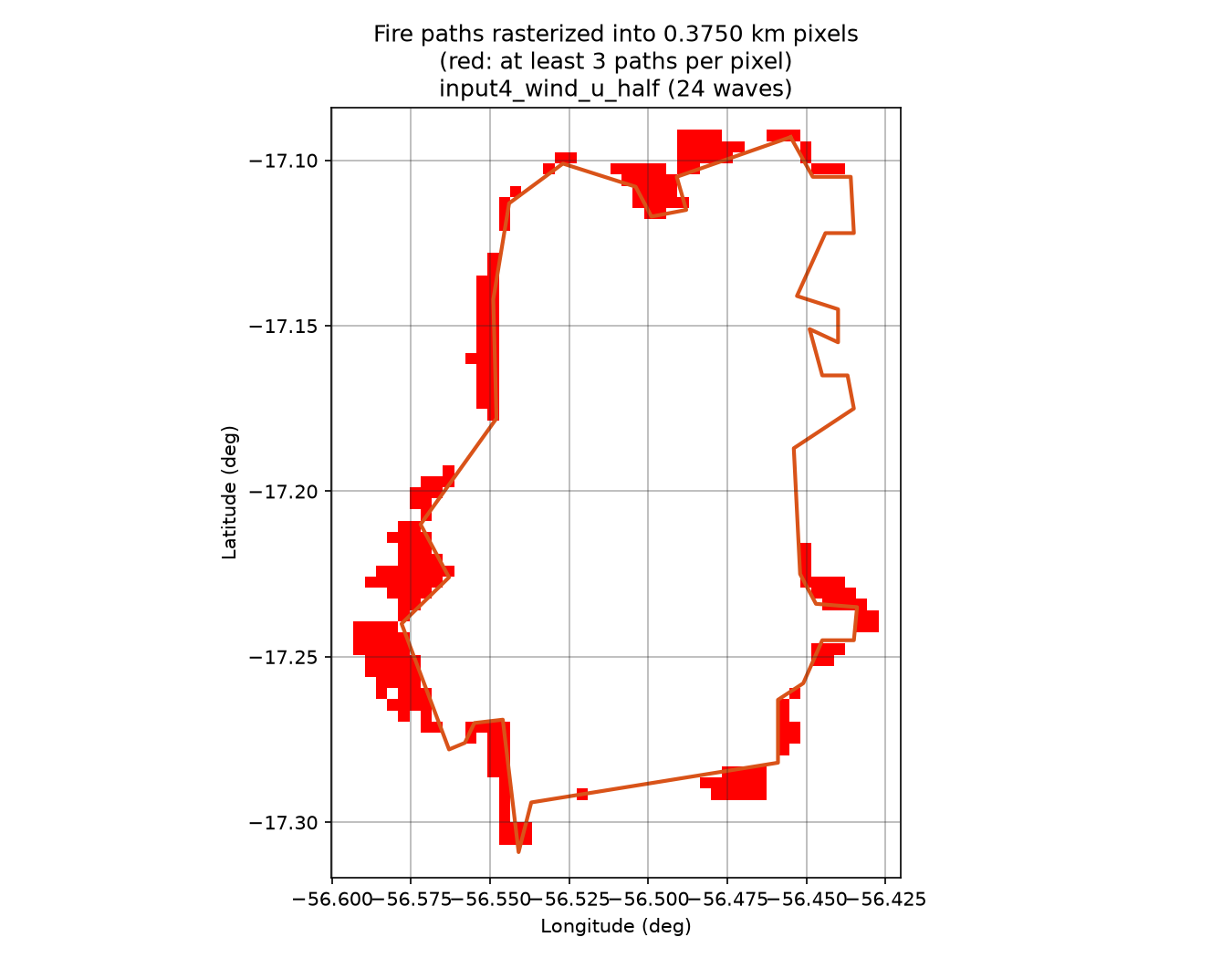}\hfill
		\includegraphics[width=0.30\textwidth]{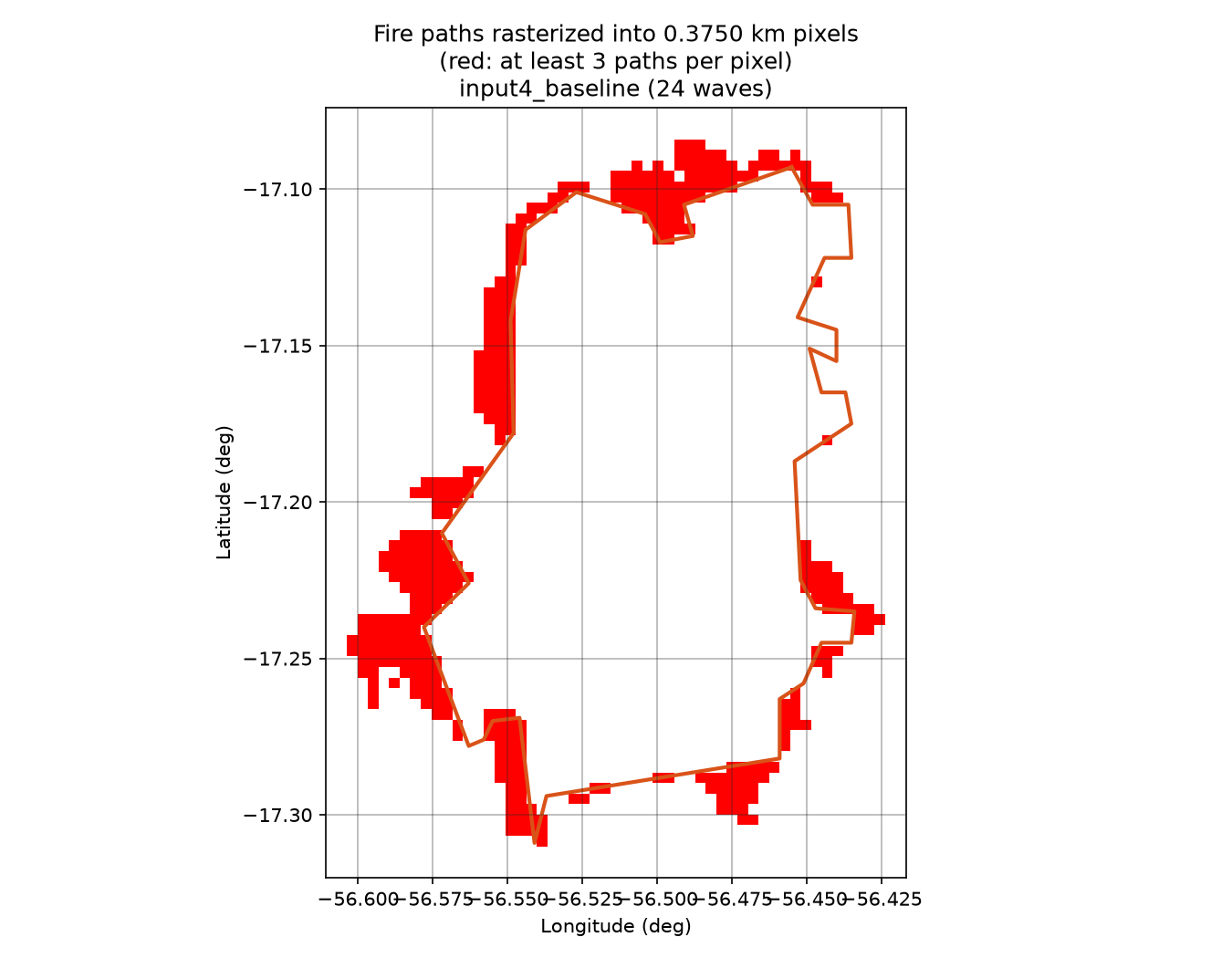}\hfill
		\includegraphics[width=0.30\textwidth]{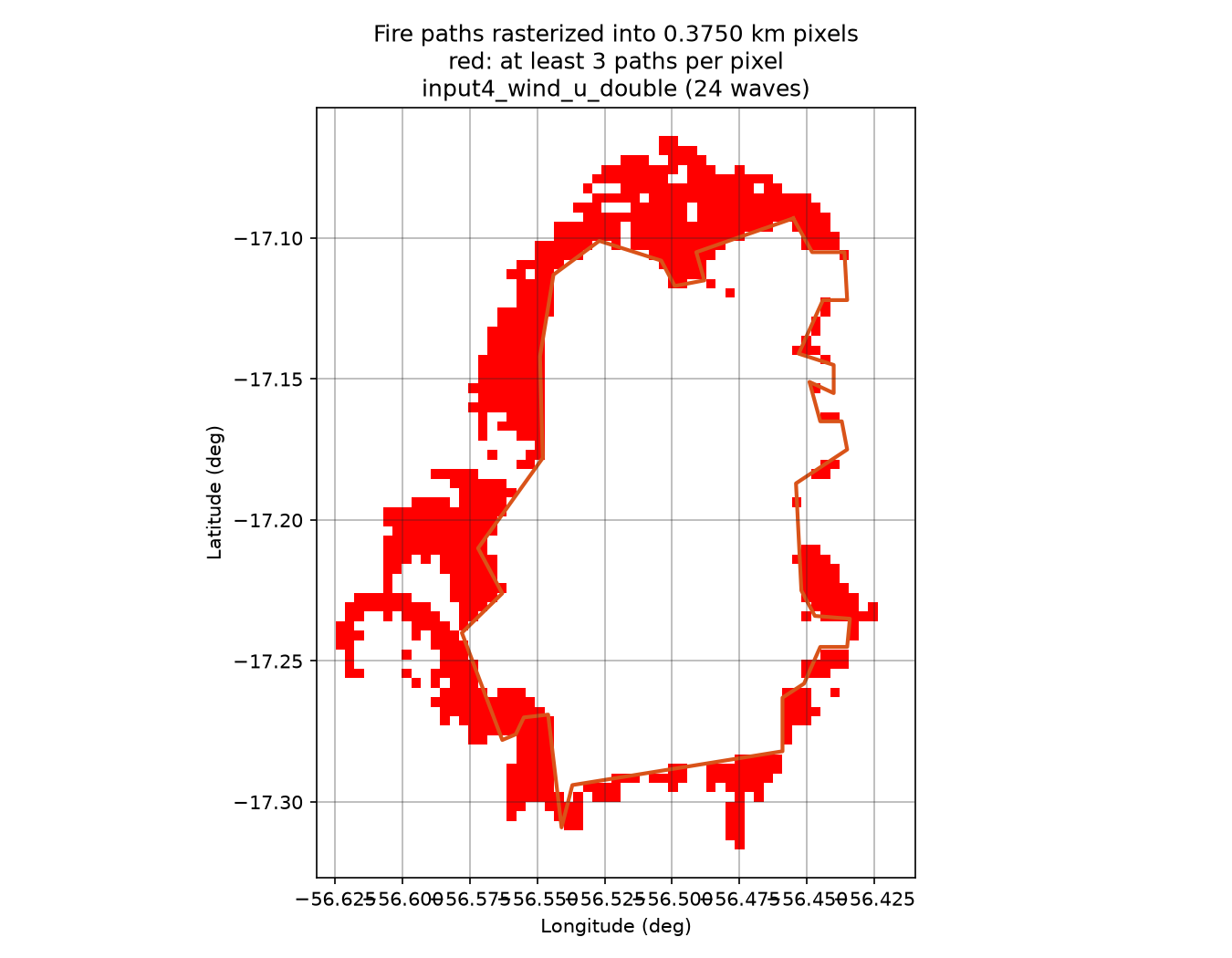}\par
		\caption{\rev{One-at-a-time sensitivity of calculated propagation for September 4, 2020. Rows show $\alpha$, $\beta$, $\gamma$, and wind speed $U$; columns show half, reference, and double values. The plotted numerical outputs are unchanged. Red pixels denote cells crossed by at least three calculated trajectories; the orange curve denotes the initial perimeter.}}
		\label{fig:pixel_sens_sep04}
	\end{figure}
	
	\begin{figure}[p]
		\centering
		\captionsetup{font=small}
		\noindent \rev{\textbf{$\alpha$}}\par
		\includegraphics[width=0.30\textwidth]{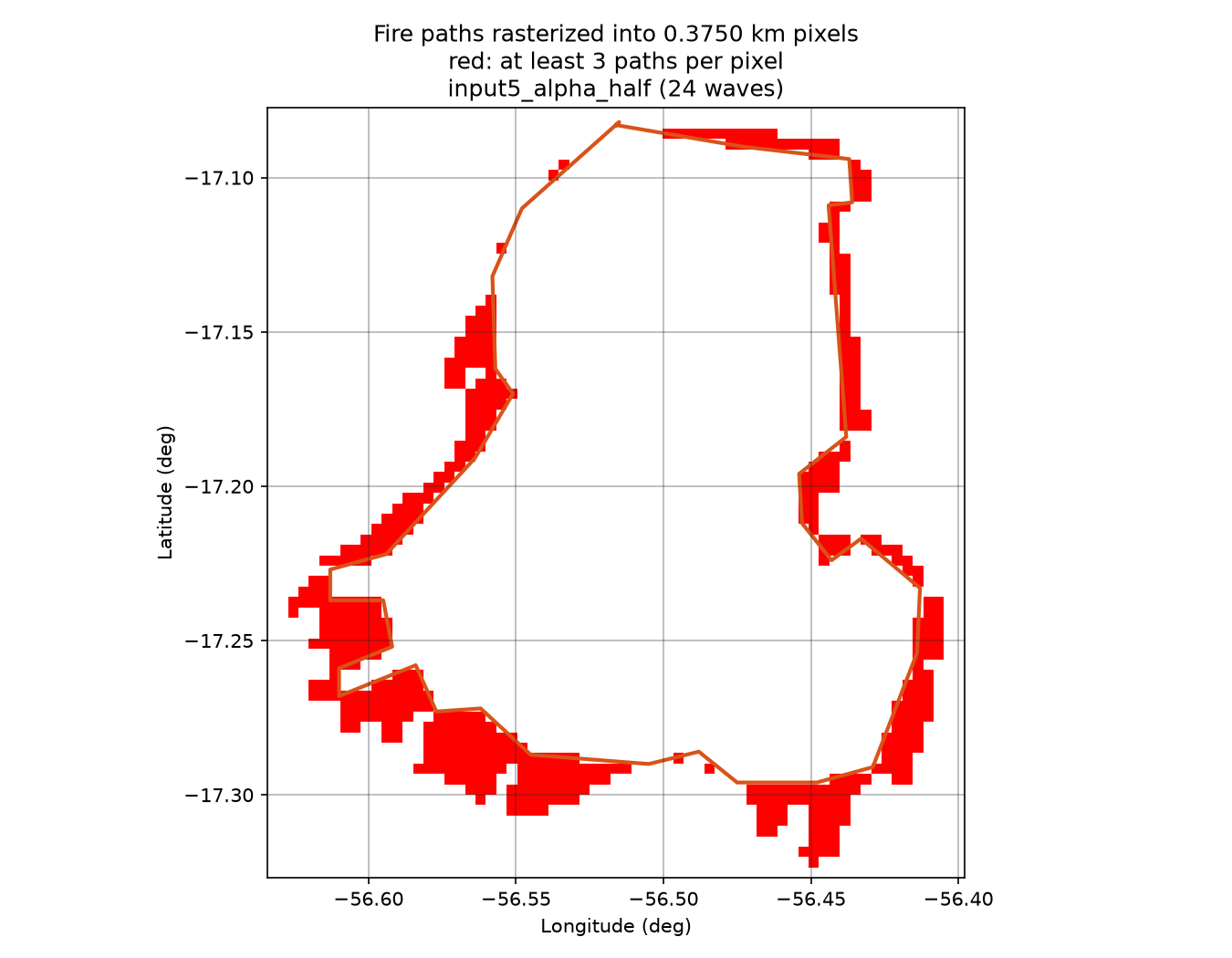}\hfill
		\includegraphics[width=0.30\textwidth]{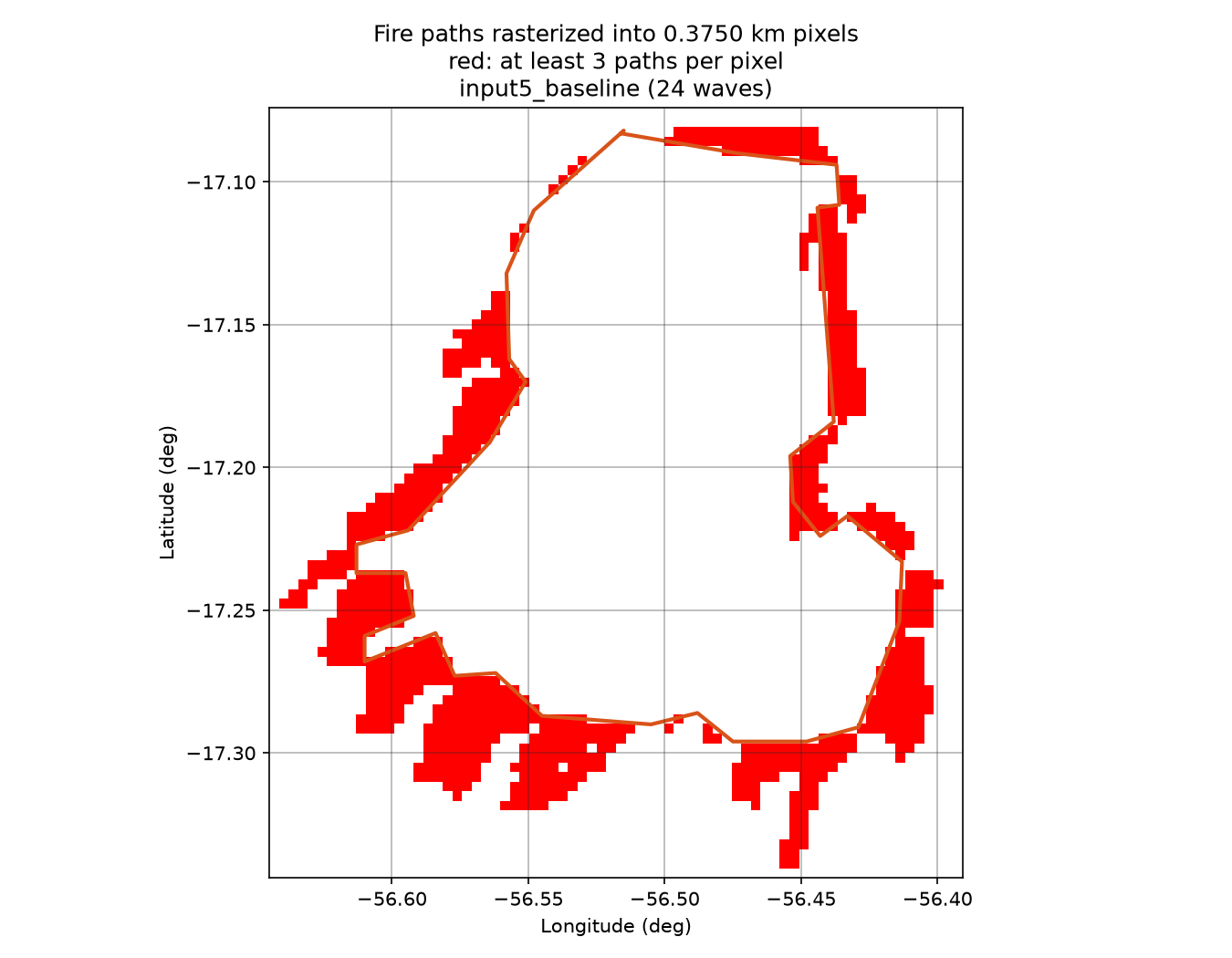}\hfill
		\includegraphics[width=0.30\textwidth]{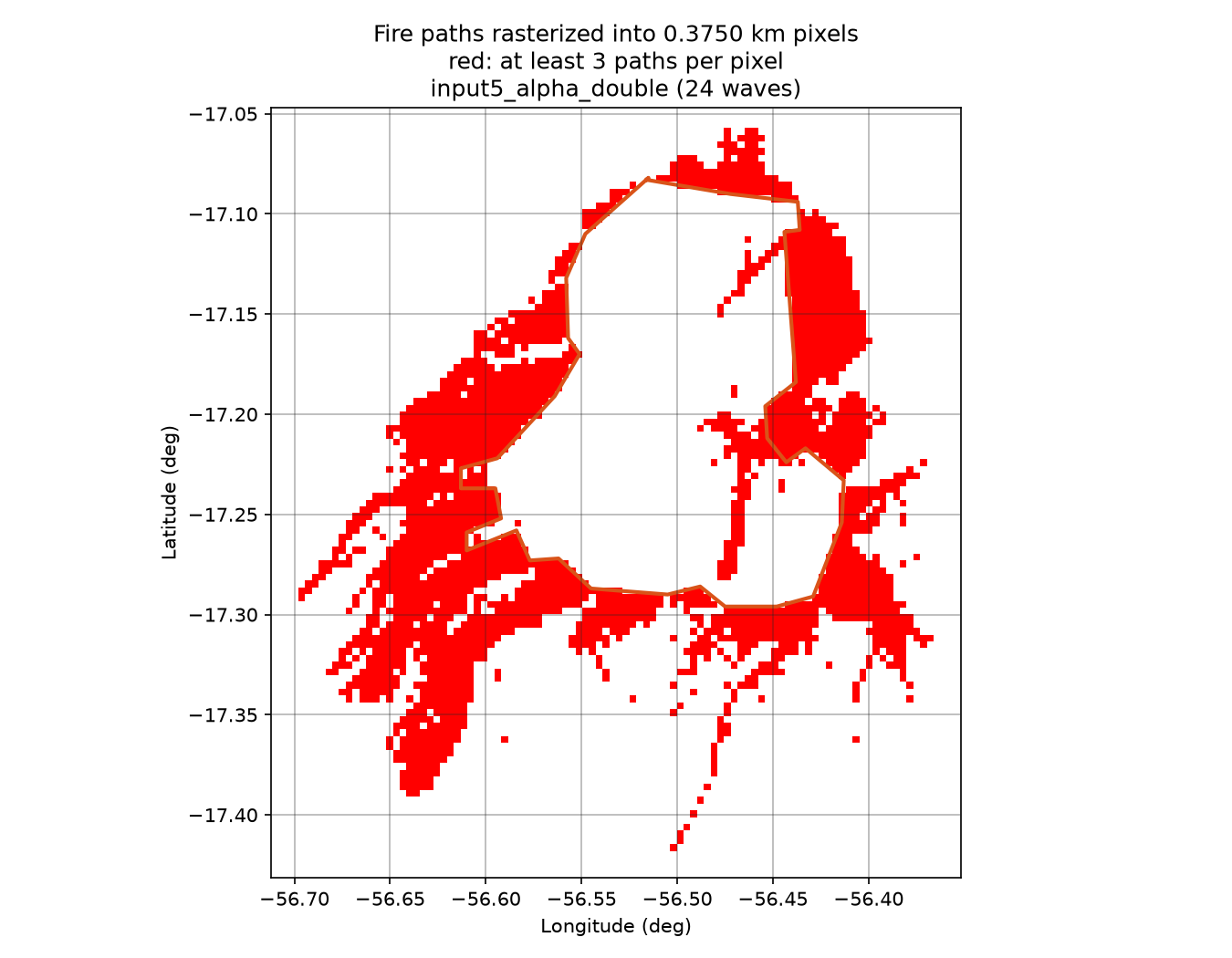}\par
		\noindent \rev{\textbf{$\beta$}}\par
		\includegraphics[width=0.30\textwidth]{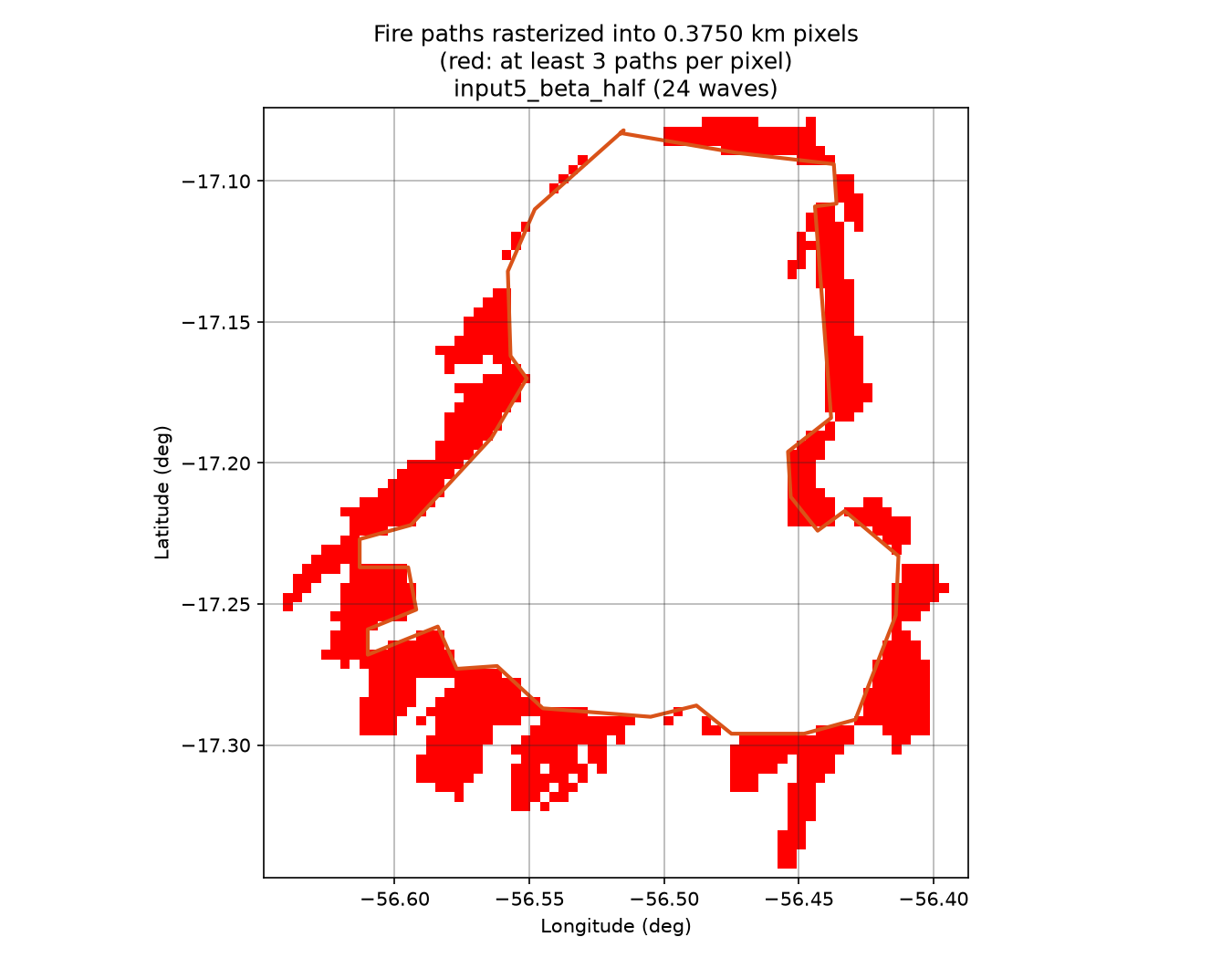}\hfill
		\includegraphics[width=0.30\textwidth]{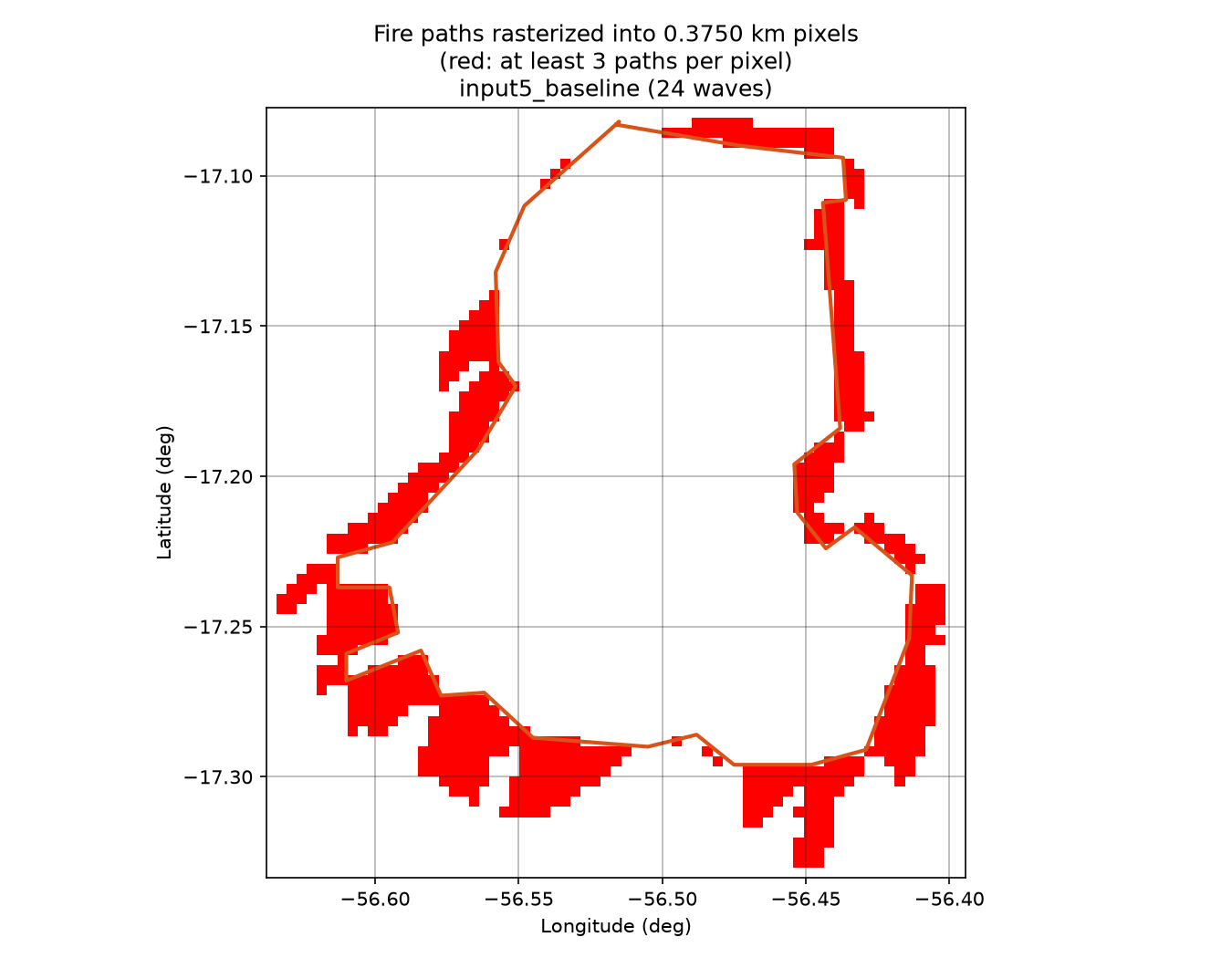}\hfill
		\includegraphics[width=0.30\textwidth]{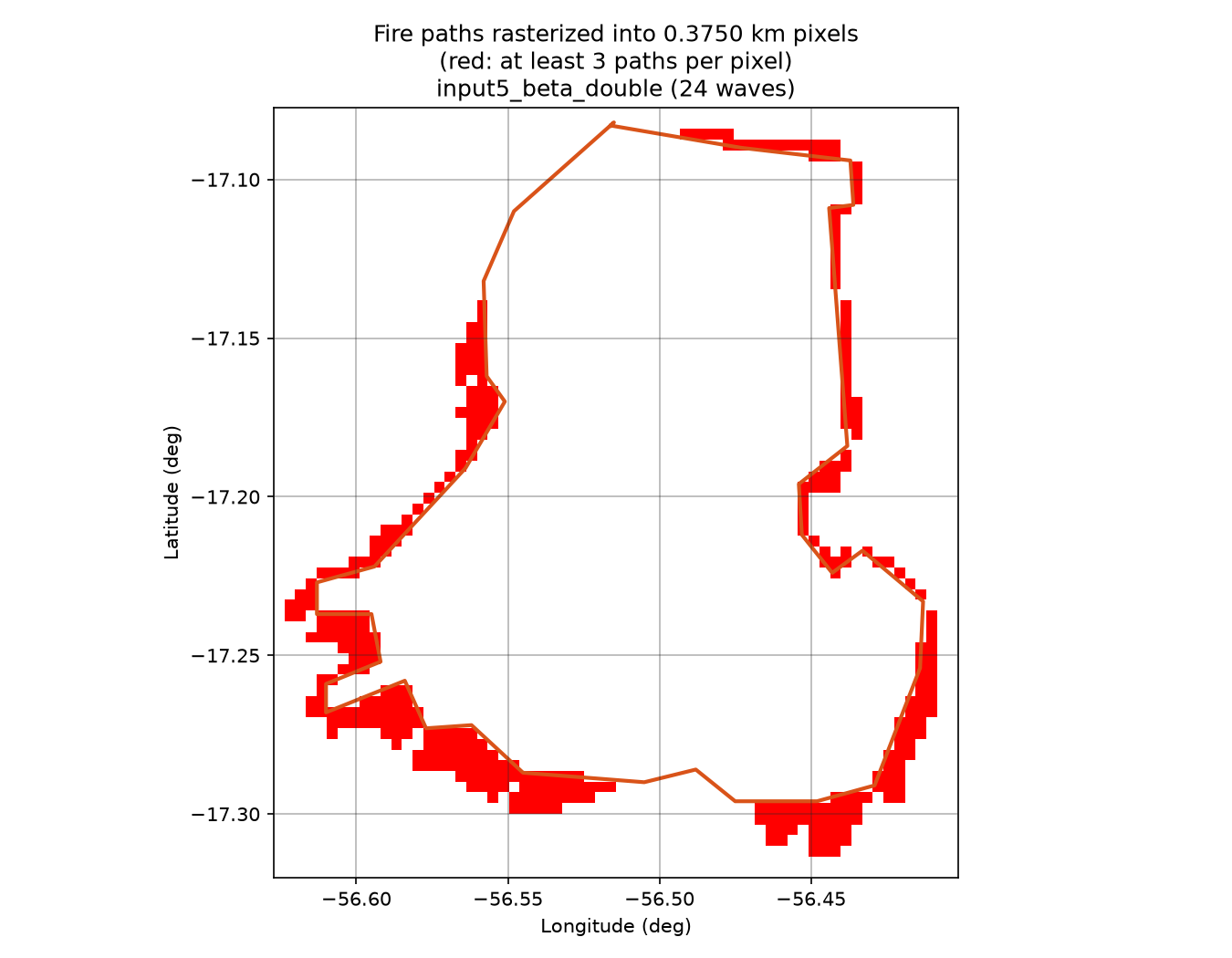}\par
		\noindent \rev{\textbf{$\gamma$}}\par
		\includegraphics[width=0.30\textwidth]{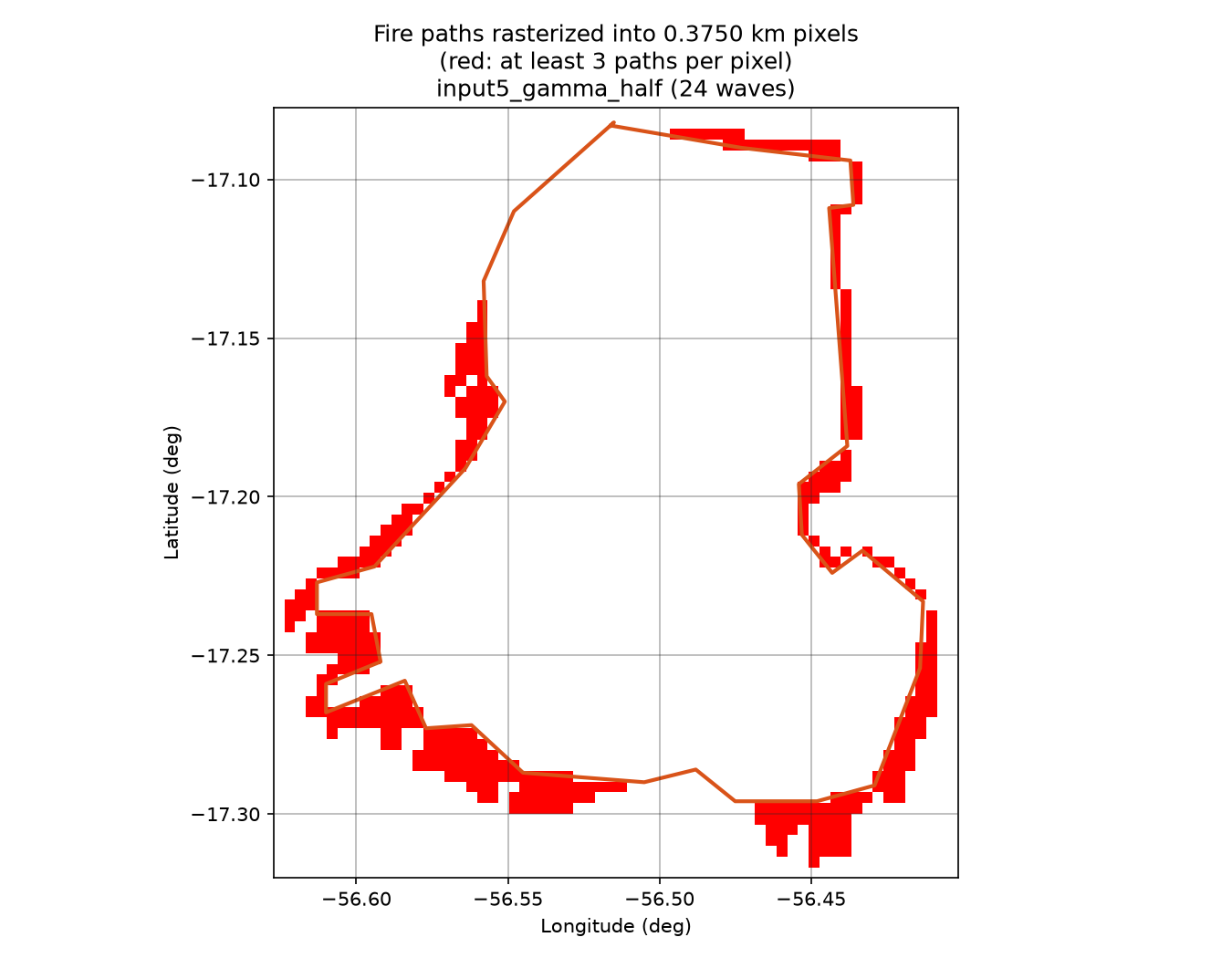}\hfill
		\includegraphics[width=0.30\textwidth]{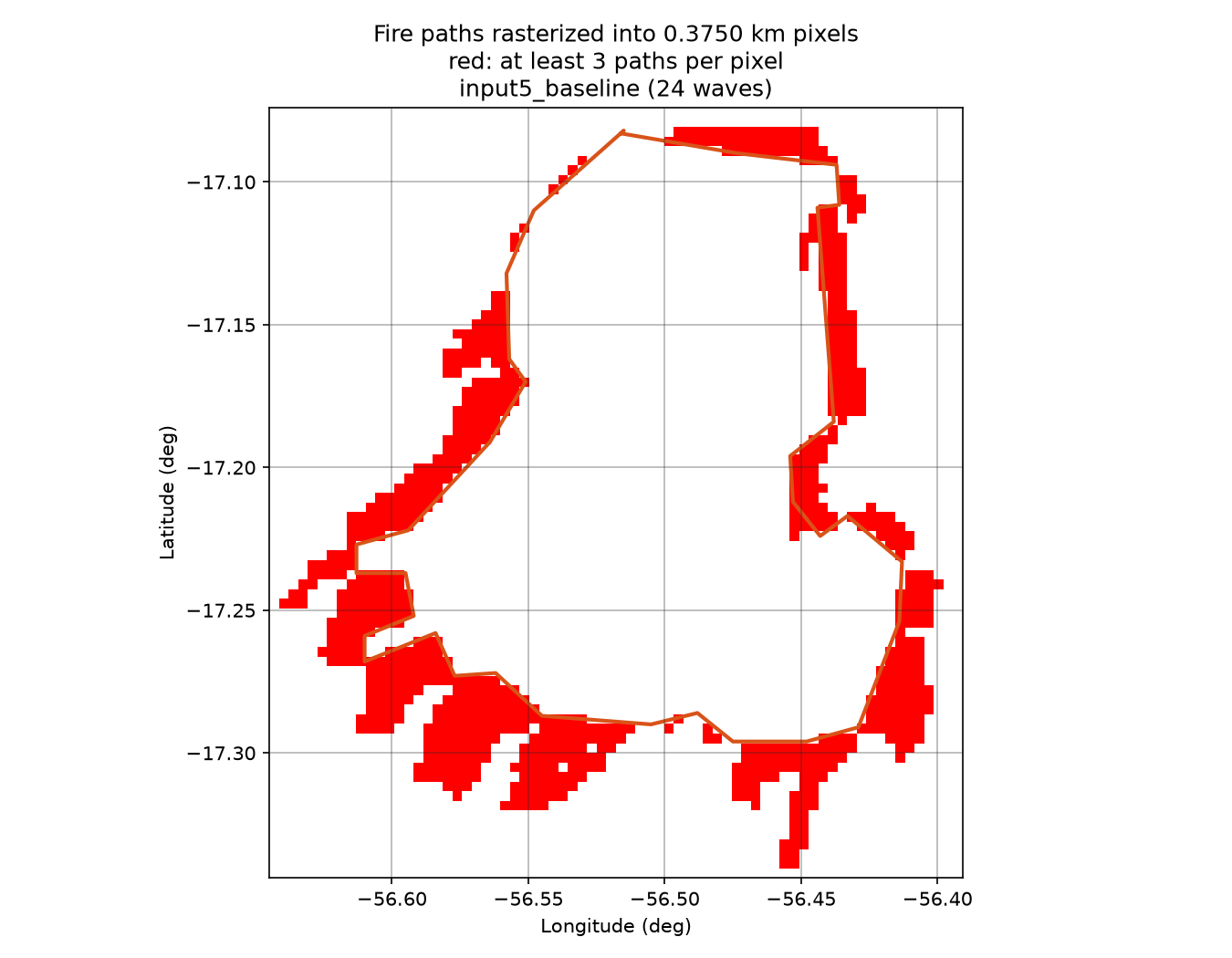}\hfill
		\includegraphics[width=0.30\textwidth]{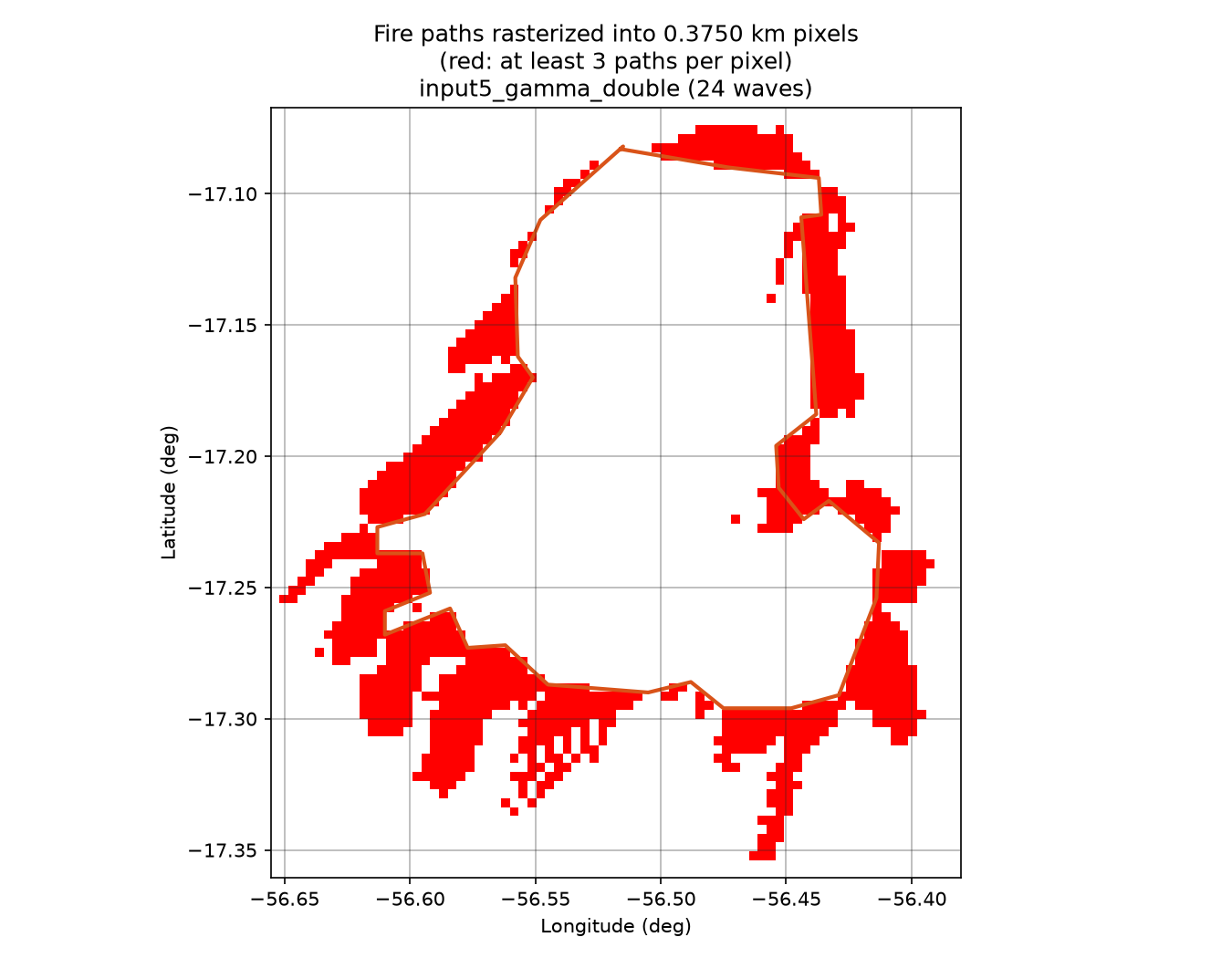}\par
		\noindent \rev{\textbf{Wind speed $U$}}\par
		\includegraphics[width=0.30\textwidth]{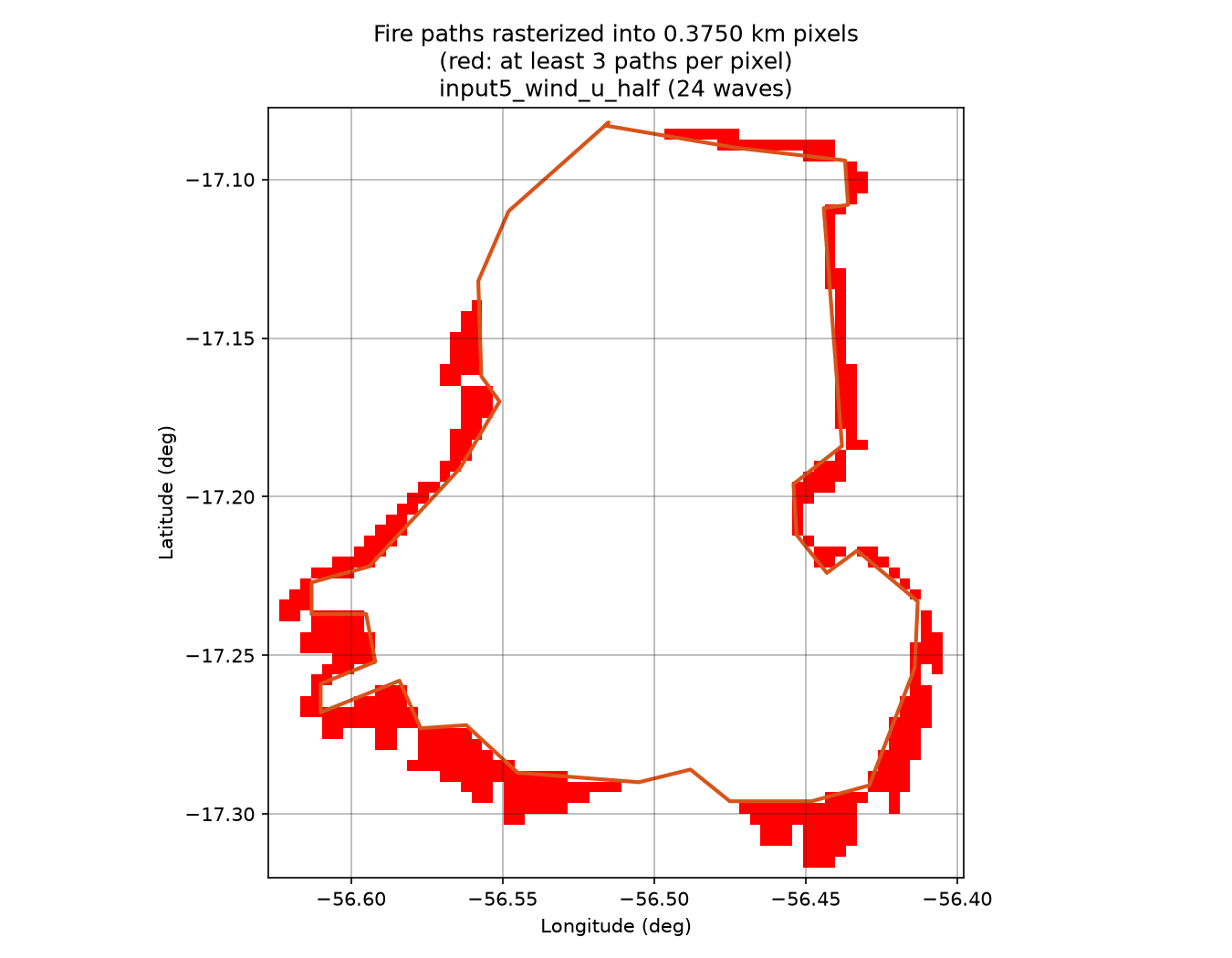}\hfill
		\includegraphics[width=0.30\textwidth]{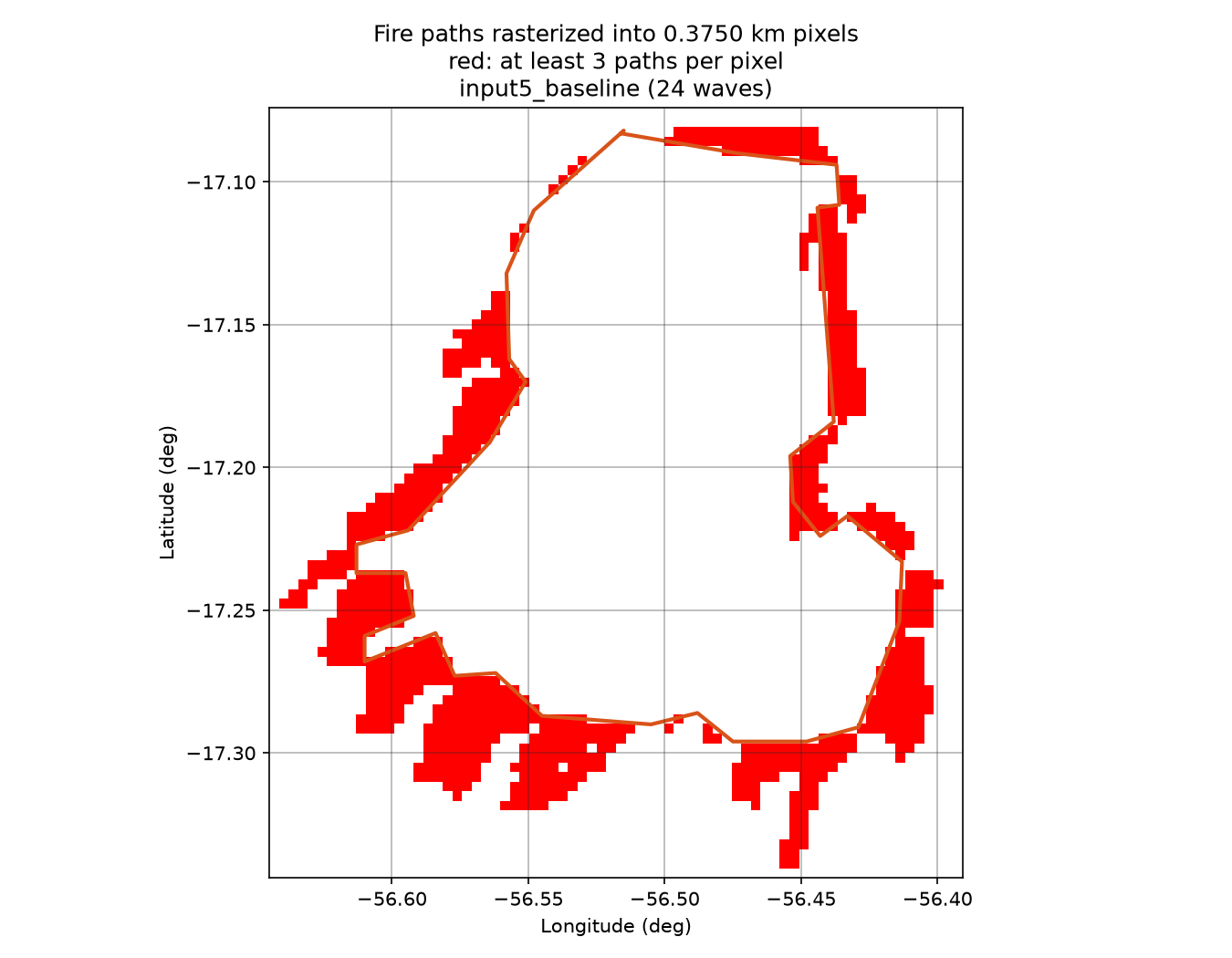}\hfill
		\includegraphics[width=0.30\textwidth]{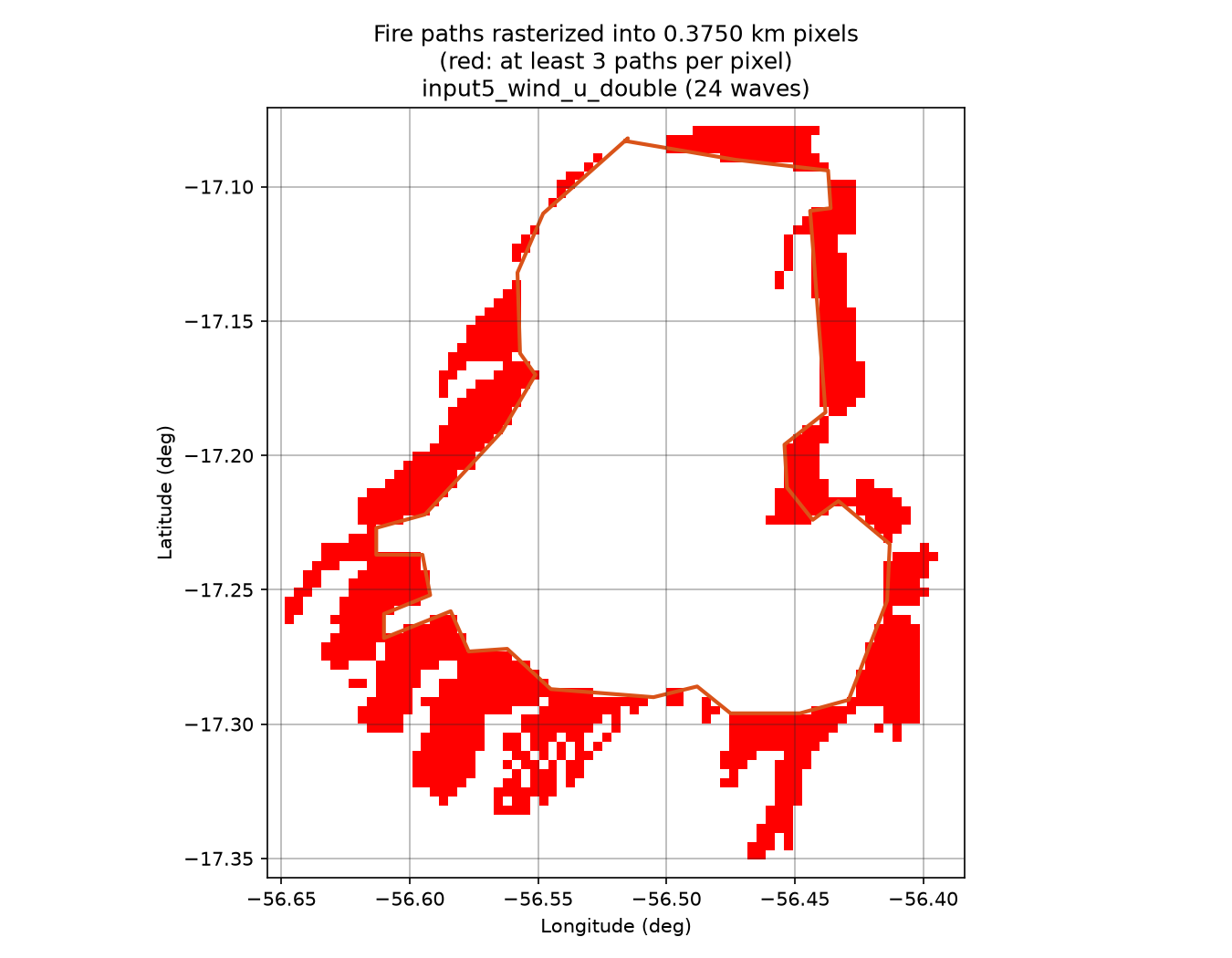}\par
		\caption{\rev{One-at-a-time sensitivity of calculated propagation for September 5, 2020. Rows show $\alpha$, $\beta$, $\gamma$, and wind speed $U$; columns show half, reference, and double values. The plotted numerical outputs are unchanged. Red pixels denote cells crossed by at least three calculated trajectories; the orange curve denotes the initial perimeter.}}
		\label{fig:pixel_sens_sep05}
	\end{figure}
	
	\rev{Figures~\ref{fig:pixel_sens}, \ref{fig:pixel_sens_sep04}, and \ref{fig:pixel_sens_sep05} show} a consistent sensitivity pattern across
	the three experiments. Changes in $\alpha$, $\gamma$, and wind speed
	mainly affect the extent of propagation, with larger values generally
	producing broader calculated fire patterns. The effect of $\beta$ is more
	moderate under the prescribed moisture condition. These results 
	\rev{characterize the response to the four magnitude-related factors. The computed patterns show a consistent qualitative response under the tested factor-two perturbations.}
	
	Because $\gamma$ is obtained through calibration, its influence was
	examined further by scaling the complete field by factors of $0.5$, $0.8$,
	$1.0$, $1.2$, and $2.0$. The corresponding Dice coefficients were
	$0.32$, $0.49$, $0.58$, $0.60$, and $0.50$ for September~3;
	$0.22$, $0.35$, $0.57$, $0.54$, and $0.23$ for September~4; and
	$0.24$, $0.38$, $0.47$, $0.38$, and $0.28$ for September~5.
	Thus, the calibrated field lies at or near the highest agreement for
	September~4 and 5, whereas September~3 shows some improvement for larger
	values of $\gamma$. These results indicate that the calibrated spatial
	coefficient provides a reasonable propagation magnitude for the studied
	event while also showing the expected uncertainty associated with its
	estimation.


	\FloatBarrier
	\section{Discussion}
	
	\subsection{Interpretation of the Spatially Varying Propagation}
	
	The three experiments demonstrate that the proposed trajectory formulation
	can represent strongly heterogeneous fire propagation along irregular
	perimeters. The framework combines VIIRS observations, hourly ERA5 wind
	forcing, ERA5-Land vegetation information, and prescribed fuel moisture.
	Spatial variability is represented through the effective coefficient
	$\gamma(x,y)$, while the local rate-of-spread relation retains the established
	dependence of fire propagation on wind, moisture, and fuel conditions. This provides a
	compact alternative when the detailed spatial fuel parameters required for
	a complete pointwise physical parameterization are unavailable.
	
	An important feature of the formulation is that the propagation equations
	are solved locally along the evolving perimeter. Consequently, spatial
	variations in $\gamma(x,y)$ and temporal variations in wind can produce
	different propagation responses along the same fire front without imposing
	a uniform rate of spread or a fixed global propagation shape. The numerical
	experiments indicate that this construction can reproduce relevant spatial
	features of the observed fire development, including the coexistence of
	sectors with substantial advance and regions with limited propagation.
	
	\subsection{Prospective Sector-Wise Risk and Management Analysis}
	
The resulting paths could be used to quantify the progress of a simulated fire front across user-specified sectors. Divide the original perimeter into sectors $S_1,\ldots,S_J$ and assign each new simulated pixel reached to a sector according to a documented nonoverlapping assignment rule. Let $A_j(T)$ be the modeled newly reached area attributed to sector $j$ at time $T$. Then, provided a denominator is nonzero, $p_j(T)=100A_j(T)/\sum_{k=1}^{J}A_k(T)$ is its percentage of the total modeled new area. These sector shares are potential simulation indicators, not observed fire percentages or probabilities. They could be used, in independent testing, to compare hypothetical risk-reduction and perimeter-management strategies, but they do not themselves prescribe firefighting interventions. No sector-share percentages or intervention effects are calculated in these experiments.
	
	\subsection{Scope and Limitations}
	
The proposed framework is presented as a mathematical model of data-assisted
fire-front propagation, and not as a complete physical description of
fire behaviour. This formulation can be used to combine the satellite
observations and environment data without the detailed
pointwise definition of all fuel and combustion properties. The effective
The $\gamma(x,y)$ coefficient offers a convenient way to account for the unresolved
The calibration is for the event in question and thus accounts for the spatial variability.
The fuel moisture is prescribed, as are processes such as fire spotting,
suppression work, small fuel discontinuities and thermal details
The transfer is beyond the scope of the present work. Topographic effects are also disregarded.
Because the study area is considering the flat terrains.
Also the available satellite temporal and spatial resolutions are
and the level of detail that can naturally be constrained by the re-analysis products
represented and evaluated.

Numerical experiments are designed to study the mathematical and computational behavior of the proposed framework in realistic environmental conditions. The 12-h calibration and 12-h continuation periods allow us to study the propagation process with a given spatial coefficient $\gamma(x,y)$ during the continuation period and provide useful information about the model behavior outside the parameter-adjustment period. The VIIRS observations are instantaneous snapshots taken by the satellite at a given time, not measurements of the two stages made independently at different times, so the comparisons are interpreted as assessments of spatial agreement within an event, not independent predictive validation. Thus, the pixel-level metrics reported provide quantitative evidence of agreement between simulated patterns and observed fire activity, while taking into account the observational characteristics of VIIRS data. Further investigations with precisely timed observations, uncertainty analyses and independent validation datasets will provide further opportunities for assessment and extension of the framework.

The current study focuses intentionally on
{one Pantanal wildfire} to derive the mathematical formulation, analyze its numerical implementation and study its ability to represent spatially heterogeneous fire propagation. One can investigate the influence of local propagation parameters, environmental heterogeneity and time-dependent wind forcing in a standard experimental setup. An important next step is to expand the analysis to multiple wildfire events or to compare it with established propagation models, especially when comparable time resolved fire-perimeter observations and environmental datasets are available. The framework proposed in this paper provides a mathematical basis to explore the fire-front propagation in different sectors and to determine the relative contributions of such sectors to the future fire growth, in addition to the present numerical experiments. This information can be used to develop and evaluate wildfire risk management and mitigation plans pending further validation and quantification of uncertainty.

\section{Conclusion}
\label{sec:conclusion}

This paper proposes a mathematical framework for the propagation of wildfire fronts in heterogeneous landscapes based on trajectories. We present our model implementation in a multi-day wildfire case study in the Brazilian Pantanal.

 The main mathematical contribution is the derivation of a local nonlinear system which uses a direction dependent propagation law. The corresponding translated ellipse construction provides an explicit and unique outward supporting velocity for positive local spread rates and a given outward unit normal. This geometric result explains why the outward branch is selected independent of the iterative root-finding procedure. The closed form support velocity sets an exact mathematical consistency benchmark for the iterative root-finding procedure. 
 
 The framework relates continuous propagation trajectories to discrete satellite observations by rasterizing the modeled advance onto an ageoreferenced $375$-m grid and registering the resulting pixel domain to VIIRS active-fire detections. The Dice coefficients for the full-run daily experiments on 3, 4 and 5 September 2020 were 0.578, 0.571 and 0.467, respectively. These statistics demonstrate the degree of within-event spatial agreement achieved through observation-guided calibration during the initial \(12\)-h window, followed by fixed-parameter continuation over the subsequent \(12\)-h period. The satellite observations provide spatial information for sector-wise parameter estimation, supporting the assessment of the framework's ability to reproduce the observed fire evolution. The sensitivity experiments are done individually to study the response of the moving front to a factor of two change in the spatial propagation coefficient \(\gamma(x,y)\), the wind-speed magnitude \(U\), and the scaling parameters \(\alpha\) and \(\beta\).
 
  The effect on simulated fire extent emphasizes the nonlinear interaction between spatially heterogeneous spread rates and time-dependent atmospheric forcing. Magnitude perturbations will be augmented in future uncertainty studies with additive wind direction rotations and temporal offset ensembles. 
  
  Beyond the present case study, the computed trajectories establish a mathematical foundation for mapping simulated fire spread onto predefined perimeter sectors and investigating a range of hypothetical risk-management scenarios. This provides a promising basis for future extensions toward operational fire-spread assessment and the quantitative evaluation of potential intervention strategies. 
  
  Future work will consider terrain slope in the anisotropic metric structure and will explore the framework on time resolved independent wildfire datasets across various landscapes.
	\section*{Data Availability}
	
	The data supporting the findings of this study are available from the
	corresponding author upon reasonable request. 
	
	\section*{Acknowledgments}
	
	The  author gratefully acknowledges support from the São Paulo Research Foundation (FAPESP, Brazil), Process No. 2022/15371-3, Process No. 2024/00923-6, and 2026/01122-2.  \rev{The author also acknowledges} the NASA FIRMS program for providing VIIRS and MODIS active fire data, and the European Centre for Medium-Range Weather Forecasts (ECMWF) for providing ERA5 reanalysis data.

	\paragraph{Declarations. }
	\paragraph{Conflicts of Interest:} \rev{The author declares no conflicts of interest related to this study.}
	\paragraph{Ethics approval and consent to participate:}  Not applicable. 
	\paragraph{Consent for publication:} \rev{The author approves the final manuscript for publication.} 
	\paragraph{Clinical trial number:} Not applicable.
	
		\paragraph{Funding.}   Not applicable.
%
%
%
%

	\paragraph{AI Contribution}
	ChatGPT and Grammarly were used to assist with paraphrasing, language refinement, and improving the clarity and readability of certain sections of the manuscript, including the Abstract, Introduction, and Conclusion. The author reviewed and approved all AI-assisted revisions and take full responsibility for the final content of the manuscript.
	\bibliographystyle{ieeetr}
	\bibliography{sample,library,refs}
	
\end{document}